%% file: article.tex
\documentclass{article}

\input{shared}

\title{\Large Taylor approximation for chance constrained optimization problems governed by partial differential equations with high-dimensional random parameters \thanks{This research was partially funded by the National Science Foundation, Division of Mathematical Sciences under award DMS-2012453; the Department of Energy, Office of Science, Office of Advanced Scientific Computing Research, Mathematical Multifaceted Integrated Capability Centers (MMICCS) program under award DE-SC0019303; and the Simons Foundation under award 560651.}}

\usepackage{authblk}
\author[1]{Peng Chen}
\author[1,2]{Omar Ghattas}
\affil[1]{Oden Institute for Computational Engineering and Sciences, The University of Texas at Austin, Austin, TX 78712. (\url{peng@oden.utexas.edu}) }
\affil[2]{Department of Mechanical Engineering, and Department of Geological Sciences, The University of Texas at Austin, Austin, TX 78712. (\url{omar@oden.utexas.edu}) }

\date{}
\begin{document}
	
	\maketitle
	
	% REQUIRED
	\begin{abstract}
	 We propose a fast and scalable optimization method to solve chance or probabilistic constrained optimization problems governed by partial differential equations (PDEs) with high-dimensional random parameters. To address the critical computational challenges of expensive PDE solution and high-dimensional uncertainty, we construct surrogates of the constraint function by Taylor approximation, which relies on efficient computation of the derivatives, low rank approximation of the Hessian, and a randomized algorithm for eigenvalue decomposition. To tackle the difficulty of the non-differentiability of the inequality chance constraint, we use a smooth approximation of the discontinuous indicator function involved in the chance constraint, and apply a penalty method to transform the inequality constrained optimization problem to an unconstrained one. Moreover, 
	we design a gradient-based optimization scheme that gradually increases smoothing and penalty parameters to achieve convergence, for which we present an efficient computation of the gradient of the approximate cost functional by the Taylor approximation. Based on numerical experiments for a problem in optimal groundwater management, we demonstrate the accuracy of the Taylor approximation, its ability to greatly accelerate constraint evaluations, the convergence of the continuation optimization scheme, and the scalability of the proposed method in terms of the number of PDE solves with increasing random parameter dimension from one thousand to hundreds of thousands.
	\end{abstract}
	
%	% REQUIRED
%	\begin{keywords}
%		chance constrained optimization, Taylor approximation, randomized algorithm, smooth approximation, penalty method, Hessian
%	\end{keywords}
	
%	% REQUIRED
%	\begin{AMS}
%		65C20, 65D32, 65N12, 49J20, 93E20
%	\end{AMS}
	
	\section{Introduction}
	\label{sec:introduction}
	
	%\textbf{things to do}
	%
	%\begin{itemize}
	%    \item implement and check the gradient of quadratic approximation
	%    \item scalability result w.r.t.\ parameter dimension
	%\end{itemize}
	%
	%\textbf{what to show:}
	%
	%\begin{itemize}
	%    \item a good (meaningful, practical) benchmark example 
	%    \item optimization strategy 
	%    \item accuracy of Taylor approximation compared to SAA in approximating the chance
	%    \item accuracy of the optimization by Taylor approximation v.s.\ SAA
	%    \item approximation errors committed by smoothing and penalty function 
	%    \item scalability w.r.t.\ parameter dimension 
	%    \item Monte Carlo correction, variance reduction? 
	%\end{itemize}
	
	Large-scale simulation in computational science and engineering is
	often carried out not only to obtain insight about a system, but also
	as a basis for decision-making. When the decision variables represent
	the design or control or data-driven inference of model
	parameters of an engineered or natural system, and the system
	is governed by partial differential equations (PDEs), the task of
	determining the optimal design, optimal control, or inversion
	parameters leads to a {PDE-constrained optimization problem}. Over
	the past several decades, research in the field of PDE-constrained
	optimization has exploded, and powerful theory and algorithms are now
	available in the case of optimization governed by {deterministic}
	PDEs (e.g., see the monographs \cite{Lions71,
		Gunzburger03, HinzePinnauUlbrichEtAl08, Troeltzsch10}). 
	However, many PDE models are characterized by {random
		parameters} 
	due to lack of knowledge or intrinsic variability. These
	include initial or boundary conditions, sources, coefficients, and
	geometry. In many of these cases, the uncertainty arises from an
	({infinite-dimensional}) random field, leading to {high-dimensional}
	random parameters after discretization. It is critical to
	incorporate this uncertainty in the optimization problem to make the
	optimal solution more reliable and robust. Optimization under uncertainty has become an important research area and received increasing attentions in recent years \cite{ShapiroDentchevaRuszczynski09, BorziSchulzSchillingsEtAl10, HouLeeManouzi11,GunzburgerLeeLee11, RosseelWells12, TieslerKirbyXiuEtAl12, KouriHeinkenschloosVanBloemenWaanders12,ChenQuarteroniRozza13, Uryasev13, ChenQuarteroni14, KunothSchwab13, NgWillcox14, ChenQuarteroniRozza16, KunothSchwab16, KouriSurowiec16, BennerOnwuntaStoll16, AlexanderianPetraStadlerEtAl17, AliUllmannHinze17, XuBoyceZhangEtAl17, YangGunzburger17, GuoXuZhang17, RoaldAndersson17, LassUlbrich17, KouriSurowiec18, LiWangZhang18, Farshbaf-ShakerHenrionHoemberg18, MaLiJiang18, KolvenbachLassUlbrich18,  AckooijHenrionPerez-Aros19, VanAckooijMalick19, AckooijPerez-Aros19,  ChenVillaGhattas19, ChenGhattas19, ZahrCarlbergKouri19, VanBarelVandewalle19, GarreisSurowiecUlbrichEtAl19, LiStadler19, MartinNobileTsilifis19, AllaHinzeKolvenbachEtAl19, GeiersbachWollner19, GeletuHoffmannSchmidtEtAl20, DeHamptonMauteEtAl19, GeiersbachScarinci20, GeiersbachLoayza-RomeroWelker20}. 
	To account for the
	uncertainty in the optimization problem, different statistical
	measures of the objective function have been studied, e.g., mean,
	variance, conditional value-at-risk, worst case scenario,
	etc., \cite{RosseelWells12, KouriSurowiec16, YangGunzburger17, AlexanderianPetraStadlerEtAl17, LassUlbrich17, KouriSurowiec18, GarreisSurowiecUlbrichEtAl19}. Moreover, the treatment of chance constraints, also known as probabilistic constraints, i.e., the probability that a certain function exceeds a threshold or below a certain level, has also been investigated \cite{ShapiroDentchevaRuszczynski09, Uryasev13, XuBoyceZhangEtAl17, GuoXuZhang17, RoaldAndersson17, Farshbaf-ShakerHenrionHoemberg18, GeletuHoffmannSchmidtEtAl20}. Several computational challenges arise in solving optimization problems under uncertainty, especially with high-dimensional random parameters and inequality chance constraints.
	
	The {first prominent challenge} is that high-fidelity
	discretizations of the (nonlinear) PDEs often lead to large-scale
	(nonlinear) algebraic systems that are extremely expensive to solve in
	practical applications. Therefore, {only a limited number of
		high-fidelity PDE solves can be afforded}. This challenge prevents
	direct application of most of the conventional numerical methods
	for computing statistics of the objective function, since they require a large number of evaluations of the
	objective function and thus the PDE solution. To tackle this
	challenge, multigrid, multilevel, and model reduction methods have
	been successfully applied to solve  stochastic
	PDE-constrained optimization problems \cite{Borzi10, NgWillcox14, ChenQuarteroni14, ChenQuarteroniRozza16, AliUllmannHinze17, LassUlbrich17, 
		YangGunzburger17, ZahrCarlbergKouri19, AllaHinzeKolvenbachEtAl19}. The multigrid discretization and multilevel
	statistical evaluation rely on a hierarchical discretization of the
	PDE model and an efficient algorithm to balance the discretization
	error and the number of samples required for statistical evaluation at
	each level. However, due to the nature of the problem, it is not always
	% straightforward or even
	possible to use multigrid discretizations or gain computational
	savings by multilevel sampling because a sufficiently fine grid must
	be used to solve the PDE model 
	% such as the turbulent jet flow,
	with reasonable accuracy (as is often the case with hyperbolic or
	multiscale problems). Meanwhile, model reduction techniques become
	problematic for highly nonlinear problems that require effective
	affine approximation or when the solution manifold becomes high-dimensional, even
	if the objective function lives in a low-dimensional manifold.
	
	The second key {challenge} arises in computing the
	statistical measures, which involves {integration of the objective
		function with respect to (w.r.t.) the probability measure of the
		high-dimensional random parameters.} A classical approach known as
	sample average approximation (SAA) or Monte Carlo quadrature is to
	take the average of the objective function at a set of samples
	randomly drawn from the probability measure of the
	random parameters. However, its convergence rate is only $O(M^{-1/2})$, where
	an often expensive PDE has to be solved {for each of the $M$
		samples}.  The
	resulting deterministic optimization problem has $M$ PDE constraints that need to be solved to evaluate the objective function,
	and as such is typically {prohibitive to solve.} In recent years,
	stochastic Galerkin and stochastic collocation based integration
	methods have been used to compute the statistical moments (e.g., mean
	and variance) of the objective function in stochastic optimization
	\cite{Borzi10, HouLeeManouzi11, GunzburgerLeeLee11, 
		TieslerKirbyXiuEtAl12,RosseelWells12,KouriHeinkenschloosVanBloemenWaanders12,ChenQuarteroniRozza13,
		ChenQuarteroniRozza13, ChenQuarteroni14,
		KunothSchwab13, ChenQuarteroniRozza16,
		KunothSchwab16}, provided that a suitable finite dimensional
	parametrization of the random parameters, such as a truncated
	Karhunen--Lo\`eve expansion, is available. These methods achieve fast
	convergence when the objective function depends smoothly on the
	low-dimensional parameters, but suffer from the so-called
	{curse of dimensionality}, i.e.\ the convergence rate quickly
	deteriorates as the  parameter dimension increases. More recent
	advances in adaptive and anisotropic sparse quadrature
	\cite{SchillingsSchwab13, Chen18} and high order quasi Monte Carlo
	methods \cite{DickLeGiaSchwab16} have been shown to achieve a
	convergence rate of $O(M^{-s})$, with $s$ potentially much larger than
	$1/2$ and independent of the nominal dimension of the random
	parameters, thus mitigating the curse of dimensionality.  
	%The convergence rate depends on the parametrization and the regularity of the objective function w.r.t.\ the parameter. 
	However, if 
	%the parametrization does not correctly capture the effective parameter dimensions, or 
	the objective function is not sufficiently smooth or sufficiently
	anisotropic in the parameter space, the convergence of these methods becomes very slow, or worse than that of Monte Carlo.

	The {third critical challenge} comes from the
	{non-differentiability of the chance constraint}, which involve
	integration of a discontinuous indicator function. The discontinuity makes the cost functional
	non-differentiable w.r.t.\ the optimization variable, so
	that rapidly convergent derivative based optimization methods, e.g., steepest descent or Newton methods, cannot be directly applied. Therefore, solving optimization problems with such non-differentiable constraints
	not only requires a large number of PDE solves at each optimization
	iteration, but also requires a large number of optimization
	iterations, especially when the optimization variable dimension 
	is high. To address this challenge, proper smoothing techniques have
	been employed to approximate the indicator function by differentiable
	functions \cite{ChenMangasarian95, QiSunZhou00, Uryasev13,
		KouriSurowiec16}, which introduces smoothing errors in the cost
	functional and optimal solution. Alternatively, a dual approach has
	been developed that reformulates the risk averse stochastic
	optimization problem as a probability constrained minimax problem
	\cite{RuszczynskiShapiro06, ShapiroDentchevaRuszczynski09, Uryasev13,
		KouriSurowiec16}. However considerable difficulties are
	still encountered in finding numerical approximations that satisfy the
	probability constraints in the presence of high-dimensional random
	parameters.

	\subsection*{Contributions}
	In this work we address the above computational challenges by proposing a Taylor approximation based continuation optimization method to solve chance constrained optimization problems governed by PDEs with high-dimensional random parameters. Following our recent work \cite{AlexanderianPetraStadlerEtAl17, ChenVillaGhattas19}, 
	we extend the Taylor approximation---including constant, linear, and quadratic approximations---of the objective function used to evaluate its mean, to approximation of the constraint function used to accelerate the evaluation of the chance/probability of the constraint function. In particular, a double-pass randomized algorithm is employed to solve a generalized eigenvalue problem, where the eigenvalues and eigenfunctions are used to construct the quadratic term of the Taylor approximation. A set of linearized PDEs are derived for the computation of the gradient and action of the Hessian of the objective and constraint functions w.r.t.\ the random parameters. To solve the optimization problem, we present a continuation BFGS algorithm, which features (1) smooth approximation of the indicator function, (2) a penalty method to transform the inequality constrained optimization problem to an unconstrained one, (3) a continuation scheme with an outer loop of increasing the smoothing and penalty parameters and an inner loop of BFGS optimization. The computation of the approximate cost functional and its gradient w.r.t.\ the optimization variable are presented in detail for both SAA and Taylor approximation. For the proposed method, we demonstrate (1) the accuracy of the Taylor approximations, (2) the efficiency on the surrogate acceleration, (3) the convergence of the continuation optimization algorithm, and (4) the independence of the number of PDE solves from increasing random parameter dimension.
	The demonstrations are carried out by numerical experiments for an example of water management in agricultural irrigation, where the PDE model is a Darcy flow equation that describes groundwater flow in the presence of an uncertain permeability field. The optimization objective is to extract water at given well locations that meets a target extraction rate, while a chance constraint is imposed on an integrated pressure field to prevent a pressure that is low enough to lead to collapse or damage of the aquifer. 
	
	\subsection*{Notations}
	% In this section we present chance constrained optimization under uncertainty in an abstract setting. To this aim, let us first introduce the notation used throughout the paper. 
	Let $\cX$ be a {Banach} space and $\cX'$ the dual space; $_\cX \langle \cdot, \cdot \rangle_{\cX'}$ then denotes the duality pairing between the spaces $\cX$ and $\cX'$. For ease of notation, we will omit specification of the subscritps $\cX$ and $\cX'$ and simply write $\langle \cdot, \cdot \rangle$ when the spaces can be inferred from the context without ambiguity. Given two {Banach} spaces $\cX$ and $\cY$ and a map $f: \cX \times \cY \mapsto \mathbb{R}$, $\partial_x f(x, y) \in \cX'$ denotes the Fr\'echet derivative of $f(x,y)$ with respect to $x$ evaluated at $(x,y)$, which satisfies
	\begin{equation}
	\lim_{\tilde{x} \to 0} \frac{f(x+\tilde{x}, y) - f(x, y) - \langle \tilde{x},  \partial_x f(x, y) \rangle}{||\tilde{x}||_\cX}  = 0.
	\end{equation}
	
	Let $\partial_{xy}f(x,y): \cY \mapsto \cX'$ denote the Fr\'echet derivative of $\partial_x f(x,y)$ with respect to $y$ evaluated at $(x,y)$, or the second order (mixed) Fr\'echet derivative of $f(x,y)$ with respect to $x$ and $y$ evaluated at $(x,y)$. Similarly, $\partial_{yx} f(x,y): \cX \mapsto \cY'$ denotes the Fr\'echet derivative of $\partial_y f(x,y)$ with respect to $x$ evaluated at $(x,y)$, which is the adjoint operator of $\partial_{xy}f(x,y)$ and satisfies 
	\begin{equation}\label{eq:adjoint_op}
	\langle \tilde{x}, \partial_{xy} f \, \hat{y} \rangle: = \;_\cX \langle \tilde{x}, \partial_{xy} f \, \hat{y} \rangle_{\cX'} = \;_\cY \langle \hat{y}, \partial_{yx} f \, \tilde{x} \rangle_{\cY'} = :\langle \hat{y}, \partial_{yx} f \, \tilde{x} \rangle, \quad \forall \tilde{x} \in \cX, \, \hat{y} \in \cY,
	\end{equation}
	where we have omitted the argument $(x,y)$ for simplicity.
	
	The rest of the paper is organized as follows. In Section \ref{sec:chance} we present the general formulation of PDE and chance constrained optimization problems, which is followed by Section \ref{sec:taylor} on SAA, Taylor approximation, a randomized algorithm for low rank approximation, and computation of the gradient and Hessian action of the objective and constraint function w.r.t.\ the random parameters. Section \ref{sec:gradientopt} is devoted to the presentation of a continuation gradient-based optimization method that involves smooth approximation of the indicator function, a penalty method for the inequality constraint, a continuation scheme to increase the smoothing and penalty parameters, and the computation of the gradient of the approximate cost functional w.r.t.\ the optimization variable. Numerical experiments and results are reported in Section \ref{sec:numerics} for the demonstration of the accuracy, efficiency, convergence, and scalability of the proposed method. Conclusions and perspectives are drawn in Section \ref{sec:conclusion}.
	
	\section{Chance constrained optimization} 
	\label{sec:chance}
	
	We consider a system to be optimized under uncertainty, which is modeled by partial differential equations presented in an abstract strong (residual) form as: find $u\in \cU$, such that
	\begin{equation}\label{eq:PDE_strong}
	\cR(u, m, z) = 0 \quad \text{ in }  \cV',
	\end{equation}
	where $m \in \cM$ is an uncertain or random parameter field that lives in a separable Banach space $\cM$, which has a probability distribution $\mu$;
	$z \in \cZ$ is an optimization variable in a separable Banach space $\cZ$; and $\cR(\cdot, m, z) : \cU  \mapsto \cV'$ denotes a (possibly nonlinear) operator from $\cU$ to $\cV'$, the dual of $\cV$, where $\cU$ and $\cV$ are two separable Banach spaces. The weak form of \eqref{eq:PDE_strong} is given by means of duality paring as: find $u\in \cU$, such that
	\begin{equation}\label{eq:PDE} 
	r(u, v, m, z) := \; _{\cV} \langle v, \cR(u, m, z) \rangle_{\cV'} = 0 \quad \forall v \in \cV,
	\end{equation} 
	where $v$ is a test variable or an adjoint variable in the optimization context. 
	By definition, $r(u, v, m, z)$ is linear with respect to (w.r.t.) $v$, and may be nonlinear w.r.t.\ $u, m$, and $z$. 
	
	By $q: \cU \times \cM \times \cZ \mapsto \bR$ and $f: \cU \times \cM \times \cZ \mapsto \bR$ we denote an objective function and a constraint function for the optimization as real-valued, continuous, and possibly nonlinear maps of $u, m, z$. Since $u$ depends on $m$ and $z$ through \eqref{eq:PDE}, for simplicity we write $q(m,z) = q(u(m,z),m,z)$ and $f(m,z) = f(u(m,z),m,z)$ by slight abuse of notation. 
	%For simplicity, we assume that $q$ and $f$ depend on $m$ and $z$ only implicitly through the solution of the state equation \eqref{eq:PDE}. The general case where $q$ and $f$ have an explicit dependence on $m$ and $z$ can be treated with slight modifications in the derivation of the cost functional and its gradient, as we address along the presentation in Sec.\ \ref{sec:taylor} and Sec.\ \ref{sec:gradientopt}. 
	
	For the optimization problem, we consider a cost functional 
	\begin{equation}
	\cJ(z) = \bE[q(\cdot, z)] + \cP(z)
	\end{equation}
	where the first term is the mean of the objective function $q$ defined as
	\begin{equation}\label{eq:mean}
	\bE[q(\cdot, z)] = \int_{\cM} q(m, z) d\mu(m),
	\end{equation}
	the second term $\cP(z)$ represents a penalty or regularization term for the optimization variable $z$. 
	Moreover, 
	we consider a chance constraint
	\begin{equation}\label{eq:probability_constraint}
	P(f(\cdot, z) \geq 0) \leq \alpha   
	\end{equation}
	for a critical chance $0< \alpha < 1$, where the probability is given by 
	\begin{equation}\label{eq:probability_definition}
	P(f(\cdot, z) \geq 0) = \bE[\bI_{[0, \infty)}(f(\cdot, z))] = \int_\cM \bI_{[0, \infty)}(f(m, z)) d\mu(m),
	\end{equation}
	where $\bI_{[0, \infty)}(f(m, z))$ is an indicator function defined as 
	\begin{equation}
	\bI_{[0, \infty)}(f(m, z)) = 
	\left\{
	\begin{array}{cc}
	1     &  \text{ if } f(m, z) \geq 0, \\
	0     &  \text{ if } f(m, z) < 0.
	\end{array}
	\right.
	\end{equation}
	Then the PDE and chance constrained optimization problem can be formulated as 
	\begin{equation}\label{eq:optimization}
	\min_{z \in \cZ}  \cJ(z),  \text{ subject to the PDE constraint } \eqref{eq:PDE} \text{ and the chance constraint } \eqref{eq:probability_constraint}.
	\end{equation}

	%%%%%%%%%%%%%%%%%%%%%%%%%%%%%%%%%%%%%%%%%%%%%%%%%
	%%%%%%%%%%%%%%%%%%%%%%%%%%%%%%%%%%%%%%%%%%%%%%%%%
	\section{Taylor approximation}
	\label{sec:taylor}
	%%%%%%%%%%%%%%%%%%%%%%%%%%%%%%%%%%%%%%%%%%%%%%%%%
	
	We first present a sample average approximation (SAA) for the optimization problem \eqref{eq:optimization}. Then we introduce an (up to quadratic) Taylor approximation for both the objective function and the constraint function, which requires an efficient eigenvalue decomposition of the Hessian of the objective and constraint functions w.r.t. the random parameter field. We present a double-pass randomized algorithm for this task, which requires only actions of the Hessian in random directions without direct assess to the entries of the Hessian matrix.
	
	% We develop the main approximation methods in this section, which include the approximation by (low-order) Taylor expansion of the control objective, the computation of the gradient and Hessian of the control objective with respect to the uncertain parameter used in the Taylor expansion, the randomized algorithm to solve a generalized eigenvalue problem for trace estimation, and the Monte Carlo correction for the remainder of the Taylor expansion.
	
	\subsection{Sample average approximation}
	\label{sec:MCinteg}
	The mean of the objective function can be evaluated by the sample average approximation (SAA) 
	\begin{equation}
	\bE[q(\cdot, z)] \approx q_M(z) :=  \frac{1}{M_q} \sum_{i=1}^{M_q} q(m_i, z),
	\end{equation}
	where $m_i$, $i = 1, \dots, M_q$, are independent identically distributed (i.i.d.) random samples drawn from the probability distribution $\mu$. Similarly, the chance constraint \eqref{eq:probability_definition} can be approximated by 
	\begin{equation}\label{eq:ChanceSAA}
	P(f(\cdot, z) \geq 0) \approx f_M(z) := \frac{1}{M_f} \sum_{i=1}^{M_f} \bI_{[0, \infty)}(f(m_i, z)),
	\end{equation}
	where $m_i$, $i = 1, \dots, M_f$, are i.i.d.\ random samples drawn from $\mu$. Note that to obtain an accurate approximation $q_M$ and especially $f_M$ for $\alpha$ close to $0$, a large number of samples are required due to the slow convergence (with rate $O(M^{-1/2})$) of the SAA approximation, thus making this approach computationally prohibitive if the PDE solve at one sample is expensive. 
	%since the number of PDE solves at each optimization step is very large.

	\subsection{Taylor approximation}
	%Given the control $z$, the state variable $u$ is determined by the realization of the uncertain parameter $m$. Therefore, by a slight abuse of notation we can denote the control objective as $q(m)$, depending only on $m$ given $z$.
	% In this section we provide expressions for the Taylor expansion of the control objective $q(m,z)$ with respect to the uncertain parameter $m$ for a fixed value of the control $z$. 
	%In what follows, we denote with the subscript $m$ the total derivative of $q(u(m,z))$ with respect to $m$ and we assume sufficient smoothness so that all Fr\'echet derivatives are well-defined. 
	% For ease of notation, we denote with $\langle \cdot, \cdot \rangle$ the duality-pairing $_\cM \langle \cdot, \cdot \rangle_{\cM'}$.
	%and, with a slight abuse of notation, we write $q(m)$ to indicate the control objective $q(u(m,z))$ evaluated at a fixed value of $z$.
	
	We assume that the objective function $q$ admits the $k$-th order Fr\'echet derivative with respect to $m$ at $\bar{m} \in \cM$, denoted as $\nabla_m^k q (\bar{m}, z)$, for $k = 1, \dots, K$. A $K$-th order Taylor expansion of the objective function $q$ evaluated at $\bar{m} \in \cM$ is given by
	\begin{equation}\label{eq:Taylor}
	T_K q(m, z) = \sum_{k = 0}^K \frac{1}{k!}\nabla_m^k q (\bar{m}, z) (m-\bar{m})^k.
	\end{equation}
	Note that $\nabla_m^k q (\bar{m}, z): \cM^k \mapsto \bR$ is a multilinear map defined in the tensor-product space 
	\begin{equation}
	\cM^k = \prod_{i = 1}^k \cM_i, \text{ where } \cM_i = \cM.
	\end{equation}
	For the Gaussian measure $m \sim \cN(\bar{m}, \cC)$, we have the analytic expression \cite{AlexanderianPetraStadlerEtAl17} %\cite{AlexanderianPetraStadlerEtAl17}
	\begin{equation}\label{eq:mean_objective}
	\bE[q(\cdot, z)] \approx \bE[T_{K}q(\cdot, z)] =
	\left\{
	\begin{array}{cc}
	q(\bar{m}, z),     &   K = 0, 1, \\[4pt]
	\displaystyle
	q(\bar{m}, z)  + \frac{1}{2} \text{trace}\left(\cC^{1/2} \nabla_m^2 q(\bar{m}, z) \cC^{1/2}\right), &  K = 2, 3.
	\end{array}
	\right.
	\end{equation}
	Let $\text{trace}(\cH_q)$ represent the trace of $\cH_q = \cC^{1/2} \nabla_m^2 q(\bar{m}, z) \cC^{1/2}$, which is the covariance-preconditioned Hessian of the objective function $q$. 
	
	Similar to \eqref{eq:Taylor}, under the assumption that the constraint function $f$ admits the $k$-th order Fr\'echet derivative with respect to $m$ at $\bar{m} \in \cM$ for $k = 1, \dots, K$, we construct a $K$-th order Taylor expansion of the constraint function $f$ at $\bar{m}$, denoted $T_K f$. 
	% One choice of $m_f$ is $\bar{m}$. To construct a more accurate Taylor approximation on the limit surface $S = \{m \in \cM: f(m, z) = 0\}$, which is important to evaluate the probability $P(f(m, Z) \leq 0)$, we can choose $m_f$ such that 
	% \begin{equation}
	% f(m_f, z) = 0.
	% \end{equation}
	Then we can approximate the probability $P(f(m, z) \geq 0)$ by SAA \eqref{eq:ChanceSAA} with the Taylor approximation $T_K f$ as 
	\begin{equation}\label{eq:ChanceTaylor}
	P(f(m, z) \geq 0) \approx f_M^K(z) : = \frac{1}{M_f} \sum_{i = 1}^{M_f} \bI_{[0, \infty)}(T_K f(m_i, z)).
	\end{equation}
	The attractiveness of the Taylor approximations of the objective and constraint, \eqref{eq:mean_objective} and \eqref{eq:ChanceTaylor}, is that once they are constructed, no further PDE solves are required. In the next section, we shall see how these Taylor approximations can b efficiently constructed.
	
	\subsection{Low-rank approximation}
	
	To compute the trace in \eqref{eq:mean_objective}, we can use
	\begin{equation}
	\text{trace}(\cH_q) = \sum_{n = 1 }^\infty \lambda_n^q,
	\end{equation}
	where $(\lambda_n^q )_{n\geq 1}$ are the eigenvalues of $\cH_q = \cC^{1/2} \nabla_m^2 q(\bar{m}, z) \cC^{1/2}$, which are equivalent to the generalized eigenvalues of $(\nabla_m^2 q(\bar{m}, z), \cC^{-1})$.  In practice, the (absolute) eigenvalues decay rapidly,  $|\lambda_n^q| \to 0$ as $n \to \infty$, as proven for some model problems and demonstrated numerically for many others in forward uncertainty quantification, Bayesian inversion and experimental design, and stochastic optimization \cite{BashirWillcoxGhattasEtAl08, FlathWilcoxAkcelikEtAl11,
		Bui-ThanhGhattas12a, Bui-ThanhGhattas13a, Bui-ThanhGhattas12,
		Bui-ThanhBursteddeGhattasEtAl12,
		Bui-ThanhGhattasMartinEtAl13,
		AlexanderianPetraStadlerEtAl16, AlexanderianPetraStadlerEtAl17,
		AlexanderianPetraStadlerEtAl14, CrestelAlexanderianStadlerEtAl17,
		PetraMartinStadlerEtAl14, IsaacPetraStadlerEtAl15,
		MartinWilcoxBursteddeEtAl12, Bui-ThanhGhattas15, ChenVillaGhattas17, ChenVillaGhattas19, ChenGhattas19a, ChenWuChenEtAl19a, ChenGhattas20, ChenHabermanGhattas20, ChenWuGhattas20, WuChenGhattas20}. 
	Given rapid decay of the eigenvalues, we can approximate the trace by the $N_q$ largest (in absolute value) eigenvalues $\lambda_n^q, \; n = 1, \dots, N_q$, i.e., 
	\begin{equation}\label{eq:TraceApprox}
	\text{trace}(\cH_q) \approx \sum_{n = 1 }^{N_q} \lambda_n^q.
	\end{equation}
	%\pc{Add further development/result on convergence analysis.}
	
	For computational efficiency, we consider the generalized eigenvalues of $(\nabla_m^2 q(\bar{m}, z), \cC^{-1})$ by solving the generalized eigenvalue problem: find $(\lambda_n^q, \psi_n^q)$, $n = 1, \dots, N_q$, such that  
	\begin{equation}\label{eq:GenEigen_q}
	\nabla_m^2 q(\bar{m}, z) \psi_n^q = \lambda_n^q \cC^{-1} \psi_n^q, \quad n = 1, \dots, N_q,
	\end{equation}
	with $|\lambda_1^q| \geq \cdots \geq |\lambda_{N_q}^q|$ corresponding to the $N_q$ largest eigenvalues in absolute value and the eigenfunctions satisfying the orthonormality condition with respect to $\cC^{-1}$, i.e.,
	\begin{equation}\label{eq:Orthon_q}
	\langle\psi_n^q, \cC^{-1}\psi_{n'}^q \rangle = \delta_{nn'}, \quad n, n' = q, \dots, N_q.
	\end{equation}
	
	To compute the Taylor approximation $T_K f(m, z)$ at $\bar{m}$ for $K = 2$, we need to evaluate $\nabla_m^2 f(\bar{m}, z)(m-\bar{m})^2$. If the eigenvalues of $\cH_f = \cC^{1/2}\nabla_m^2 f(\bar{m}, z) \cC^{1/2}$ decay rapidly to $0$, we can we approximate $\nabla_m^2 f(\bar{m}, z)(m-\bar{m})^2$ by a low-rank approximation as
	% \begin{equation}\label{eq:HessianApprox}
	% \nabla_m^2 f(\bar{m}) \approx  \cC^{-1} \left(\sum_{n = 1}^{N_f} \lambda_n^f \psi_n^f \otimes \psi_n^f \right) \cC^{-1},
	% \end{equation}
	\begin{equation}\label{eq:HessianApprox}
	\nabla_m^2 f(\bar{m}, z) (m-\bar{m})^2 \approx   \sum_{n = 1}^{N_f} \lambda_n^f \langle m-\bar{m},   \cC^{-1} \psi_n^f \rangle^2,
	\end{equation}
	where $(\lambda_n^f, \psi_n^f)$, $n = 1, \dots, N_f$, are the solution of the generalized eigenvalue problem 
	\begin{equation}\label{eq:GenEigen_f}
	\nabla_m^2 f(\bar{m},z) \psi_n^f = \lambda_n^f \cC^{-1} \psi_n^f, \quad n = 1, \dots, N_f,
	\end{equation}
	with  $|\lambda_1^f| \geq \cdots \geq |\lambda_{N_f}^f|$ corresponding to the $N_f$ largest eigenvalues in absolute value and the eigenfunctions satisfying the orthonormality condition with respect to $\cC^{-1}$, i.e., 
	\begin{equation}\label{eq:Orthon_f}
	\langle\psi_n^f, \cC^{-1}\psi_{n'}^f \rangle = \delta_{nn'}, \quad n, n' = 1, \dots, N_f.
	\end{equation}
	With the low-rank (LR) decomposition of the Hessian $\nabla_m^2 f(\bar{m})$ in \eqref{eq:HessianApprox}, we define the quadratic Taylor approximation of $f$ corresponding to \eqref{eq:Taylor} as 
	\begin{equation}\label{eq:TaylorLR}
	T_2^{\text{LR}} f(m,z) := f(\bar{m},z) + \nabla_m f(\bar{m},z) (m - \bar{m}) + \frac{1}{2}  \sum_{n = 1}^{N_f} \lambda_n^f \langle m-\bar{m},  \cC^{-1}\psi_n^f\rangle^2.
	\end{equation}
	
	%\subsection{Variance reduction for $\bE[q]$}
	%
	%\subsection{Multifidelity approximation by a-posteriori error indicator}
	%
	%We observe that to obtain a faithful evaluation of the indicator function $\bI_{[0, \infty)}(f(m_i,z))$, we can use $\bI_{[0, \infty)}(T_K f(m_i,z))$ as long as the Taylor approximation satisfies 
	%\begin{equation}
	%|f(m_i,z) - T_K f(m_i,z)| \leq |T_K f(m_i,z)|.
	%\end{equation}
	%
	%In fact, if $T_K f(m_i,z) < 0$ and $|f(m_i,z) - T_K f(m_i,z)| \leq - T_K f(m_i,z) $, then we have $f(m_i,z) < 0$, so that $\bI_{[0, \infty)}(f(m_i,z)) = \bI_{[0, \infty)}(T_K f(m_i,z)) = 0$. On the other hand, if $T_K f(m_i,z) \geq 0$ and $|f(m_i,z) - T_K f(m_i,z)| \leq T_K f(m_i,z)$, then we have $f(m_i,z) \geq 0$, so that $\bI_{[0, \infty)}(f(m_i,z)) = \bI_{[0, \infty)}(T_K f(m_i,z)) = 1$. \pc{to be further developed.}
	
	% \subsection{Randomized SVD}
	
	%\begin{algorithm}
	%	\caption{Randomized SVD for generalized eigenvalue decomposition.}
	%	\label{alg:RandomizedSVD}
	%	\begin{algorithmic}[1]
	%		\STATE{to be added }
	%	\end{algorithmic}
	%\end{algorithm}
	
	% Let $H$ and $C$ denote the discrete version of the Hessian $\cH$ and covariance $\cC$
	
	To solve the generalized eigenvalue problems \eqref{eq:GenEigen_q} and \eqref{eq:GenEigen_f}
	for the dominant eigenvalues, we apply a double-pass randomized algorithm \cite{HalkoMartinssonTropp11, SaibabaLeeKitanidis16}, presented in Algorithm \ref{alg:RandomizedSVD}.  Here, by $H$ and $C^{-1}$ of size $N_h \times N_h$ each, we denote discrete approximations of the Hessians $\nabla_m^2 q$ and $\nabla_m^2 f$ and the covariance $\cC^{-1}$, e.g., by a finite element method. In Algorithm \ref{alg:RandomizedSVD} only the action of $H$ and $C$ on a given vector is required, which does not require access to the entries of $H$ and $C$.
	
	\begin{algorithm}[!htb]
		\caption{Double-pass randomized eigensolver for $(H, C^{-1})$}
		\label{alg:RandomizedSVD}
		\begin{algorithmic}
			\STATE{\textbf{Input: } the number of desired eigenpairs $N$, an oversampling factor $c \leq 10$.}
			\STATE{\textbf{Output: } $(\Lambda_N, \Psi_N)$ with $\Lambda_N = \text{diag}(\lambda_1, \dots, \lambda_N)$ and $\Psi_N = (\psi_1, \dots, \psi_N)$.}
			\STATE{1. Draw a Gaussian random matrix $\Omega \in \bR^{N_h \times (N+c)}$.}
			\STATE{2. Compute $Y = C (H \Omega)$.}
			\STATE{3. Compute $QR$-factorization $Y = QR$ such that $Q^\top C^{-1} Q = I_{N+c}$.}
			\STATE{4. Form $T = Q^\top H Q \in \bR^{(N+c)\times (N+c)}$ and compute eigendecomposition $T = S \Lambda S^\top$.}
			\STATE{5. Extract $\Lambda_N = \Lambda(1:N, 1:N)$ and $\Psi_N = QS_N$ with $S_N = S(:,1:N)$.}
		\end{algorithmic}
	\end{algorithm}
	
	In the next section, we present the computation of the Hessian actions $H \Omega$ and $H Q$, which dominates the cost in Algorithm \ref{alg:RandomizedSVD}, by solving linearized PDEs. 
	%that involve $4(k+c)$ linearized PDE solves, including $2(k+c)$ for the incremental state equations and $2(k+c)$ for the incremental adjoint equations. 
	%The error of the eigenvalues $\lambda_j$, $j = 1, \dots, k$, are bounded by the remaining ones $\lambda_j$, $j > k$, which is small if they decay fast \cite{avron2011randomized, saibaba2016randomized}. Moreover, the oversampling factor $c$ is typically small, $c \leq 10$, to achieve high computational accuracy with high probability \cite{avron2011randomized}. 
	As observed in \cite{ChenVillaGhattas19, ChenGhattas19a, ChenHabermanGhattas20} and in the numerical results in Section \ref{sec:numerics}, the advantages of Algorithm \ref{alg:RandomizedSVD} are: (i) the error in the approximation of the eigenvalues $\lambda_n$, $n = 1, \dots, N$, is bounded by the remaining ones $\lambda_n$, $n > N$, which is small if they decay rapidly;
	(ii) the computational cost is dominated by $2(N+c)$ Hessian actions, where the application of $C$ on a vector is inexpensive, e.g., it takes only $O(N_h)$ operations by a multigrid solver for $C$ discretized from a differential operator; (iii) it is scalable in terms of the number of PDEs to solve, because $N$ typically does not change when $N_h$ increases;
	%and $O((k+c)^2n)$ operations if needed for computing $\Psi_k$;
	and (iv) computing the Hessian actions $H\Omega$ and $H Q$ can be asynchronously parallelized.
	
	\subsection{Computation of the gradient and Hessian action}
	% of $q$ and $f$ w.r.t.\ the uncertain variable $m$ at its mean $\bar{m}$ 
	For a given optimization variable $z$, we first compute $u(\bar{m})$ by solving the state equation \eqref{eq:PDE} at $\bar{m}$, which can be equivalently written as: find $u \in \cU$ such that 
	\begin{equation}\label{eq:state}
	\langle \tilde{v}, \partial_v \bar{r} \rangle = 0 \quad \forall \tilde{v} \in \cV,
	\end{equation}
	where for ease of notation we use $\bar{r}$ to represent the weak form \eqref{eq:PDE} at $\bar{m}$, i.e.,
	\begin{equation}
	\bar{r} = r(u,v,\bar{m},z).
	\end{equation}
	Then we can evaluate the objective function $q(\bar{m},z) = q(u(\bar{m},z),\bar{m},z)$ and the constraint function $f(\bar{m},z) = f(u(\bar{m},z),\bar{m},z)$.
	The gradient and Hessian of $q$ is computed the same way as for $f$, so we present only the derivation for $q$ and then state the result for $f$.
	
	% To compute $q(\bar{m})$, we first solve the state problem \eqref{eq:PDE} at $m = \bar{m}$, and then compute $q(\bar{m}) = q(u(\bar{m}))$. To compute the gradient $\nabla_m q(\bar{m})$, we form a Lagrangian as 
	% \begin{equation}
	% \cL(u, v, m, z) = q(u) + r(u, v, m, z),
	% \end{equation}
	% where the adjoint variable $v$ plays the role of a Lagrangian multiplier. By taking the first order variation of the Lagrangian at $\bar{m}$ w.r.t.\ $v$ to zero, we obtain the state $u \in \cU$: find $v \in \cV$ such that 
	% \begin{equation}
	% \langle \tilde{u}, \partial_u \bar{r}(u, v, z) \rangle = - \langle \tilde{u}, \partial_u q \rangle, \quad \forall \tilde{u} \in \cU.
	% \end{equation}
	% Then the gradient of $q$ can be evaluated as
	% \begin{equation}
	% \langle \tilde{m}, \nabla_m q(\bar{m}) \rangle = \langle \tilde{m}, \partial_m \bar{r}(u, v, z) \rangle, \quad \forall \tilde{m} \in \cM.
	% \end{equation}
	
	% In this section we use the Lagrangian formalism to derive the expressions for the gradient $q_m$ and for the action of the Hessian $q_{mm}$ in given direction $\hat{m}$ evaluated at $\bar{m}$. 
	
	%In the trace approximation and the variance reduction, we need to compute the gradient and the Hessian of the control objective with respect to the uncertain parameter acting in given directions. 
	%We employ a Lagrangian multiplier method for their computations, with the Lagrangian defined as
	We use the Lagrangian formalism to derive the expressions for the gradient $\nabla_m q(\bar{m},z)$ evaluated at $\bar{m}$ and for the action of the Hessian $\nabla_m^2 q(\bar{m},z)$ evaluated at $\bar{m}$ in given direction $\hat{m}$. 
	First, we define the Lagrangian functional
	\begin{equation}\label{eq:Lagrangian Inner}
	\cL(u,v^q,m,z) := q(u,m,z) + r(u,v^q,m,z),
	\end{equation}
	where the adjoint $v^q$ is the Lagrange multiplier for the state equation in the computation w.r.t.\ the objective function $q$. 
	In what follows, for ease of notation, we define
	\begin{equation}\label{eq:rq_bar}
	\bar{r}_q = r(u,v^q,\bar{m},z), \text{ and } \bar{q} = q(u(\bar{m}, z),\bar{m},z).
	\end{equation}
	% By requiring the first order variation of \eqref{eq:Lagrangian Inner} at $\bar{m}$ with respect to the adjoint $v^q$ to vanish, we obtain the state problem: find $u \in \cU$, such that 
	% \begin{equation}\label{eq:state}
	% \langle \tilde{v}, \partial_v \bar{r}_q\rangle = 0, \quad \forall \tilde{v} \in \cV,
	% \end{equation}
	% Similarly, by 
	By setting the first order variation of \eqref{eq:Lagrangian Inner} at $\bar{m}$ with respect to the state $u$ to zero, 
	we obtain the adjoint problem: find $v^q \in \cV$, such that
	\begin{equation}\label{eq:adjoint_q}
	\langle \tilde{u} , \partial_u \bar{r}_q\rangle = - \langle \tilde{u}, \partial_u \bar{q} \rangle, \quad \forall \tilde{u} \in \cU.
	\end{equation} 
	Then, the gradient of $q$ at $\bar{m}$ acting in direction $\tilde{m}$ is given by  
	\begin{equation}\label{eq:gradient_q}
	\langle \tilde{m}, \nabla_m \bar{q} \rangle  =  \langle \tilde{m}, \partial_m \cL \rangle =\langle \tilde{m}, \pc{\partial_m \bar{q} +} \partial_m \bar{r}_q\rangle, \quad \forall \tilde{m} \in \cM,
	\end{equation}
	where $u$ solves the state problem \eqref{eq:state} and $v$ solves the adjoint problem \eqref{eq:adjoint_q}.
	%which can be computed by first solving the the state problem \eqref{eq:state} and then the adjoint problem \eqref{eq:adjoint_q}, followed by evaluating \eqref{eq:gradient_q} in direction $\tilde{m}$. Note that $\langle \tilde{m}, \partial_m \bar{q} \rangle$ should be added to \eqref{eq:gradient} if $q$ explicitly depends on $m$.
	Similarly, for the computation of the Hessian of $q$ acting in direction $\hat{m}^q$, we form the Lagrangian
	\begin{equation}\label{eq:meta_lag}
	\cL^H(u,v^q, \bar{m},z; \hat{u}^q, \hat{v}^q, \hat{m}^q) := \langle \hat{m}^q, \pc{\partial_m \bar{q} + }\partial_m \bar{r}_q\rangle + \langle \hat{v}^q, \partial_v \bar{r}_q\rangle + \langle \hat{u}^q , \partial_u \bar{r}_q + \partial_u \bar{q} \rangle,
	\end{equation}
	where $\hat{u}^q$ and $\hat{v}^q$ denote the incremental state and incremental adjoint, respectively.
	%More precisely, by variation of the state problem, 
	%to hold the following relations
	%\begin{equation}
	%\left(
	%\begin{array}{lll}
	%\partial_{vv} \cL &  \partial_{vu} \cL & \partial_{vm} \cL \\
	%\partial_{uv} \cL &  \partial_{uu} \cL & \partial_{um} \cL \\
	%\partial_{mv} \cL &  \partial_{mu} \cL & \partial_{mm} \cL 
	%\end{array}
	%\right)
	%\left(
	%\begin{array}{lll}
	%\hat{v} \\
	%\hat{u} \\ 
	%\hat{m}
	%\end{array}
	%\right)
	%= 
	%\left(
	%\begin{array}{ccc}
	%0 \\
	%0 \\
	%q_{mm}(\bar{m}) \hat{m}
	%\end{array}
	%\right)
	%\end{equation}
	%we obtain the incremental state problem: find $\hat{u} \in \cU$ such that
	By taking the variation of \eqref{eq:meta_lag} with respect to the adjoint $v^q$ and using \eqref{eq:adjoint_op}, we obtain the incremental state problem: find $\hat{u}^q \in \cU$ such that
	\begin{equation}\label{eq:incfwd_q}
	\langle \tilde{v},  \partial_{vu} \bar{r}_q\, \hat{u}^q \rangle = - \langle \tilde{v}, \partial_{vm} \bar{r}_q\, \hat{m}^q\rangle, \quad \forall \tilde{v} \in \cV,
	\end{equation} 
	where the derivatives $\partial_{vu} \bar{r}_q: \cU \mapsto \cV'$ and $\partial_{vm} \bar{r}_q: \cM \mapsto \cV'$ are linear operators.
	%More in general, $\partial_{ab}f: \cB \mapsto \cA'$ and $\partial_{ba} f: \cA \mapsto \cB'$ are linear operators with $a \in \cA$, and $b \in \cB$ representing any of $u \in \cU, v \in \cV, m \in \cM, z \in \cZ$, and $f$ representing one of $\bar{r}_q$ and $\bar{q}$, which satisfy 
	%\begin{equation}
	%\langle a, \partial_{ab} f \, b \rangle: = \;_\cA \langle a, \partial_{ab} f \, b \rangle_{\cA'} = \;_\cB \langle b, \partial_{ba} f \, a \rangle_{\cB'} = :\langle b, \partial_{ba} f \, a \rangle.
	%\end{equation}
	The incremental adjoint problem, obtained by taking variation of \eqref{eq:meta_lag} with respect to the state $u$ and using \eqref{eq:adjoint_op}, reads: find $\hat{v}^q \in \cV$ such that
	\begin{equation}\label{eq:incadj_q}
	%\begin{split}
	\langle \tilde{u}, \partial_{uv} \bar{r}_q \,\hat{v}^q \rangle  = - \langle \tilde{u}, \partial_{uu} \bar{r}_q\,\hat{u}^q + \partial_{uu} \bar{q}\,\hat{u}^q  + \partial_{um} \bar{r}_q\,\hat{m}^q \pc{+\partial_{um} \bar{q} \, \hat{m}^q}\rangle, \quad \forall \tilde{u} \in \cU,
	%\end{split}
	\end{equation}
	where $\partial_{uv} \bar{r}_q: \cV \mapsto \cU'$ is the adjoint of $\partial_{vu} \bar{r}_q: \cU \mapsto \cV'$ in the sense of \eqref{eq:adjoint_op}.
	%We remark that an extra term $-\langle \tilde{u}, \partial_{um} \bar{q} \,\hat{m}^q\rangle$ is added on the right hand side if $q$ explicitly depends on $m$.
	The Hessian of $q$ at $\bar{m}$ acting in direction $\hat{m}^q$ can then be computed by taking variation of \eqref{eq:meta_lag} with respect to $m$ and using \eqref{eq:adjoint_op} as 
	\begin{equation}\label{eq:Hessian_q}
	\begin{split}
	\langle \tilde{m}, \nabla_m^2 \bar{q}\,\hat{m}^q \rangle  & = \langle \tilde{m}, \partial_m \cL^H \rangle \\
	& =  \langle \tilde{m}, \partial_{mv}\bar{r}_q \, \hat{v}^q + \partial_{mu}\bar{r}_q\,\hat{u}^q \pc{+ \partial_{mu}\bar{q}\,\hat{u}^q } + \partial_{mm} \bar{r}_q\,\hat{m}^q \pc{+ \partial_{mm} \bar{q}\,\hat{m}^q}\rangle, \; \forall \tilde{m} \in \cM,
	\end{split}
	\end{equation} 
	where the incremental state $\hat{u}^q$ and adjoint $\hat{v}^q$ solve \eqref{eq:incfwd_q} and \eqref{eq:incadj_q}, respectively.
	%Note that $\langle \tilde{m}, \partial_{mu}\bar{q}\,\hat{u}^q + \partial_{mm} \bar{q}\,\hat{m}^q\rangle$ is also added if $q$ explicitly depends on $m$.
	
	Similarly, we can compute the gradient and Hessian action of $f$. For ease of notation we define
	\begin{equation}\label{eq:rf_bar}
	\bar{r}_f = r(u,v^f,\bar{m},z), \text{ and } \bar{f} = f(u(\bar{m}, z),\bar{m},z).
	\end{equation}
	By solving the adjoint problem: 
	find $v^f \in \cV$, such that
	\begin{equation}\label{eq:adjoint_f}
	\langle \tilde{u} , \partial_u \bar{r}_f\rangle = - \langle \tilde{u}, \partial_u \bar{f} \rangle, \quad \forall \tilde{u} \in \cU,
	\end{equation}
	we obtain the gradient of $f$ at $\bar{m}$ as 
	\begin{equation}\label{eq:gradient_f}
	\langle \tilde{m}, \nabla_m \bar{f} \rangle  =\langle \tilde{m}, \pc{\partial_m \bar{f}} + \partial_m \bar{r}_f\rangle, \quad \forall \tilde{m} \in \cM.
	\end{equation}
	Then by solving the incremental state problem: 
	find $\hat{u}^f \in \cU$ such that
	\begin{equation}\label{eq:incfwd_f}
	\langle \tilde{v},  \partial_{vu} \bar{r}_f\, \hat{u}^f \rangle = - \langle \tilde{v}, \partial_{vm} \bar{r}_f\, \hat{m}^f\rangle, \quad \forall \tilde{v} \in \cV,
	\end{equation} 
	and the incremental adjoint problem: find $\hat{v}^f \in \cV$ such that
	\begin{equation}\label{eq:incadj_f}
	%\begin{split}
	\langle \tilde{u}, \partial_{uv} \bar{r}_f \,\hat{v}^f \rangle  = - \langle \tilde{u}, \partial_{uu} \bar{r}_f\,\hat{u}^f + \partial_{uu} \bar{f}\,\hat{u}^f  + \partial_{um} \bar{r}_f\,\hat{m}^f \pc{ + \partial_{um} \bar{f}\,\hat{m}^f}\rangle, \quad \forall \tilde{u} \in \cU,
	%\end{split}
	\end{equation}
	we obtain the Hessian action of $f$ at $\bar{m}$ in direction $\hat{m}^f$ as
	\begin{equation}\label{eq:Hessian_f}
	% \begin{split}
	\langle \tilde{m}, \nabla_m^2 \bar{f}\,\hat{m}^f\rangle =  \langle \tilde{m}, \partial_{mv}\bar{r}_f \, \hat{v}^f + \partial_{mu}\bar{r}_f\,\hat{u}^f \pc{+ \partial_{mu}\bar{f}\,\hat{u}^f} + \partial_{mm} \bar{r}_f\,\hat{m}^f \pc{+ \partial_{mm} \bar{f}\,\hat{m}^f}\rangle, \quad \forall \tilde{m} \in \cM.
	% \end{split}
	\end{equation} 
	
	\section{Gradient-based optimization}
	\label{sec:gradientopt}
	
	In this section, we develop a gradient based optimization method to solve the chance constrained optimization problem \eqref{eq:optimization}. The method employs (1) a smooth approximation of the indicator function involved in the chance evaluation, (2) an exterior penalty method for the inequality chance constraint, (3) a continuation scheme to refine the smooth approximation and the penalty for inequality constraint, and (4) an approximate cost functional, with both SAA and Taylor approximations, and their gradients with respect to the optimization variable.
	
	\subsection{Smooth approximation of the indicator function} 
\begin{figure}[!htb]
		\centering
		\includegraphics[width=0.45\textwidth]{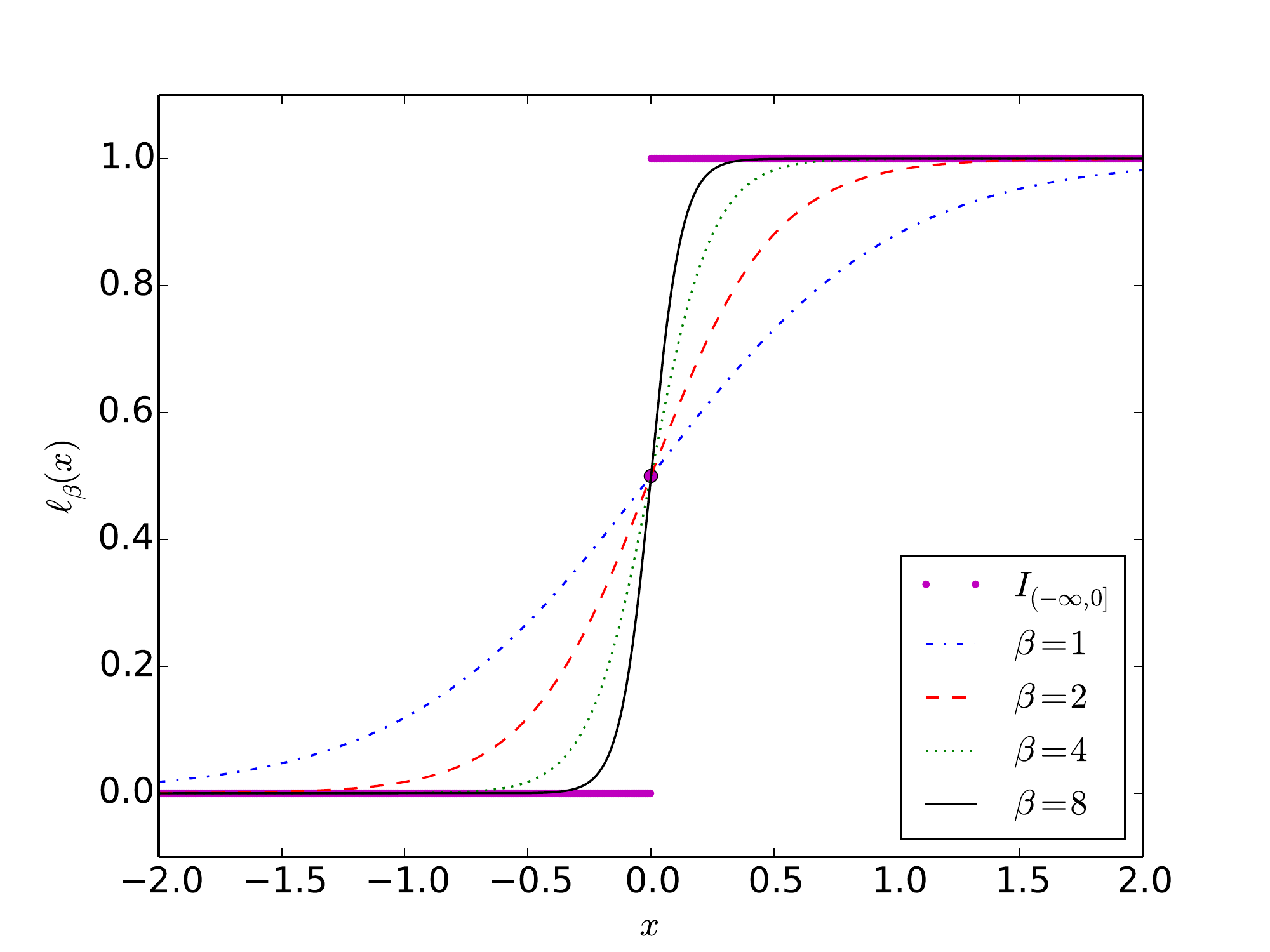}
		\includegraphics[width=0.45\textwidth]{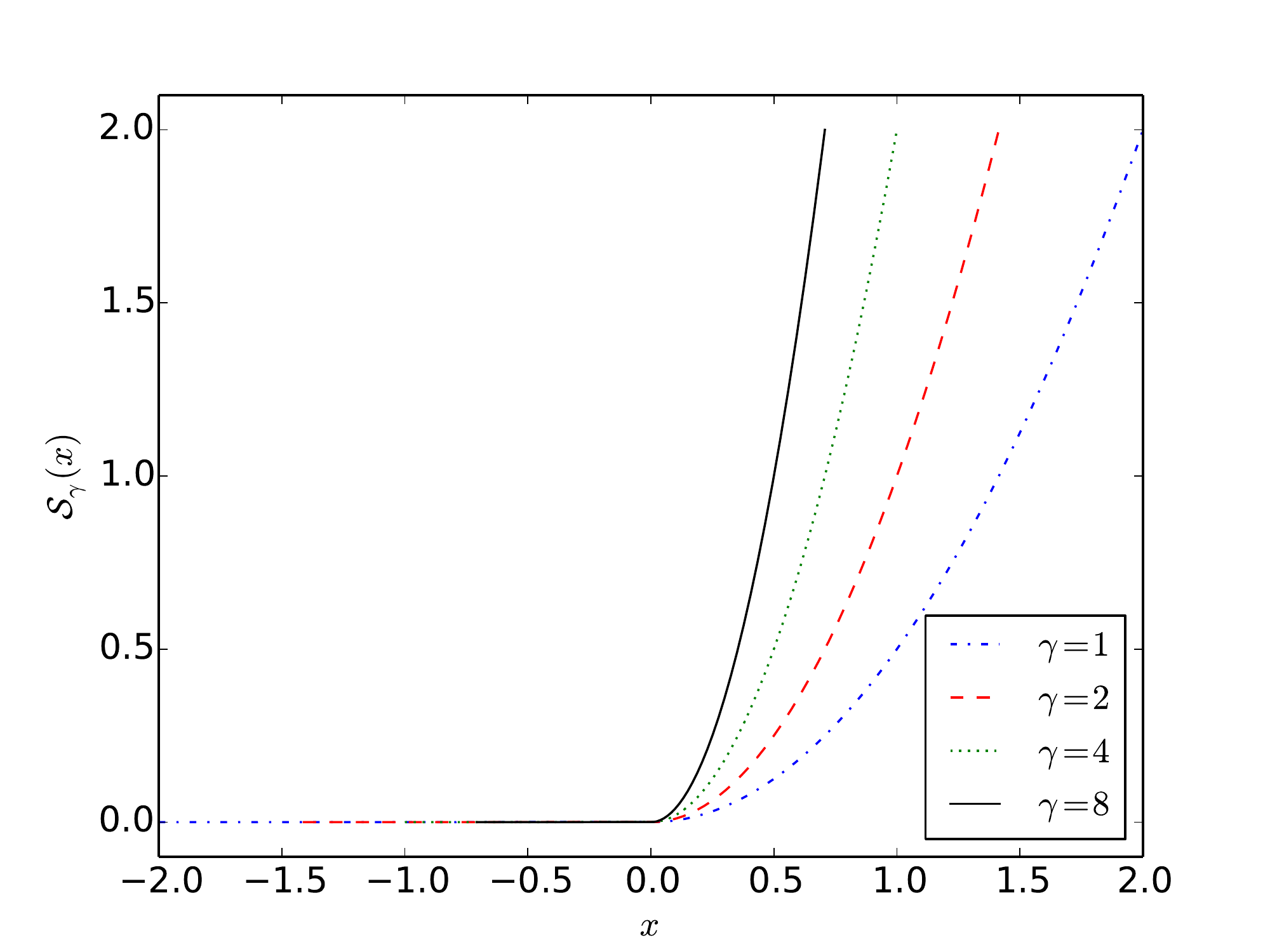}
		\caption{Left: Smooth approximation of the discontinuous indicator function $\bI_{[0, \infty)}(x)$ by a logistic function $\ell_\beta(x) = \frac{1}{1+ e^{-2\beta x}}$ with $\beta > 0$. Right: Penalty function $\cS_\gamma(x) = \frac{\gamma}{2}(\max\{0, x\})^2$ with $\gamma > 0$.}
		\label{fig:logistic}
	\end{figure}

The evaluation of the probability \eqref{eq:probability_definition} involves the indicator function $\bI_{[0, \infty)}(f(m, z))$, which is discontinuous at $f(m, z) = 0$. To use a gradient-based optimization method, we consider a smooth approximation of the indicator function by a logistic function 
	\begin{equation}\label{eq:logistic}
	\bI_{[0, \infty)}(x) \approx \ell_\beta(x) = \frac{1}{1+e^{-2\beta x}},
	\end{equation}
	where a larger $\beta$ corresponds to a sharper transition at $x = 0$, as shown in Figure \ref{fig:logistic}. With the definition $\bI_{[0, \infty)}(0) = \frac{1}{2}$, we have the convergence 
	\begin{equation}\label{eq:logistic_grad}
	\lim_{\beta \to \infty} \ell_\beta(x) = \bI_{[0, \infty)}(x) \text{ and }
	\lim_{\beta \to \infty} \nabla \ell_\beta(x) = \frac{2\beta e^{-2\beta x}}{(1+ e^{-2\beta x})^2} = \nabla \bI_{[0, \infty)}(x).
	\end{equation}
	
	\subsection{A penalty method for inequality constraint optimization}
	To solve the optimization problem \eqref{eq:optimization} with the inequality constraint \eqref{eq:probability_constraint}, we employ a quadratic penalty method \cite{NocedalWright06} by first defining 
	% \begin{equation}
	% g(f(z)) := \alpha - P(f(\cdot,z) \leq 0),
	% \end{equation}
	% which is approximated by 
	% \begin{equation}
	% g_\beta(f(z)) := \alpha - \bE[\ell_\beta(f(\cdot, z))]
	% \end{equation}
	% under the smooth approximation of the indicator function \eqref{eq:logistic}.
	the exterior penalty function
	\begin{equation}\label{eq:penalty}
	\cS_\gamma(x) := \frac{\gamma}{2}\left(\max\{0, x\}\right)^2, \text{ with } \nabla \cS_\gamma(x) = \gamma \max\{0, x\},
	\end{equation}
	for a constant $\gamma > 0$ controlling the weight of the penalty. 
	% The gradient of $\cS_\gamma(g_k(z))$ w.r.t.\ $z$ is given by
	% \begin{equation}
	% \nabla_z \cS_\gamma(g_k(z)) = \gamma \max\{0, g_k(z)\} \nabla_z g_k(z).
	% \end{equation} 
	%, so that the penalty function $\cS_\gamma(g_k(z))$ is continuous w.r.t.\ $z$.
	Then the chance constrained optimization problem \eqref{eq:optimization} can be approximated by the unconstrained problem 
	\begin{equation}\label{eq:OptUnC}
	\min_{z \in \cZ} \cJ(z) + \cS_\gamma( \bE[\ell_\beta(f)]  - \alpha).
	\end{equation}
	
	For general inequality constrained optimization problems, convergence of the optimization variable as $\gamma \to \infty$ by the penalty method is studied in, e.g., \cite{NocedalWright06}. 
	
	%\pc{to be further developed on convergence analysis when $\gamma, \beta \to \infty$}
	
	\subsection{Adaptive BFGS optimization}
	Let $\cE(z)$ denote an approximate cost functional, which is an approximation of the cost functional in \eqref{eq:OptUnC} by SAA or Taylor approximations. Let $\nabla_z \cE(z)$ denote its gradient with respect to the optimization variable, which is computed in next section. We present a continuation optimization scheme by increasing the smoothing parameter $\beta$ and the penalty parameter $\gamma$ in an outer loop and applying a gradient-based quasi Newton optimization algorithm, BFGS \cite{NocedalWright06, MoralesNocedal11}, to solve the optimization problem \eqref{eq:OptUnC} (with possible bound constraint on the optimization variable $z$) using the approximate cost functional $\cE(z)$ and gradient $\nabla_z \cE(z)$ in an inner loop. The optimization procedure is presented in Algorithm \ref{alg:adaptiveBFGS}. For the continuation, we specify an initial smoothing parameter $\beta_0$ and penalty parameter $\gamma_0$, and scale them with the power parameters $\sigma_\beta$ and $\sigma_\gamma$  in Line 6. We stop the outer loop if the maximum number of iterations is reached or the difference between the critical chance  $\alpha$  and the approximate chance, defined in \eqref{eq:ChanceSAA} for SAA or \eqref{eq:ChanceTaylor} for the Taylor approximation of the constraint function, is smaller than a tolerance; see Line 4. The inner loop of BFGS optimization is stopped if the maximum number of iteration is reached or the gradient of the approximate cost functional is smaller than a tolerance; see Line 5. 
	% Note that if a bound is present on the optimization variable $z$, as in our numerical experiment, then $z_k + \alpha p_k$ is projected back to the feasible region for $\alpha$ obtained in Line 9 when the step goes outside the feasible region.
	
	\begin{algorithm}[!htb]
		\caption{Continuation BFGS for chance constrained optimization}
		\label{alg:adaptiveBFGS}
		\begin{algorithmic}[1]
			\STATE{\textbf{Input}: initial value $z_0$, bound constraint $[z_{\text{min}}, z_{\text{max}}]$, smoothing parameter $\beta_0$, penalty parameter $\gamma_0$, power parameters $\sigma_\beta$ and $\sigma_\gamma$, maximum number of iterations $k_{\text{max}}$ and $l_{\text{max}}$, and tolerances $\varepsilon_{\text{in}}$ and $\varepsilon_{\text{out}}$ for the inner and outer loops.}
			\STATE{\textbf{Output}: the optimal variable $z_{\text{opt}}$.}
			\STATE{Set $l = 0$, compute the approximate chance $\hat{f}$ at $z_0$, where $\hat{f} = f_M(z_0)$ in \eqref{eq:ChanceSAA} for SAA or $\hat{f} = f_M^K(z_0)$ in \eqref{eq:ChanceTaylor} for the Taylor approximation of the constraint function.}
			\WHILE{$l \leq l_{\text{max}}$ or $|\hat{f}- \alpha| \geq \varepsilon_{\text{out}}$}
			%\STATE{Set $k = 0$, compute $\cE(z_k)$ and $\nabla_z \cE(z_k)$.}
			%\STATE{Set an initial guess of the Hessian $B_0 = \nabla_z^2 \cP(z_0)$.}
			%\WHILE{$k \leq k_{\text{max}}$ or $||\nabla_z \cE(z_k)|| \geq \varepsilon_{\text{in}}$}
			%\STATE{Solve the linear system $p_k = - B_k^{-1} \nabla_z \cE(z_k)$ to get $p_k$.}
			%\STATE{Perform a line search to get \ $\alpha_k = \argmin_{\alpha} \cE(z_k + \alpha p_k)$.}
			%\STATE{Compute the gradient $\nabla_z \cE(z_{k}+\alpha_k p_k)$.}
			%\STATE{Set $s_k = \alpha_k p_k$, $z_{k+1} = z_k + s_k$, $y_k = \nabla_z \cE(z_{k+1}) - \nabla_z \cE(z_{k})$.}
			%\STATE{
			%Update the inverse of the Hessian as 
			%$$
			%B_{k+1}^{-1} = B_k^{-1} + \frac{(s_k^Ty_k + y_k^TB_k^{-1}y_k)(s_ks_k^T)}{(s_k^T y_k)^2} - \frac{B_k^{-1} y_k s_k^T + s_k y_k^TB_k^{-1}}{s_k^T y_k}.
			%$$
			%}
			%\STATE{Set $k \leftarrow k+1.$}
			%\ENDWHILE
			\STATE{Run an inner loop optimization by the Algorithm L-BFGS-B \cite{MoralesNocedal11}
				$$
				z_{\text{opt}} = \text{L-BFGS-B}(\cE(z), z_0, \nabla_z \cE(z), k_{\text{max}}, \varepsilon_{\text{in}}, [z_{\text{min}}, z_{\text{max}}]),
				$$
				with objective $\cE(z)$ and gradient $\nabla_z \cE(z)$ given by either SAA or Taylor approximation.
			}
			\STATE{Set $\beta_{l+1} = \sigma_\beta^{l+1} \beta_l$, $\gamma_{l+1} = \sigma_\gamma^{l+1} \gamma_l$, $l \leftarrow l+1$.}
			\STATE{Set $z_0 = z_{\text{opt}}$, and compute the approximate chance $\hat{f}$ at $z_{\text{opt}}$.}
			%\RETURN{$z_{\text{opt}} = z_{k+1}$.}
			\ENDWHILE
			% \RETURN{$z_{\text{opt}} = z_{k}$.}
		\end{algorithmic}
	\end{algorithm}
	
	\subsection{Computation of the approximate cost functional and its gradient}
	
	Algorithm \ref{alg:adaptiveBFGS} requires the computation of the approximate cost functional $\cE(z)$ and its gradient $\nabla_z \cE(z)$ at a given optimization variable $z$, constrained by the state PDE \eqref{eq:PDE}. In this section, we present their computation for both SAA and Taylor approximations.
	
	\subsubsection{Sample average approximation}
	
	By sample average approximation of the mean and chance presented in Section \ref{sec:MCinteg},
	the cost functional in \eqref{eq:OptUnC} can be approximated as
	\begin{equation}\label{eq:CostApproxSAA}
	\begin{split}
	& \cJ(z) + \cS_\gamma(\alpha - \bE[\ell_\beta(f(\cdot, z))]) \\
	& = \bE[q(\cdot, z)] + \cP(z) + \cS_\gamma( \bE[\ell_\beta(f(\cdot, z))] - \alpha)
	\\
	& \approx \frac{1}{M_q}\sum_{i = 1}^{M_q} q(m_i^q, z) + \cP(z) + \cS_\gamma
	\left(
	\frac{1}{M_f} \sum_{i=1}^{M_f} \ell_\beta(f(m_i^f,z)) - \alpha 
	\right) =: \cE(z),
	\end{split}
	\end{equation}
	where $M_q$ and $M_f$ independent random samples are taken such that the SAA errors of the two approximate terms are balanced. A practical approach to determine $M_q$ and $M_f$ is to first evaluate the variances of $q$ and $\ell_\beta(f)$ and then set the two numbers with the ratio equal to that of the variances. For simplicity, one can use the same $M_q = M_f$ random samples. Note that to compute $q(m_i^q, z) = q(u(m_i^q), m_i^q, z)$ and $f(m_i^f, z) = f(u(m_i^f), m_i^f, z)$, we need to solve the state equation \eqref{eq:PDE} at $m_i^q$ and $m_i^f$, respectively.
	
	To compute the gradient of this approximate cost functional, we define the Lagrangian
	\begin{equation}\label{eq:LagGradzSAA}
	\begin{split}
	& \cL_{SAA} \big(z, \{u_i^q\}, \{u_i^f\}, \{v_i^q\}, \{v_i^f\}\big) \\
	& := \frac{1}{M_q}\sum_{i = 1}^{M_q} q(m_i^q, z) + \cP(z) + \cS_\gamma
	\left(
	\frac{1}{M_f} \sum_{i=1}^{M_f} \ell_\beta(f(m_i^f,z)) - \alpha 
	\right) \\
	& + \sum_{i=1}^{M_q} \langle v_i^q,  \partial_v r(u_i^q, v, m_i^q, z)\rangle +  \sum_{i=1}^{M_f} \langle v_i^f,  \partial_v r(u_i^f, v, m_i^f, z)\rangle,
	\end{split}
	\end{equation}
	where $ \{v_i^q\}$ and $\{v_i^f\}$ are the Lagrangian multipliers, which can be computed by solving the adjoint problems obtained by setting the variation of the Lagrangian w.r.t.\ $u_i^q$ to zero, i.e., for each $i = 1, \dots, M_q$, find $v_i^q \in \cV$ such that 
	\begin{equation}
	\langle \tilde{u}, \partial_{uv}r(u_i^q, v, m_i^q, z) v_i^q \rangle = - \frac{1}{M_q}\langle \tilde{u}, \partial_u q(u_i^q, m_i^q, z) \rangle, \quad \forall \tilde{u} \in \cU,
	\end{equation}
	where we recall that $q(u_i^q, m_i^q, z) = q(m_i^q, z)$ by a slight abuse of notation. Similarly, for each $i = 1, \dots, M_f$, find $v_i^f \in \cV$ such that 
	\begin{equation}
	\langle \tilde{u}, \partial_{uv}r(u_i^f, v, m_i^f, z) v_i^f \rangle = - \frac{1}{M_f}\langle \tilde{u}, 
	\nabla \cS_\gamma \nabla \ell_\beta \partial_u f(u_i^f, m_i^f, z)
	\rangle, \quad \forall \tilde{u} \in \cU. 
	\end{equation}
	Then the gradient of the approximate cost functional can be evaluated as 
	\begin{equation}
	\begin{split}
	\langle \tilde{z}, \nabla \cE(z) \rangle &= \langle \tilde{z}, \partial_z \cL_{SAA} \rangle
	\\
	& = \frac{1}{M_q}\sum_{i = 1}^{M_q} \langle \tilde{z}, \partial_z q(m_i^q, z) \rangle +  \langle \tilde{z}, \nabla \cP(z) \rangle + \frac{1}{M_f} \sum_{i=1}^{M_f} \langle \tilde{z}, \nabla \cS_\gamma \nabla \ell_\beta \partial_z f (m_i^f, z) \rangle \\
	& +  \sum_{i=1}^{M_q} \langle \tilde{z}, \partial_{zv} r(u_i^q, v, m_i^q, z) v_i^q\rangle +  \sum_{i=1}^{M_f} \langle \tilde{z}, \partial_{zv} r(u_i^f, v, m_i^f, z) v_i^f\rangle.
	\end{split}
	\end{equation}
	In summary, $M_q + M_f$ state PDEs and $M_q + M_f$ linearized (adjoint) PDEs have to be solved to evaluate the SAA of $\cE(z)$ and its gradient $\nabla \cE(z)$ at any given $z$.
	
	\subsubsection{Taylor approximation}
	For simplicity, we present only the approximate cost functional and its gradient by the quadratic Taylor approximation $T_2$ of the objective and constraint functions. The computation corresponding to the constant and linear approximations $T_0$ and $T_1$, which are contained within the quadratic approximation, are omitted.
	%The cost functional with constant approximation is given by 
	%\begin{equation}\label{eq:CostApproxT0}
	%\begin{split}
	% \cJ(z) + \cS_\gamma(\alpha - \bE[\ell_\beta(f(\cdot, z))]) & \approx \bE[T_0 q(\cdot, z)] + \cP(z) + \cS_\gamma( \bE[\ell_\beta(T_0 f(\cdot, z))] - \alpha)\\
	%& \approx \bar{q} + \cP(z) + \cS_\gamma
	%\left(
	%\ell_\beta(\bar{f}) - \alpha
	%\right)
	%\end{split}
	%\end{equation}
	%
	%The Lagrangian is written as 
	%
	%\begin{equation}\label{eq:LagGradzT0}
	%\begin{split}
	%& \cL_0\big(z, u, v^*\big) \\
	% & = \;\bar{q} + \cP(z) + \cS_\gamma\left(
	%\ell_\beta(\bar{f}) - \alpha
	%\right) \\
	%& + \langle v^*, \partial_v \bar{r}\rangle 
	%\end{split}
	%\end{equation}
	%
	%
	%The cost functional with linear Taylor approximation is given by 
	%\begin{equation}\label{eq:CostApproxT1}
	%\begin{split}
	% \cJ(z) + \cS_\gamma(\alpha - \bE[\ell_\beta(f(\cdot, z))]) & \approx \bE[T_1 q(\cdot, z)] + \cP(z) + \cS_\gamma( \bE[\ell_\beta(T_1f(\cdot, z))] - \alpha)\\
	%& \approx \bar{q} + \cP(z) + \cS_\gamma
	%\left(
	%g_{\beta,M_f}(T_1^{\text{LR}}f(\cdot, z))
	%\right)
	%\end{split}
	%\end{equation}
	%
	%The Lagrangian is written as 
	%
	%\begin{equation}\label{eq:LagGradzT1}
	%\begin{split}
	%& \cL_1\big(z, u, v^*, v^q,  (u^q)^*, v^f, (u^f)^*\big) \\
	% & = \;\bar{q} + \cP(z) + \cS_\gamma\left(
	%g_{\beta,M_f}(T_1^{\text{LR}}f)
	%\right) \\
	%& + \langle v^*, \partial_v \bar{r}\rangle \\
	%& + \langle (u^q)^* , \partial_u \bar{r}_q+ \partial_u \bar{q} \rangle \\
	%& + \langle (u^f)^* , \partial_u \bar{r}_f + \partial_u \bar{f} \rangle 
	%\end{split}
	%\end{equation}
	
	Using the quadratic Taylor approximations of the objective function $q$ and the constraint function $f$ in \eqref{eq:Taylor}, the low-rank approximations of the trace in \eqref{eq:TraceApprox} and the Hessian in \eqref{eq:HessianApprox}, as well as the sample average approximation for the probability \eqref{eq:ChanceSAA},  
	the cost functional in the unconstrained optimization problem \eqref{eq:OptUnC} becomes 
	\begin{equation}\label{eq:CostApprox}
	\begin{split}
	\cJ(z) + \cS_\gamma(\alpha - \bE[\ell_\beta(f(\cdot, z))]) & \approx \bE[T_2 q(\cdot, z)] + \cP(z) + \cS_\gamma(\bE[\ell_\beta(T_2f(\cdot, z))] - \alpha)\\
	& \approx \bar{q} + \frac{1}{2}\sum_{n=1}^{N_q} \lambda_n^q + \cP(z) + \cS_\gamma
	\left(
	g_{\beta,M_f}(T_2^{\text{LR}}f(\cdot, z))
	\right) \\
	& =: \cE(z)
	\end{split}
	\end{equation}
	where for simplicity we denote
	\begin{equation}
	g_{\beta,M_f}(T_2^{\text{LR}}f(\cdot, z)) =  \frac{1}{M_f}\sum_{i = 1}^{M_f} \ell_\beta(T_2^{\text{LR}} f(m_i,z)) - \alpha,
	\end{equation}
	where the quadratic Taylor approximation with low-rank decomposition $T_2^{\text{LR}} f (m_i,z)$ is given in \eqref{eq:TaylorLR}.
	Note that $\bar{q}$ and $\lambda_n^q$, $n = 1, \dots, N_q$, which are part of the approximation of $\bE[q]$, as well as $\bar{f}$, $\nabla_m \bar{f}$, $\lambda_n^f$ and $\psi_n^f$, $n = 1, \dots, N_f$, which are part of the approximation of $\bE[\ell_\beta(f)]$, all implicitly (possibly also explicitly) depend on the optimization variable $z$ through the state equation \eqref{eq:state}, the adjoint equations \eqref{eq:adjoint_q} and \eqref{eq:adjoint_f}, the generalized eigenvalue problems \eqref{eq:GenEigen_q} and \eqref{eq:GenEigen_f} with orthonormal conditions \eqref{eq:Orthon_q} and \eqref{eq:Orthon_f}, the incremental state and adjoint equations \eqref{eq:incfwd_q} and \eqref{eq:incadj_q} for the incremental state $\hat{u}^q = \hat{u}_n^q$ and adjoint $\hat{v}^q = \hat{v}_n^q$ at $\hat{m}^q = \psi_n^q$ needed by the Hessian action $\nabla_m^2 \bar{q} \, \psi_n^q$ in \eqref{eq:GenEigen_q} through \eqref{eq:Hessian_q}, $n = 1, \dots, N_q$, as well as the incremental state and adjoint equations \eqref{eq:incfwd_f} and \eqref{eq:incadj_f} for the incremental state $\hat{u}^f = \hat{u}_n^f$ and adjoint $\hat{v}^f = \hat{v}_n^f$ at $\hat{m}^f = \psi_n^f$ needed by the Hessian action $\nabla_m^2 \bar{f} \, \psi_n^f$ in \eqref{eq:GenEigen_f} through \eqref{eq:Hessian_f}, $n = 1, \dots, N_f$. To derive the gradient of the approximate cost functional \eqref{eq:CostApprox}, we define a Lagrangian $\cL_2$ to enforce all of the PDE constraint equations as follows:
	
	\allowdisplaybreaks
	\begin{align}\label{eq:LagGradz}
	%\begin{split} %\label{eq:LagGradz}
	& \cL_2\big(z, u, v^*, v^q,  (u^q)^*, \{\lambda_n^q\}, \{\psi_n^q\}, \{\hat{u}_n^q\}, \{\hat{v}_n^q\}, \{(\lambda_{nn'}^q)^*\}, \{(\psi_n^q)^*\}, \{(\hat{u}_n^q)^*\}, \{(\hat{v}_n^q)^*\}, \\
	& \qquad \qquad \quad v^f, (u^f)^*, \{\lambda_n^f\}, \{\psi_n^f\}, \{\hat{u}_n^f\}, \{\hat{v}_n^f\},  \{(\lambda_{nn'}^f)^*\}, \{(\psi_n^f)^*\}, \{(\hat{u}_n^f)^*\}, \{(\hat{v}_n^f)^*\}\big) \nonumber\\
	& := \;\bar{q} + \frac{1}{2}\sum_{n=1}^{N_q} \lambda_n^q + \cP(z) + \cS_\gamma\left(
	g_{\beta,M_f}(T_2^{\text{LR}}f)
	\right) \nonumber\\
	& + \langle v^*, \partial_v \bar{r}\rangle \nonumber\\
	& + \langle (u^q)^* , \partial_u \bar{r}_q+ \partial_u \bar{q} \rangle \nonumber\\
	& + \sum_{n = 1}^{N_q}\langle (\psi_n^q)^*,  \nabla_m^2 \bar{q} \, \psi_n^q - \lambda_n^q \cC^{-1} \psi_n^q \; \rangle \nonumber\\
	& + \sum_{n, n'= 1}^{N_q} (\lambda_{nn'}^q)^* (\langle\psi_n^q, \cC^{-1}\psi_{n'}^q \rangle - \delta_{nn'}) \nonumber\\
	& + \sum_{n = 1}^{N_q} \langle (\hat{v}_n^q)^*,  \partial_{vu} \bar{r}_q\, \hat{u}_n^q + \partial_{vm} \bar{r}_q\, \psi_n^q \rangle \nonumber\\
	& + \sum_{n = 1}^{N_q} \langle (\hat{u}_n^q)^*, \partial_{uv} \bar{r}_q \,\hat{v}^q_n + \partial_{uu} \bar{r}_q\,\hat{u}^q_n + \partial_{uu} \bar{q}\,\hat{u}^q_n  + \partial_{um} \bar{r}_q\, \psi_n^q \pc{+ \partial_{um} \bar{q}\, \psi_n^q}\rangle \nonumber\\[8pt]
	& + \langle (u^f)^* , \partial_u \bar{r}_f + \partial_u \bar{f} \rangle \nonumber\\
	& + \sum_{n = 1}^{N_f}\langle (\psi_n^f)^*,  \nabla_m^2 \bar{f} \, \psi_n^f - \lambda_n^f \cC^{-1} \psi_n^f \; \rangle \nonumber\\
	& + \sum_{n,n' = 1}^{N_f} (\lambda_{nn'}^f)^* (\langle\psi_n^f, \cC^{-1}\psi_{n'}^f \rangle - \delta_{nn'}) \nonumber\\
	& + \sum_{n = 1}^{N_f} \langle (\hat{v}_n^f)^*,  \partial_{vu} \bar{r}_f\, \hat{u}_n^f + \partial_{vm} \bar{r}_f\, \psi_n^f\rangle \nonumber\\
	& + \sum_{n = 1}^{N_f} \langle (\hat{u}_n^f)^*, \partial_{uv} \bar{r}_f \,\hat{v}^f_n + \partial_{uu} \bar{r}_f\,\hat{u}^f_n + \partial_{uu} \bar{f}\,\hat{u}^f_n  + \partial_{um} \bar{r}_f\, \psi_n^f \pc{+ \partial_{um} \bar{f}\, \psi_n^f}\rangle.
	%\end{split}
	\end{align}
	% % where for simplicity we denote $\{\lambda_n^q\} = \{\lambda_n^q, n = 1, \dots, N_q\}$, $\{\lambda_n^f\} = \{\lambda_n^f, n = 1, \dots, N_f\}$, and similarly for other quantities in the bracket $\{\cdot\}$.
	% % where $v^*, (u^q)^*, \{(\lambda_n^q)^*\}, \{(\psi_n^q)^*\}, \{(u_n^q)^*\}, \{(v_n^q)^*\}, (u^f)^*, \{(\lambda_n^f)^*\}, \{(\psi_n^f)^*\}, \{(u_n^f)^*\}, \{(v_n^f)^*\}$ are the Lagrange multipliers. 
	
	% Here the dominating eigenvalues $\{\lambda_n^q\}$ (and $\{\lambda_n^f\}$) typically decay rapidly, so that we can assume that they are not repeated, i.e., $\lambda_n^q \neq \lambda_{n'}^q$ for $n \neq n'$, $n, n' = 1, \dots, N_q$. Therefore, the constraints $\langle \psi_n^q, \cC^{-1} \psi_n^q \rangle = 1$, $n = 1, \dots, N_q$, in \eqref{eq:LagGradz}, are sufficient to guarantee the orthonormality $\langle \psi_n, \cC^{-1} \psi_{n'} \rangle = \delta_{nn'}$ for any $n,n' = 1, \dots, N_q$ in \eqref{eq:Orthon_q}. In fact, for any $n \neq n'$, by definition we have 
	% \begin{equation}
	% \langle \psi_{n'}, \nabla_m^2 q(\bar{m}) \psi_{n} \rangle = \langle \psi_{n'}, \lambda_{n} C^{-1} \psi_{n} \rangle \text{ and } \langle \psi_{n}, \nabla_m^2 q(\bar{m}) \psi_{n'} \rangle = \langle \psi_{n}, \lambda_{n'} \cC^{-1} \psi_{n'} \rangle, 
	% \end{equation} 
	% which, by the symmetry of $\nabla_m^2 q(\bar{m})$ and $\cC^{-1}$, leads to 
	% \begin{equation}
	% (\lambda_{n} - \lambda_{n'}) \langle \psi_{n}, \cC^{-1} \psi_{n'} \rangle = 0, \text{ so that } \langle \psi_{n}, \cC^{-1} \psi_{n'} \rangle = 0 \text{ if }\lambda_{n} \neq \lambda_{n'}.
	% \end{equation}
	
	We compute all the Lagrange multipliers by setting the variation of the Lagrangian $\cL_2$ to zero. Specifically, by setting the  variation of the Lagrangian $\cL_2$ 
	w.r.t.\ $\lambda_n^q$, $n = 1, \dots, N_q$, to zero, we obtain 
	\begin{equation}\label{eq:PhiStar_q}
	(\psi_n^q)^* = \frac{1}{2} \psi_n^q, \quad n = 1, \dots, N_q.
	\end{equation}
	Subsequently, for each $n = 1, \dots, N_q$, by setting the variation of $\cL_2$ w.r.t.\ $\hat{v}_n^q$ to zero, using the Hessian action \eqref{eq:Hessian_q} and \eqref{eq:adjoint_op}, we have: find $(\hat{u}_n^q)^* \in \cU$ such that 
	\begin{equation}
	\langle \tilde{v}, \partial_{vu} \bar{r}_q (\hat{u}_n^q)^* \rangle = - \langle \tilde{v}, \partial_{vm} \bar{r}_q (\psi_n^q)^*\rangle, \quad \forall \tilde{v} \in \cV, 
	\end{equation}
	which, together with \eqref{eq:PhiStar_q} and \eqref{eq:incfwd_q} with $(\hat{u}^q, \hat{m}^q) = (\hat{u}_n^q, \psi_n^q)$, leads to 
	\begin{equation}\label{eq:IncStateStar_q}
	(\hat{u}_n^q)^* = \frac{1}{2} \hat{u}_n^q, \quad n = 1, \dots, N_q.
	\end{equation}
	Similarly, for each $n = 1, \dots, N_q$, by setting the variation of $\cL_2$ w.r.t.\ $\hat{u}_n^q$ to zero, using the Hessian action \eqref{eq:Hessian_q} and \eqref{eq:adjoint_op}, we have: find $(\hat{v}_n^q)^* \in \cU$ such that 
	\begin{equation}
	\langle \tilde{u}, \partial_{uv} \bar{r}_q \,(\hat{v}_n^q)^* \rangle  = - \langle \tilde{u}, \partial_{uu} \bar{r}_q\,(\hat{u}_n^q)^* + \partial_{uu} \bar{q}\,(\hat{u}_n^q)^*  + \partial_{um} \bar{r}_q\,(\psi_n^q)^* \pc{+ \partial_{um} \bar{q}\,(\psi_n^q)^*}\rangle, \quad \forall \tilde{u} \in \cU,
	\end{equation}
	which, together with \eqref{eq:PhiStar_q}, \eqref{eq:IncStateStar_q}, and \eqref{eq:incadj_q} with $(\hat{v}^q, \hat{u}^q, \hat{m}^q) = (\hat{v}_n^q, \hat{u}_n^q, \psi_n^q)$, leads to 
	\begin{equation}\label{eq:IncAdjointStar_q}
	(\hat{v}_n^q)^* = \frac{1}{2} \hat{v}_n^q, \quad n = 1, \dots, N_q.
	\end{equation}
	Then by setting the variation of $\cL_2$ w.r.t.\ $v^q$ to zero, we obtain: find $(u^q)^* \in \cU$ such that 
	\begin{equation}\label{eq:uStar_q}
	\begin{split}
	\langle \tilde{v}, \partial_{vu} \bar{r}_q (u^q)^* \rangle  = & - \sum_{n=1}^{N_q} \langle \tilde{v}, \partial_{vmu} \bar{r}_q \hat{u}_n^q (\psi_n^q)^* + \partial_{vmm} \bar{r}_q \psi_n^q (\psi_n^q)^* \rangle \\
	& - \sum_{n=1}^{N_q} \langle \tilde{v}, \partial_{vuu} \bar{r}_q \hat{u}_n^q (\hat{u}_n^q)^* + \partial_{vum} \bar{r}_q \psi_n^q (\hat{u}_n^q)^* \rangle, \quad \tilde{v} \in \cV.
	\end{split}
	\end{equation}
	
	By setting the $\cL_2$ w.r.t.\ $\lambda_n^f$, $n = 1, \dots, N_f$, to zero, we obtain 
	\begin{equation}\label{eq:PhiStar_f}
	(\psi_n^f)^* = c^f_n \psi_n^f, \quad n = 1, \dots, N_f,
	\end{equation}
	where the constant $c^f_n$ is given by 
	\begin{equation}
	c_n^f =  \nabla \cS_\gamma(g_{\beta, M_f}(T_2^{\text{LR}}f))  \frac{1}{2M_f}
	\sum_{i = 1}^{M_f} \nabla \ell_\beta (T_2^{\text{LR}}f(m_i,z))  \langle m_i-\bar{m},  \cC^{-1} \psi_n^f\rangle^2,
	\end{equation}
	where $\nabla \cS_\gamma$ and $\nabla \ell_\beta$ are defined in \eqref{eq:logistic_grad} and \eqref{eq:penalty}, respectively. 
	Subsequently, for each $n = 1, \dots, N_f$, by setting the variation of $\cL_2$ w.r.t.\ $\hat{v}_n^f$ to zero, and using the Hessian action \eqref{eq:Hessian_f} and \eqref{eq:adjoint_op}, we have: find $(\hat{u}_n^f)^* \in \cU$ such that 
	\begin{equation}
	\langle \tilde{v}, \partial_{vu} \bar{r}_f (\hat{u}_n^f)^* \rangle = - \langle \tilde{v}, \partial_{vm} \bar{r}_f (\psi_n^f)^*\rangle, \quad \forall \tilde{v} \in \cV, 
	\end{equation}
	which, together with \eqref{eq:PhiStar_f} and \eqref{eq:incfwd_f} with $(\hat{u}^f, \hat{m}^f) = (\hat{u}_n^f, \psi_n^f)$, leads to 
	\begin{equation}\label{eq:IncStateStar_f}
	(\hat{u}_n^f)^* = c_n^f \hat{u}_n^f, \quad n = 1, \dots, N_f.
	\end{equation}
	Similarly, for each $n = 1, \dots, N_f$, by setting the variation of $\cL_2$ w.r.t.\ $\hat{u}_n^f$ to zero, using the Hessian action \eqref{eq:Hessian_f} and \eqref{eq:adjoint_op}, we have: find $(\hat{v}_n^f)^* \in \cU$ such that
	\begin{equation}
	\langle \tilde{u}, \partial_{uv} \bar{r}_f \,(\hat{v}_n^f)^* \rangle  = - \langle \tilde{u}, \partial_{uu} \bar{r}_f\,(\hat{u}_n^f)^* + \partial_{uu} \bar{q}\,(\hat{u}_n^f)^*  + \partial_{um} \bar{r}_f\,(\psi_n^f)^* \pc{+ \partial_{um} \bar{f}\,(\psi_n^f)^*} \rangle, \quad \forall \tilde{u} \in \cU,
	\end{equation}
	which, together with \eqref{eq:PhiStar_f}, \eqref{eq:IncStateStar_f}, and \eqref{eq:incadj_f} with $(\hat{v}^f, \hat{u}^f, \hat{m}^f) = (\hat{v}_n^f, \hat{u}_n^f, \psi_n^f)$, leads to 
	\begin{equation}\label{eq:IncAdjointStar_f}
	(\hat{v}_n^f)^* = c_n^f \hat{v}_n^f, \quad n = 1, \dots, N_f.
	\end{equation}
	Then by setting the variation of $\cL_2$ w.r.t.\ $v^f$ to zero, we obtain: find $(u^f)^* \in \cU$ such that 
	\begin{equation}\label{eq:uStar_f}
	\begin{split}
	\langle \tilde{v}, \partial_{vu} \bar{r}_f (u^f)^* \rangle  = & - \sum_{n=1}^{N_f} \langle \tilde{v}, \partial_{vmu} \bar{r}_f \hat{u}_n^f (\psi_n^f)^* + \partial_{vmm} \bar{r}_f \psi_n^f (\psi_n^f)^* \rangle \\
	& - \sum_{n=1}^{N_f} \langle \tilde{v}, \partial_{vuu} \bar{r}_f \hat{u}_n^f (\hat{u}_n^f)^* + \partial_{vum} \bar{r}_f \psi_n^f (\hat{u}_n^f)^* \rangle\\
	& - \langle \tilde{v}, \partial_{vm}\bar{r}_f m^f \rangle, \quad \tilde{v} \in \cV,
	\end{split}
	\end{equation}
	where the last term is due to the gradient \eqref{eq:gradient_f}, which appears in the quadratic Taylor approximation \eqref{eq:TaylorLR} that is used in $\cS_\gamma
	(g_{\beta,M_f}(T_2^{\text{LR}}f )$, with $m^f$ given by 
	\begin{equation}
	m^f = \nabla \cS_\gamma(g_{\beta, M_f}(T_2^{\text{LR}}f)) \frac{1}{M_f} 
	\sum_{i = 1}^{M_f} \nabla \ell_\beta (T_2^{\text{LR}}f(m_i,z)) (m_i-\bar{m}).
	\end{equation}
	
	% \pc{the constant $c_n^f$ also depends on $v^f$} 
	
	Finally, by setting the variation of $\cL_2$ w.r.t.\ $u$ to zero, we obtain: find $v^* \in \cV$ such that 
	\begin{equation}\label{eq:vStar}
	\begin{split}
	& \langle \tilde{u}, \partial_{uv} \bar{r} v^* \rangle  
	= \\
	& - \langle \tilde{u}, \partial_u \bar{q}  + d^f \partial_u \bar{f} + \partial_{um} \bar{r}_f m^f \pc{+ \partial_{um} \bar{f} m^f}\rangle 
	\\
	& - \langle \tilde{u},
	\partial_{uu} \bar{r}_q (u^q)^* + \partial_{uu} \bar{q} (u^q)^*
	\rangle
	\\
	& -  \sum_{n=1}^{N_q} \langle \tilde{u},
	\partial_{umv} \bar{r}_q \hat{v}^q_n (\psi_n^q)^* + \partial_{umu} \bar{r}_q \hat{u}^q_n (\psi_n^q)^* \pc{+ \partial_{umu} \bar{q} \hat{u}^q_n (\psi_n^q)^*} + \partial_{umm} \bar{r}_q \psi_n^q (\psi_n^q)^* \pc{+ \partial_{umm} \bar{q} \psi_n^q (\psi_n^q)^*}
	\rangle 
	\\
	& -  \sum_{n=1}^{N_q} \langle \tilde{u},
	\partial_{uvu} \bar{r}_q \hat{u}_n^q (\hat{v}_n^q)^* + \partial_{uvm} \bar{r}_q \psi_n^q (\hat{v}_n^q)^* 
	\rangle
	\\
	& -  \sum_{n=1}^{N_q} \langle \tilde{u},
	\partial_{uuv} \bar{r}_q \hat{v}_n^q (\hat{u}_n^q)^* + \partial_{uuu} \bar{r}_q \hat{u}_n^q (\hat{u}_n^q)^* + \partial_{uuu} \bar{q} \hat{u}_n^q (\hat{u}_n^q)^* + \partial_{uum} \bar{r}_q \psi_n^q (\hat{u}_n^q)^* \pc{+ \partial_{uum} \bar{q} \psi_n^q (\hat{u}_n^q)^*}
	\rangle
	\\
	& - \langle \tilde{u},
	\partial_{uu} \bar{r}_f (u^f)^* + \partial_{uu} \bar{f} (u^f)^*
	\rangle
	\\
	& -  \sum_{n=1}^{N_f} \langle \tilde{u},
	\partial_{umv} \bar{r}_f \hat{v}^f_n (\psi_n^f)^* + \partial_{umu} \bar{r}_f \hat{u}^f_n (\psi_n^f)^* \pc{+ \partial_{umu} \bar{f} \hat{u}^f_n (\psi_n^f)^*} + \partial_{umm} \bar{r}_f \psi_n^f (\psi_n^f)^* \pc{+ \partial_{umm} \bar{f} \psi_n^f (\psi_n^f)^*}
	\rangle 
	\\
	& -  \sum_{n=1}^{N_f} \langle \tilde{u},
	\partial_{uvu} \bar{r}_f \hat{u}_n^f (\hat{v}_n^f)^* + \partial_{uvm} \bar{r}_f \psi_n^f (\hat{v}_n^f)^*
	\rangle
	\\
	& -  \sum_{n=1}^{N_f} \langle \tilde{u},
	\partial_{uuv} \bar{r}_f \hat{v}_n^f (\hat{u}_n^f)^* + \partial_{uuu} \bar{r}_f \hat{u}_n^f (\hat{u}_n^f)^* + \partial_{uuu} \bar{f} \hat{u}_n^f (\hat{u}_n^f)^* + \partial_{uum} \bar{r}_f \psi_n^f (\hat{u}_n^f)^* \pc{+ \partial_{uum} \bar{f} \psi_n^f (\hat{u}_n^f)^*}
	\rangle,  \; \forall \tilde{u} \in \cU,
	\end{split}
	\end{equation}
	where the constant $d^f$ in the first line is given by 
	\begin{equation}
	d^f = \nabla \cS_\gamma(g_{\beta, M_f}(T_2^{\text{LR}}f))  \frac{1}{M_f}
	\sum_{i = 1}^{M_f} \nabla \ell_\beta (T_2^{\text{LR}}f(m_i,z)).
	\end{equation}
	
	% \pc{the constants $c_n^f$ and $m^f$ also depend on $u$}
	
	With all the Lagrange multipliers computed above, we can evaluate the gradient of the cost functional \eqref{eq:CostApprox} as 
	
	\allowdisplaybreaks
	\begin{align}\label{eq:gradient_z}
	%\begin{split} %\label{eq:gradient_z}
	& \langle \tilde{z}, \nabla_z \cE(z) \rangle = \langle \tilde{z}, \partial_z \cL_2 \rangle 
	\nonumber\\
	& = \langle \tilde{z}, \pc{\partial_z \bar{q} + }\nabla_z \cP(z) + \pc{d^f \partial_z \bar{f} + } \partial_{zm} \bar{r}_f m^f \pc{+ \partial_{zm} \bar{f} m^f}\rangle
	\nonumber\\
	& + \langle \tilde{z}, 
	\partial_{zv} \bar{r} v^*
	\rangle 
	\nonumber\\
	& + \langle \tilde{z}, 
	\partial_{zu} \bar{r}_q (u^q)^* \pc{+ \partial_{zu} \bar{q} (u^q)^*}
	\rangle 
	\nonumber\\
	& + \sum_{n=1}^{N_q} \langle \tilde{z},
	\partial_{zmv} \bar{r}_q \hat{v}^q_n (\psi_n^q)^* + \partial_{zmu} \bar{r}_q \hat{u}^q_n (\psi_n^q)^* + \partial_{zmm} \bar{r}_q \psi_n^q (\psi_n^q)^* \pc{+ \partial_{zmu} \bar{q} \hat{u}^q_n (\psi_n^q)^* + \partial_{zmm} \bar{q} \psi_n^q (\psi_n^q)^*}
	\rangle 
	\nonumber\\
	& +  \sum_{n=1}^{N_q} \langle \tilde{z},
	\partial_{zvu} \bar{r}_q \hat{u}_n^q (\hat{v}_n^q)^* + \partial_{zvm} \bar{r}_q \psi_n^q (\hat{v}_n^q)^*
	\rangle
	\nonumber\\
	& +  \sum_{n=1}^{N_q} \langle \tilde{z},
	\partial_{zuv} \bar{r}_q \hat{v}_n^q (\hat{u}_n^q)^* + \partial_{zuu} \bar{r}_q \hat{u}_n^q (\hat{u}_n^q)^* + \partial_{zum} \bar{r}_q \psi_n^q (\hat{u}_n^q)^* \pc{+ \partial_{zuu} \bar{q} \hat{u}_n^q (\hat{u}_n^q)^* + \partial_{zum} \bar{q} \psi_n^q (\hat{u}_n^q)^*}
	\rangle
	\nonumber\\
	& + \langle \tilde{z},
	\partial_{zu} \bar{r}_f (u^f)^* \pc{+ \partial_{zu} \bar{f} (u^f)^* }
	\rangle
	\nonumber\\
	& +  \sum_{n=1}^{N_f} \langle \tilde{z},
	\partial_{zmv} \bar{r}_f \hat{v}^f_n (\psi_n^f)^* + \partial_{zmu} \bar{r}_f \hat{u}^f_n (\psi_n^f)^* + \partial_{zmm} \bar{r}_f \psi_n^f (\psi_n^f)^* \pc{+ \partial_{zmu} \bar{f} \hat{u}^f_n (\psi_n^f)^* + \partial_{zmm} \bar{f} \psi_n^f (\psi_n^f)^*}
	\rangle 
	\nonumber\\
	& +  \sum_{n=1}^{N_f} \langle \tilde{z},
	\partial_{zvu} \bar{r}_f \hat{u}_n^f (\hat{v}_n^f)^* + \partial_{zvm} \bar{r}_f \psi_n^f (\hat{v}_n^f)^*
	\rangle
	\nonumber\\
	& +  \sum_{n=1}^{N_f} \langle \tilde{z},
	\partial_{zuv} \bar{r}_f \hat{v}_n^f (\hat{u}_n^f)^* + \partial_{zuu} \bar{r}_f \hat{u}_n^f (\hat{u}_n^f)^* + \partial_{zum} \bar{r}_f \psi_n^f (\hat{u}_n^f)^*\pc{+ \partial_{zuu} \bar{f} \hat{u}_n^f (\hat{u}_n^f)^* + \partial_{zum} \bar{f} \psi_n^f (\hat{u}_n^f)^*}
	\rangle.
	\nonumber
	%\end{split}
	\end{align}
	
	% \pc{the constants $c_n^f$ and $m^f$ also depend on $z$} 
	
	We remark that the constraints of the orthonormal conditions in the Lagrangian $\cL_2$ do not explicitly depend on the optimization variable $z$, so that for computing the gradient $\nabla_z \cE(z)$ we do not need the Lagrange multipliers $\{(\lambda_{nn'}^q)^*\}$ and $\{(\lambda_{nn'}^f)^*\}$, neither of which are used in computing all other Lagrange multipliers. To solve the chance constrained optimization problem under uncertainty with the unconstrained penalty formulation \eqref{eq:OptUnC} and its approximation in \eqref{eq:CostApprox}, we apply a gradient-based BFGS algorithm, where the evaluation of the cost functional $\cE(z)$ and its gradient $\nabla_k \cE(z)$ at $z$, and their computational cost in terms of PDE solves, are summarized in \eqref{alg:BFGS4CCOUU}. In summary, one evaluation of the cost functional takes $1$ state PDE solve and $2(N_q+c) + 2(N_f+c) + 2$ linearized PDE solves, while one evaluation of its gradient takes $2N_q + 2N_f + 3$ linearized PDE solves, where $N_q$ and $N_f$ are the ranks in \eqref{eq:TraceApprox} and \eqref{eq:HessianApprox}, $c$ is a small oversampling parameter in Algorithm \ref{alg:RandomizedSVD}.
	
	\renewcommand{\algorithmiccomment}[1]{\bgroup \hfill \footnotesize//~#1\egroup}
	
	\begin{algorithm}
		\caption{{Compute} $\cE(z)$ and its gradient $\nabla_z \cE(z)$}
		\label{alg:BFGS4CCOUU}
		\begin{algorithmic}[1]
			\STATE{\textbf{Compute} $\cE(z)$:}\label{lst:CostApprox}
			\STATE{solve the state equation \eqref{eq:state} for $u$;}  \algorithmiccomment{1 state PDE} 
			\STATE{solve the linearized PDE \eqref{eq:adjoint_q} for $v^q$;} \algorithmiccomment{1 linear PDE}
			\STATE{solve the linearized PDE \eqref{eq:adjoint_f} for $v^f$;} \algorithmiccomment{1 linear PDE}
			\STATE{solve the generalized eigenvalue problems \eqref{eq:GenEigen_q}  for $(\lambda_n^q, \psi_n^q)_{n=1}^{N_q}$;} \algorithmiccomment{$2(N_q+c)$ linear PDEs} 
			\STATE{solve the generalized eigenvalue problems \eqref{eq:GenEigen_f} for $(\lambda_n^f, \psi_n^f)_{n=1}^{N_f}$;} \algorithmiccomment{$2(N_f+c)$ linear PDEs}
			\STATE{compute the cost functional $\cE(z_k)$ by \eqref{eq:CostApprox}.}
			\STATE{\textbf{Compute} $\nabla_z \cE(z)$: }\label{lst:gradient}
			\STATE{solve the linearized PDEs \eqref{eq:incfwd_q} and \eqref{eq:incadj_q} at $(\psi_n^q)_{n=1}^{N_q}$ for $(\hat{u}_n^q, \hat{v}_n^q)_{n=1}^{N_q}$;} \algorithmiccomment{$2N_q$ linear PDEs}
			% at $(\psi_n^q)_{n=1}^{N_q}$ 
			\STATE{solve the linearized PDEs \eqref{eq:incfwd_f} and \eqref{eq:incadj_f} at $(\psi_n^f)_{n=1}^{N_f}$ for $(\hat{u}_n^f, \hat{v}_n^f)_{n=1}^{N_f}$;} \algorithmiccomment{$2N_f$ linear PDEs}
			\STATE{set $((\psi_n^q)^*, (\hat{u}_n^q)^*, (\hat{v}_n^q)^*)_{n=1}^{N_q}$ by \eqref{eq:PhiStar_q}, \eqref{eq:IncStateStar_q}, and \eqref{eq:IncAdjointStar_q};}
			\STATE{solve the linearized PDE  \eqref{eq:uStar_q} for  $(u^q)^*$;} \algorithmiccomment{1 linear PDE}
			\STATE{set $((\psi_n^f)^*, (\hat{u}_n^f)^*, (\hat{v}_n^f)^*)_{n=1}^{N_f}$ by \eqref{eq:PhiStar_f}, \eqref{eq:IncStateStar_f}, and \eqref{eq:IncAdjointStar_f};}
			\STATE{solve the linearized PDE  \eqref{eq:uStar_f}for  $(u^f)^*$;} \algorithmiccomment{1 linear PDE}
			\STATE{solve the linearized PDE \eqref{eq:vStar} for $v^*$;} \algorithmiccomment{1 linear PDE}
			\STATE{compute the gradient $\nabla_z \cE(z_k)$ by \eqref{eq:gradient_z}.}
		\end{algorithmic}
	\end{algorithm}

	\section{Numerical examples}
	\label{sec:numerics}
	
\begin{figure}[!htb]
		\centering
		\includegraphics[width=0.48\textwidth]{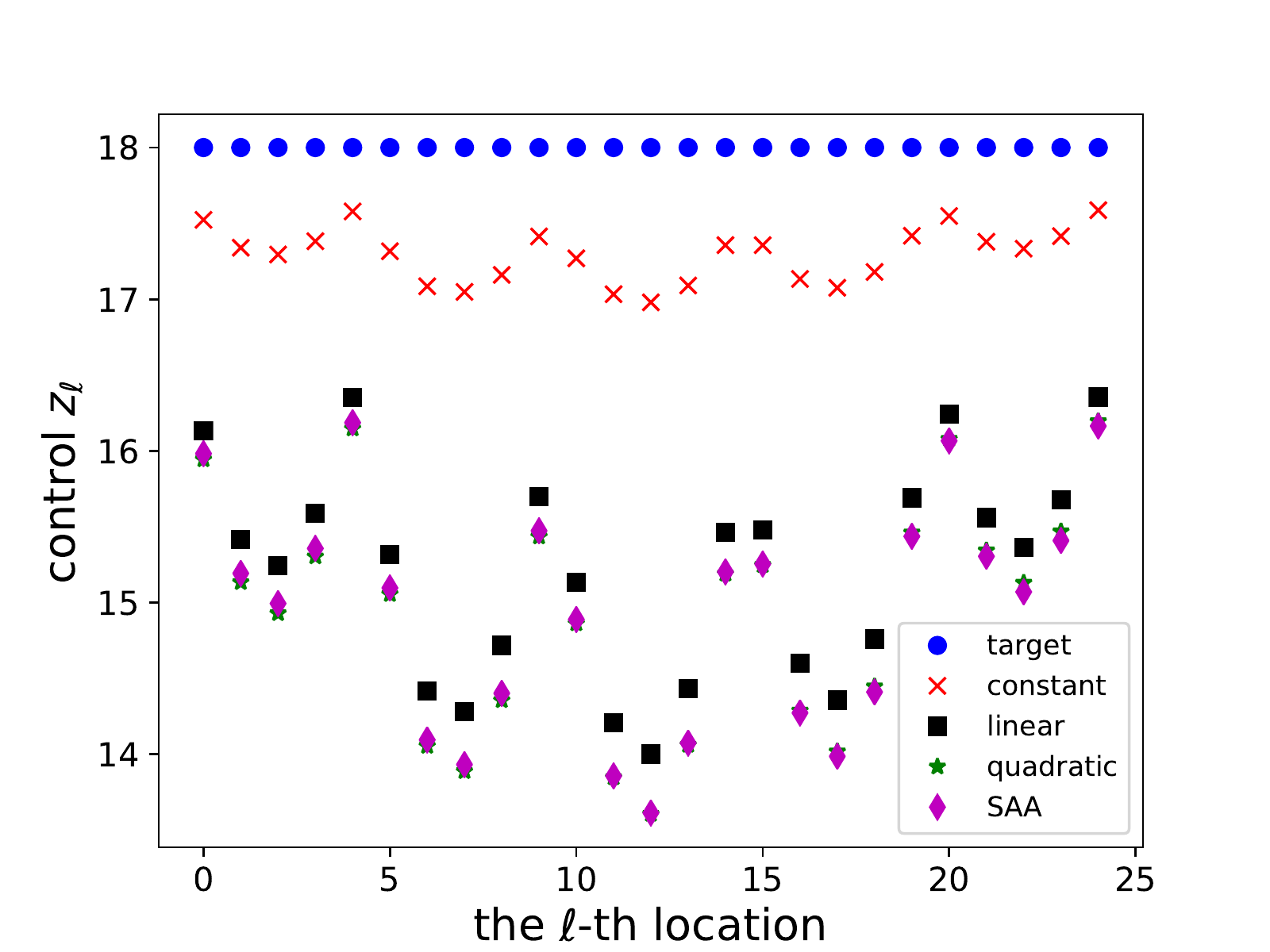}
		\includegraphics[width=0.48\textwidth]{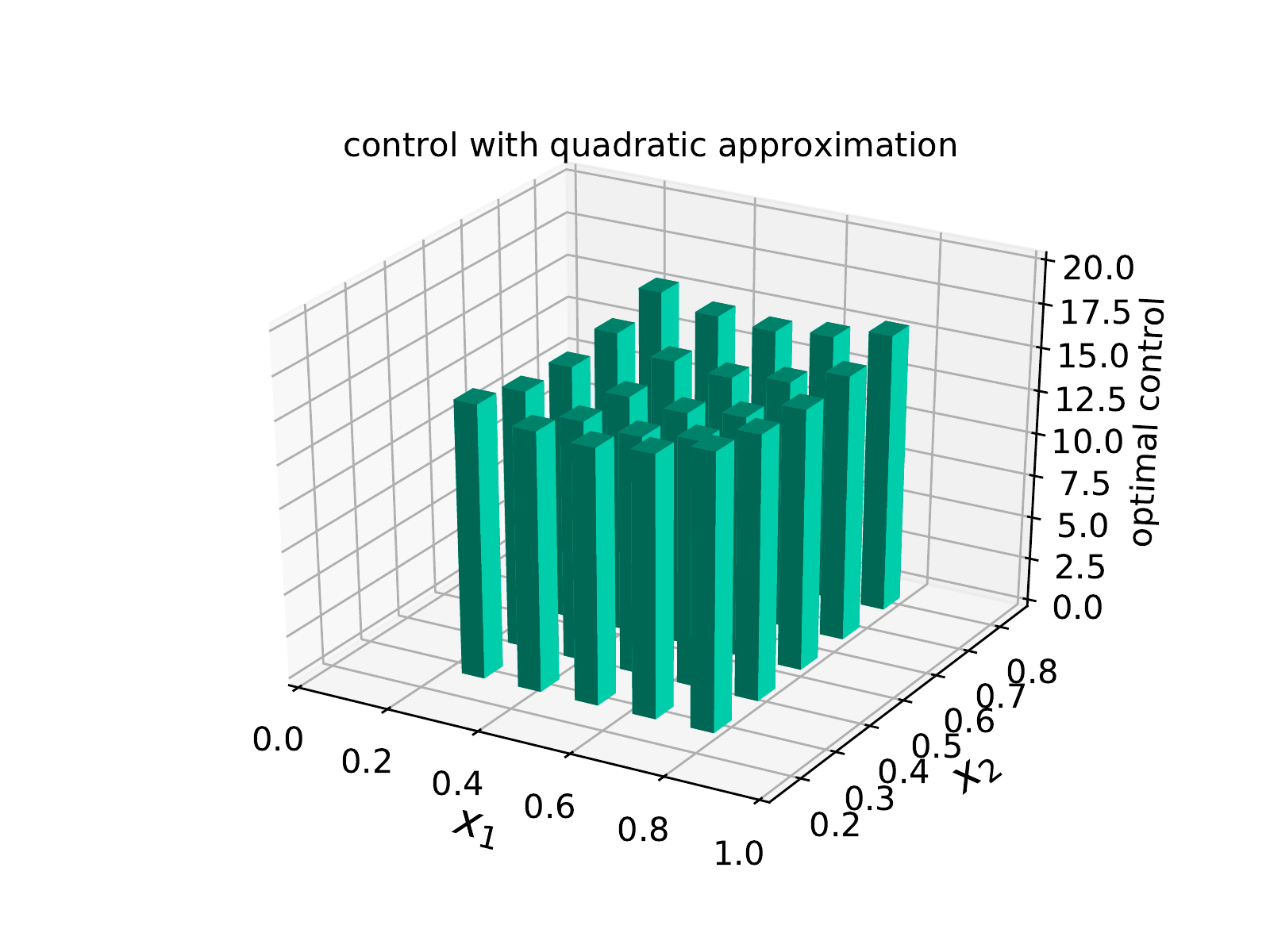}
		\caption{Left: comparison of optimal variables $z_\ell$ at the $\ell$-th location obtained by different approximations of the constraint function $f$,  including Taylor approximations $T_K f$ for $K = 0$ (constant), $K = 1$ (linear), $K = 2$ (quadratic), and sample average approximation (SAA) with $1024$ samples. Right: distribution of the optimal variable $z^*$ obtained with quadratic approximation of the constraint function.
		}\label{fig:optimal-controls}
	\end{figure}

	We consider the following PDEs that model a steady state Darcy flow,
	\begin{equation}\label{eq:darcy}
	\begin{split}
	\mathbf{v} + \frac{e^m}{\mu} \nabla u  & = 0\quad \text{ in } D, \\
	\nabla \cdot \mathbf{v}  & = h \quad \text{ in } D, 
	\end{split}
	\end{equation}
	where a homogeneous Dirichlet boundary condition for the pressure $u$ is imposed along the boundary $\partial D$ of a physical domain $D = (0, 1)^2$.  $e^m$ represents a random permeability field, while $\mu$ is the fluid viscosity. For simplicity we specify $\mu = 1$ in a dimensionless setting. The source term $h$ depends on an $L$-dimensional (we take $L = 25$ in the numerical test) optimization variable $z = (z_1, \dots, z_L)$, given by 
	\begin{equation}
	h(z) = -\sum_{\ell = 1}^L z_\ell h_\ell, 
	\end{equation}
	where $z_\ell$ is a pointwise optimization variable with bound $z_\ell \in [z_{\text{min}}, z_{\text{max}}]$, where we take $z_{\text{min}} = 0$ and $z_{\text{max}} = 36$; $h_\ell$ is a smooth mollifier function defined at point $x_\ell \in D$ as 
	\begin{equation}
	h_\ell = \exp\left(-\frac{1}{\varepsilon^2}||x-x_\ell||^2\right)
	\end{equation}
	for a positive number $\varepsilon > 0$, which we take $\varepsilon = 0.1$.
	\begin{figure}[!htb]
		\centering
		\includegraphics[trim=20 20 40 40, clip, width=0.5\textwidth]{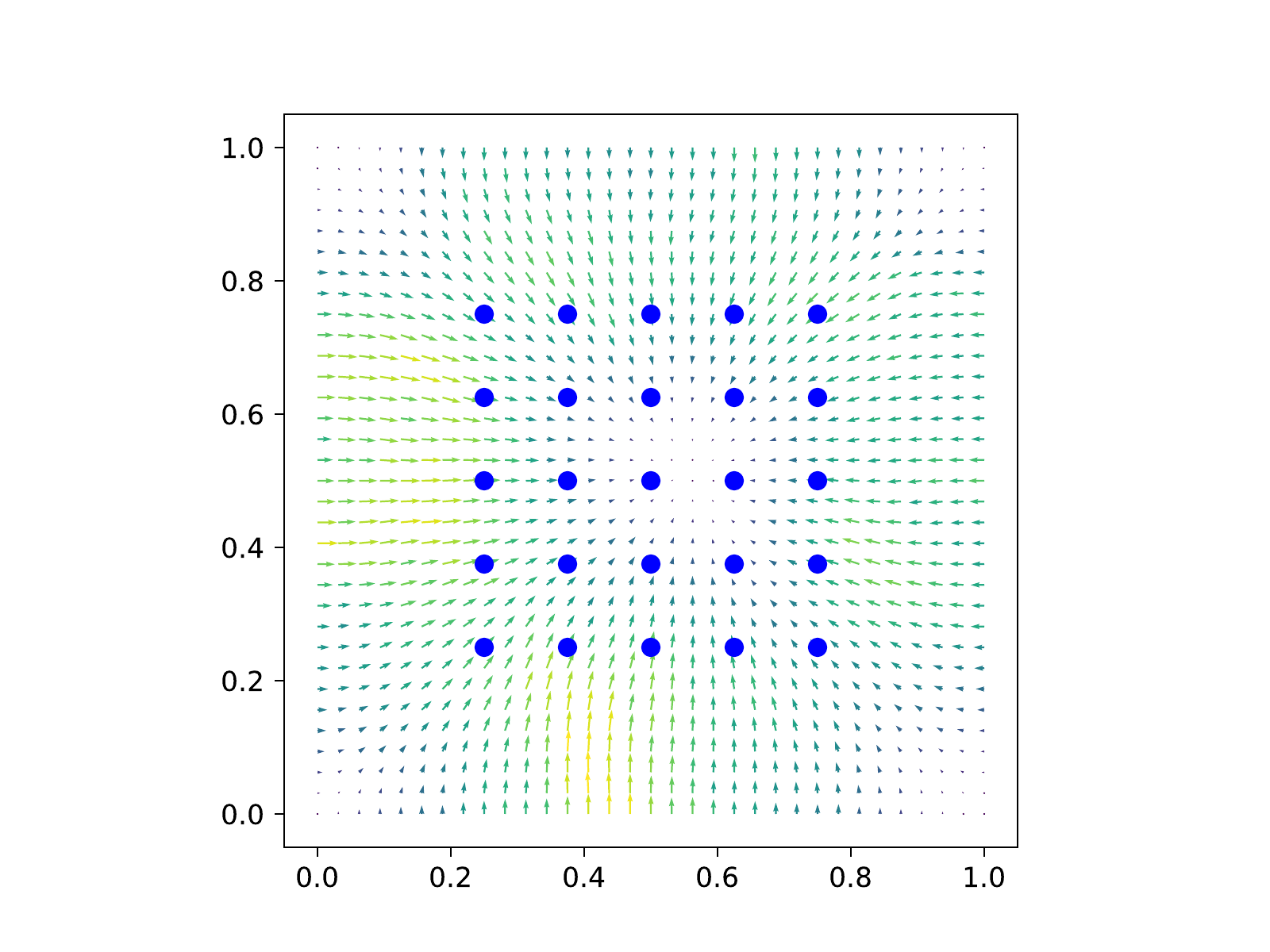}
		\caption{Groundwater flow in physical domain $(0, 1)^2$. Blue dots stand for the location of the extraction wells. A velocity field $\mathbf{v}$ is shown at the mean of the log-permeability $m = \bar{m}$ and optimal variable $z^*$, which is obtained with quadratic approximation of the constraint function; see Figure \ref{fig:optimal-controls} for the value of $z^*$.}\label{fig:flow}
	\end{figure}
	The system \eqref{eq:lognormal-pde} models steady state groundwater flow, $z_\ell$ represents the water extraction rate at location $x_\ell$ of $L$ wells. Figure \ref{fig:flow} illustrate the groundwater flow, where the velocity field $\mathbf{v}$ at the mean $\bar{m}$ and an optimal extraction rate (with value given in Figure \ref{fig:optimal-controls}) is shown.
	Note that by eliminating the velocity field $\mathbf{v}$ from \eqref{eq:darcy}, we obtain a single equation for the pressure 
	\begin{equation}\label{eq:lognormal-pde}
	-\nabla \cdot \left(\frac{e^m}{\mu} \nabla u\right) = h \quad \text{ in } D.
	\end{equation}
	
	The objective of the optimization problem is to achieve a target groundwater extraction rate $\bar{z}_\ell$ at each well, which can be represented by 
	\begin{equation}
	q(z) = \frac{1}{L} \sum_{\ell = 1}^L (z_\ell - \bar{z}_\ell)^2.
	\end{equation}
	We use a penalty term $\cP(z) = \frac{\eta}{2} ||z||^2_2$ with $\eta = 10^{-5}$ for the optimization variable $z$, representing the cost of the extraction. 
	To prevent excessive extraction leading to potential collapse of the aquifer, we consider the constraint function for the state (pressure field) $u$ 
	\begin{equation}
	f(u) = \int_{D_o} u^2(x) dx - f_c,
	\end{equation}
	where $D_o \subset D$ is a region of interest, for which we take $D_o = (0.25, 0.75)^2$; and $f_c > 0$ is a critical value, which we take $f_c = 2$. We consider the chance or probability constraint \eqref{eq:probability_constraint},  i.e., $P(f \geq 0) \leq \alpha$, such that the probability of $f$ greater than or equal to zero should be less than or equal to a given value $\alpha > 0$ (we take $\alpha = 0.05$), where the probability is defined with respect to the probability distribution of the random field $m$.

	\begin{figure}[!htb]
		\centering
		\includegraphics[trim=80 80 80 80, clip, width=0.32\textwidth]{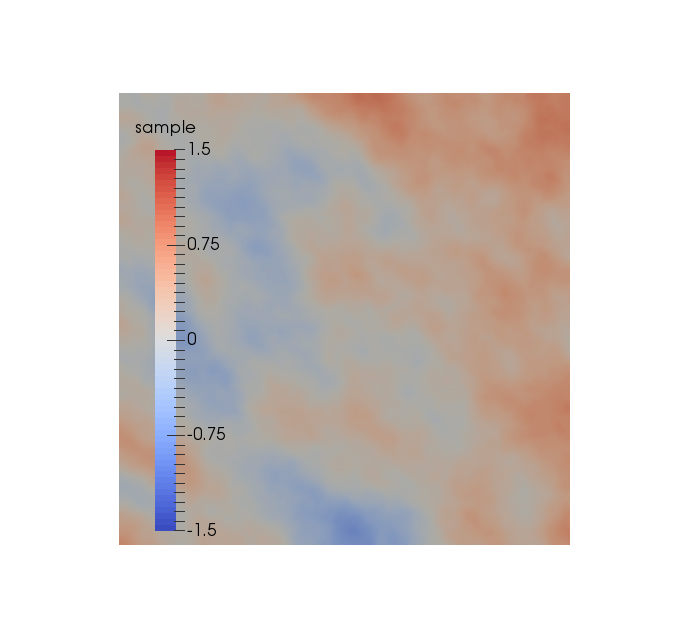}
		\includegraphics[trim=80 80 80 80, clip, width=0.32\textwidth]{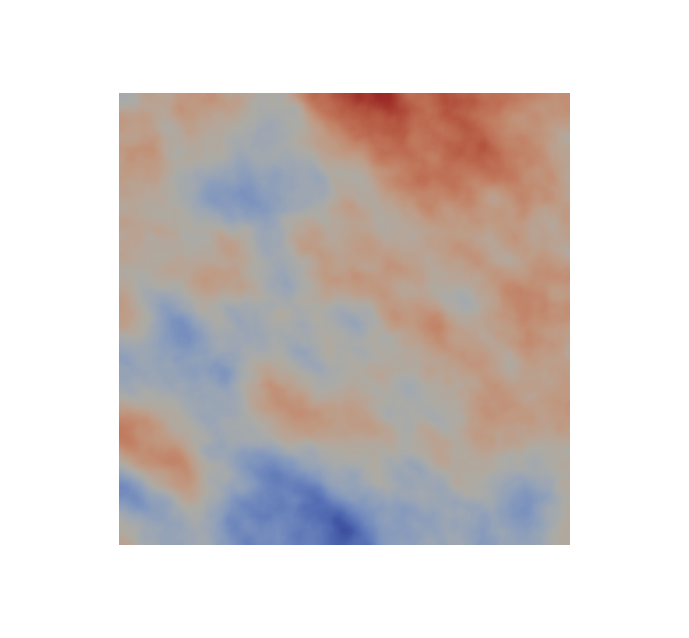}
		\includegraphics[trim=80 80 80 80, clip, width=0.32\textwidth]{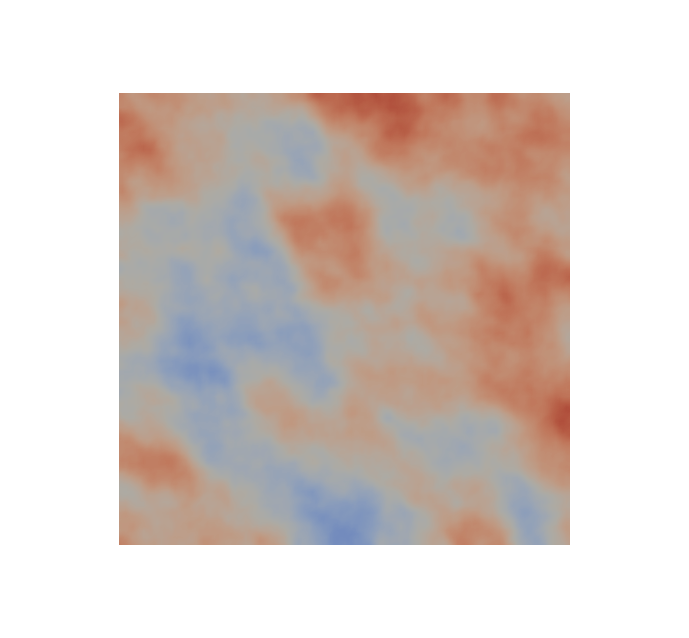}
		\caption{Mean of the log-permeability field $\bar{m}$ (left) and two random samples (middle and right).}\label{fig:sample}
	\end{figure}
	We use a Gaussian random field $m \sim \cN(\bar{m}, \cC)$ with mean $\bar{m}$ and covariance field $\cC$, which is represented as the square of the inverse of an elliptic operator \cite{Stuart10, LindgrenRueLindstroem11},
	\begin{equation}\label{eq:cov}
	\cC = (- \eta \Delta + \delta I )^{-2},
	\end{equation}
	where $I$ is the identity and $\Delta$ is the Laplacian operator. $\delta > 0$ and $\eta > 0$ are two parameters that control the variance and correlation of the random field. Sampling from this Gaussian distribution is equivalent to solving the elliptic equation 
	\begin{equation}
	- \eta \Delta m + \delta m =\dot{W} \quad \text{ in } D
	\end{equation}
	with homogeneous Neumann boundary condition along $\partial D$, where $\dot{W}$ is a spatial white noise with unit pointwise variance. The mean $\bar{m}$ is given in Figure \ref{fig:sample}, which is obtained as a random sample from $\cN(0, \cC)$ with $\eta = 0.1$, $\delta = 10$ for the covariance $\cC$ in \eqref{eq:cov}. Two random samples drawn from $\cN(\bar{m}, \cC)$ are also shown in Figure \ref{fig:sample}.
	
	\begin{figure}[!htb]
		\centering
		
		\includegraphics[width=0.48\textwidth]{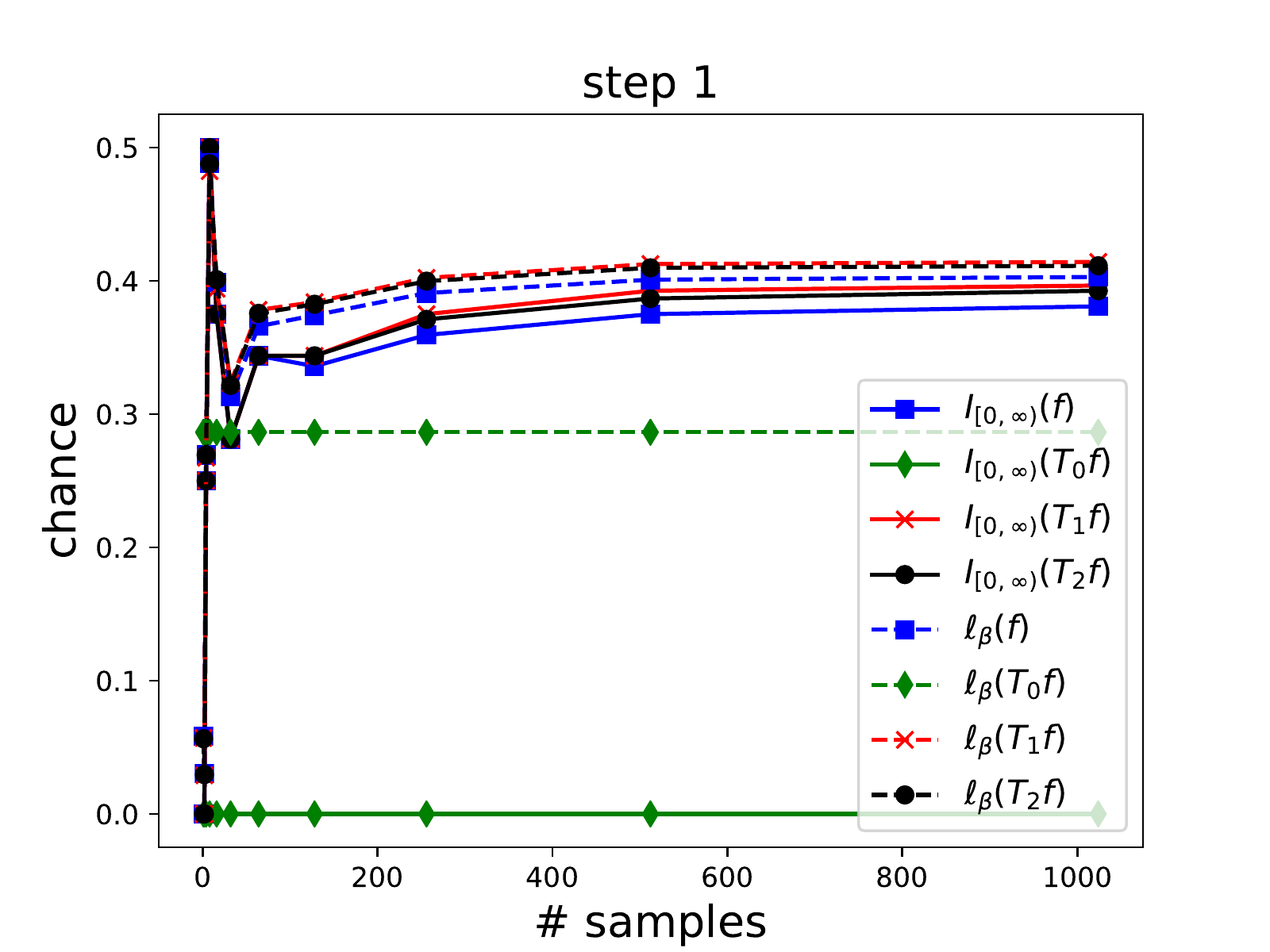}
		\includegraphics[width=0.48\textwidth]{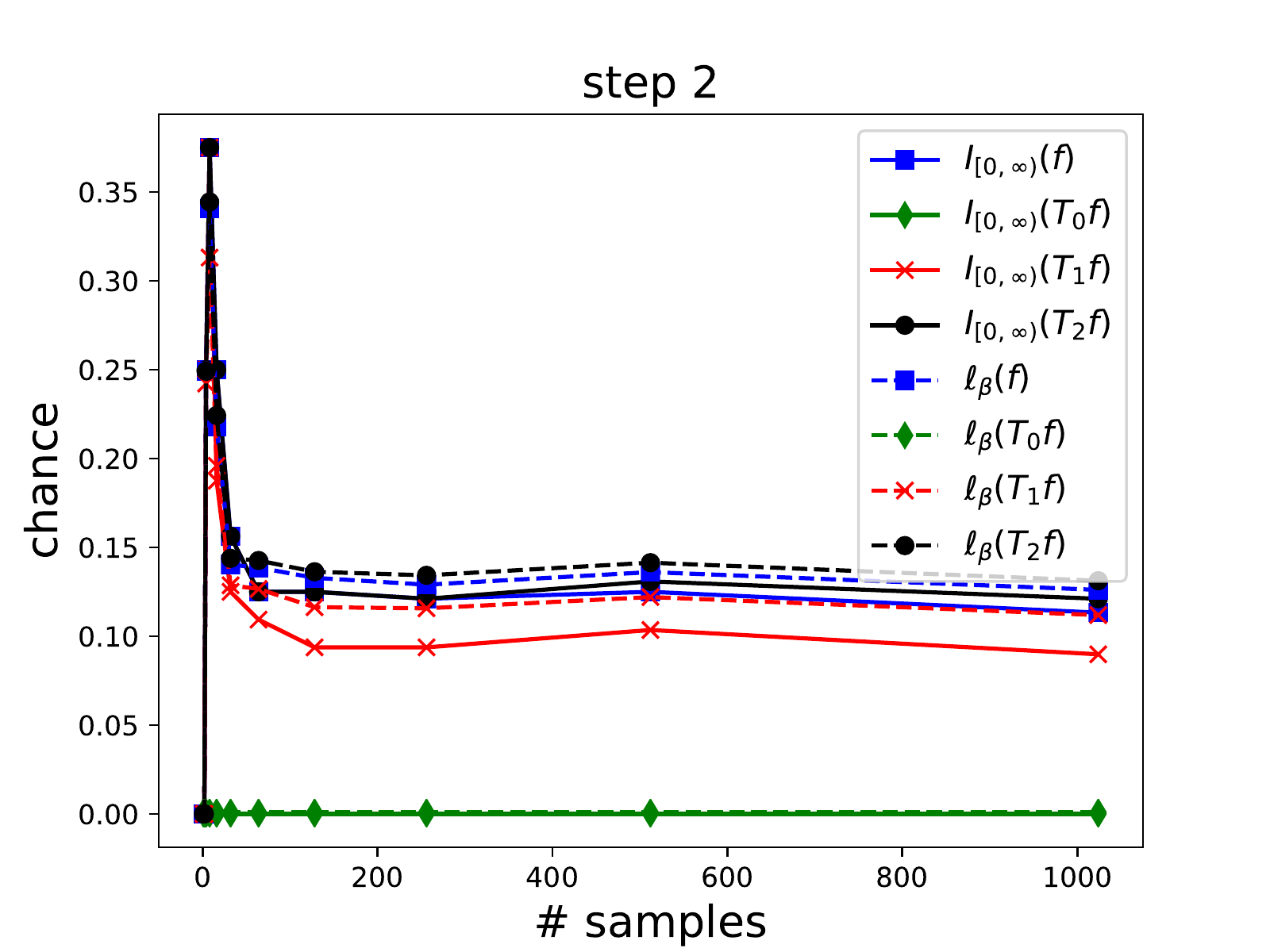}
		
		\includegraphics[width=0.48\textwidth]{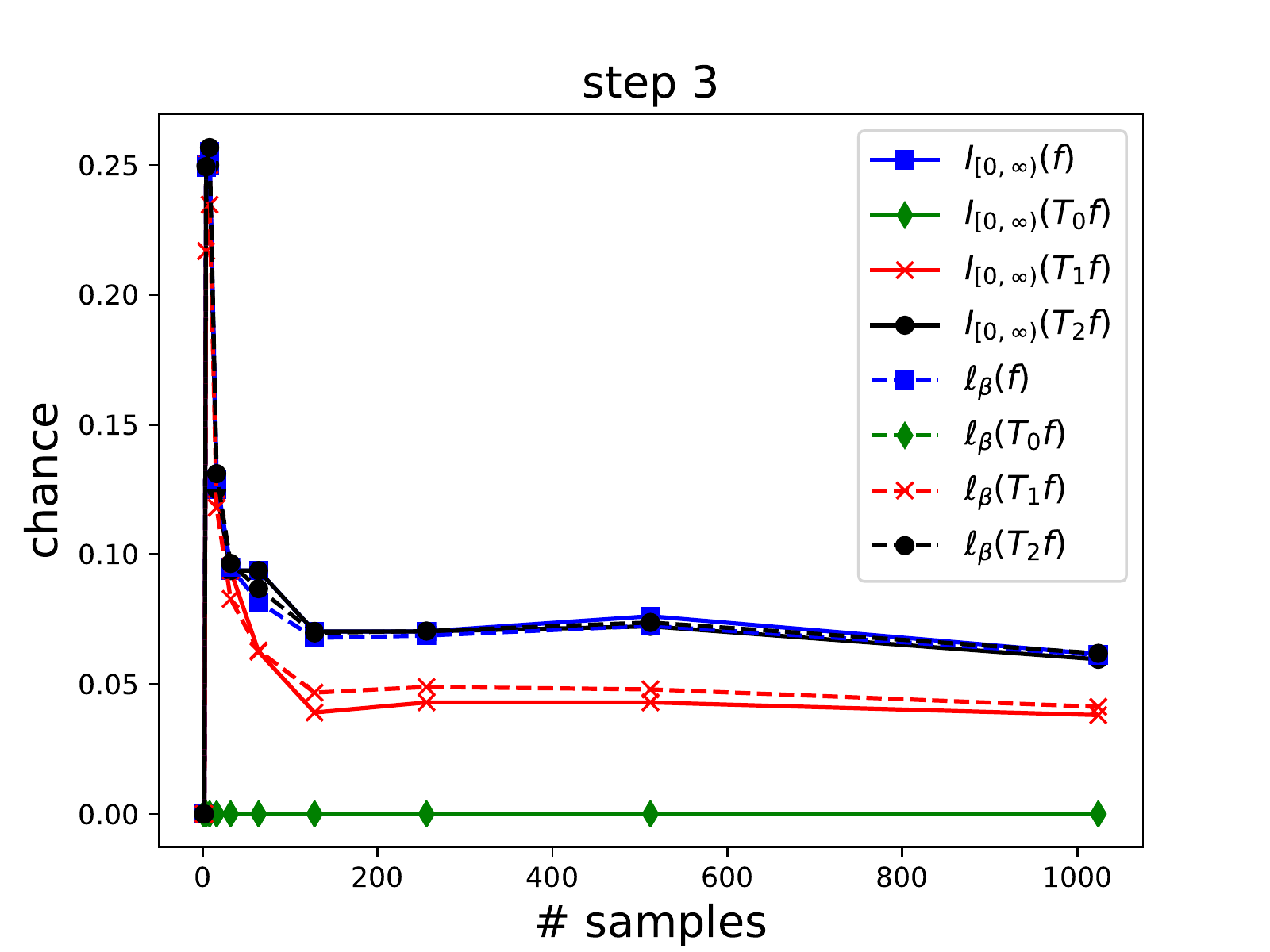}
		\includegraphics[width=0.48\textwidth]{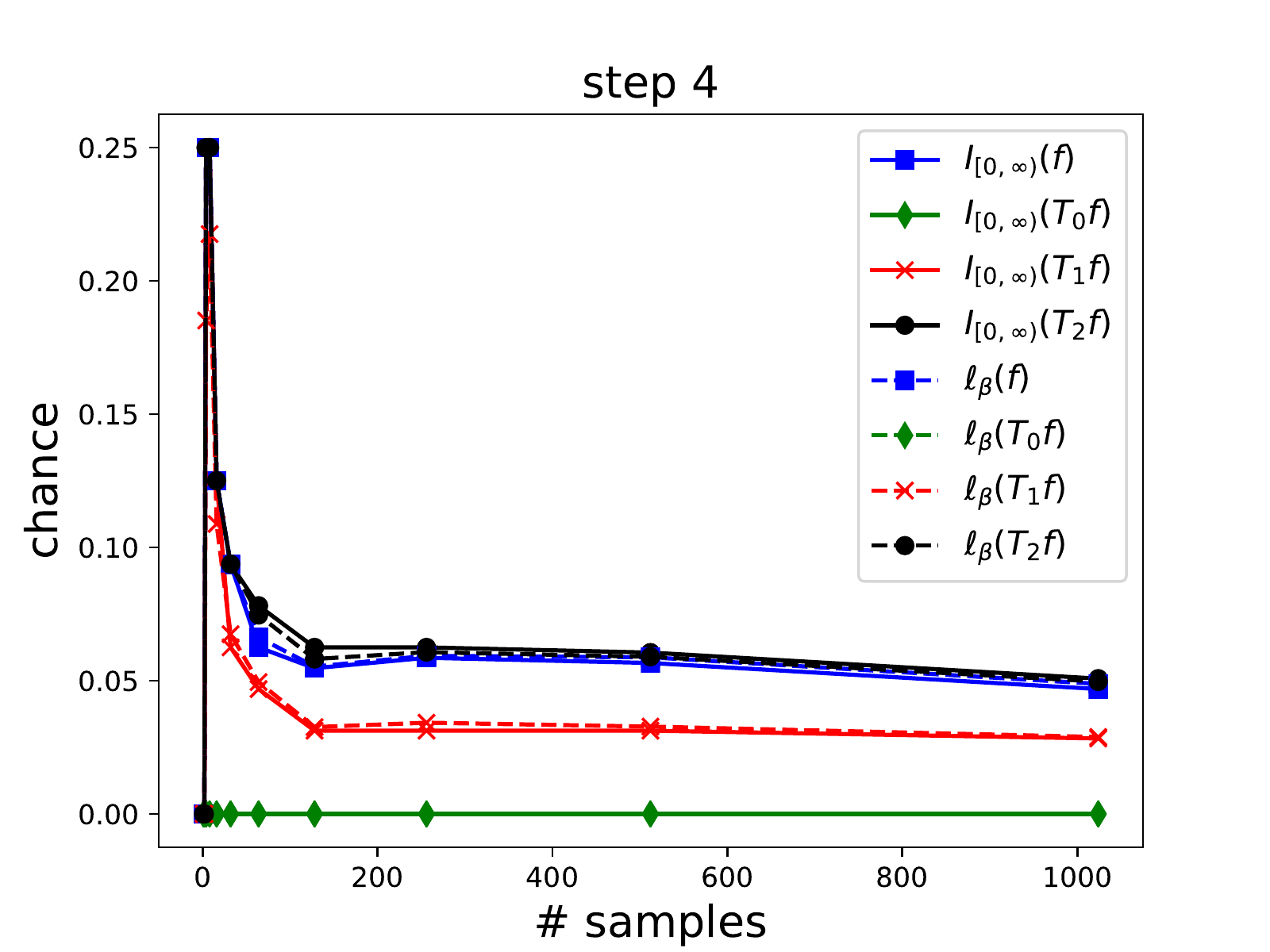}
		\caption{The chance evaluated with different numbers of samples, by different Taylor approximations, and using different smoothing and penalty parameters $(\beta, \gamma) = (2^{n+2}, 10^{n+2})$ for step $n = 1, 2, 3, 4$.
		}\label{fig:chance}
	\end{figure}
	
	We use a finite element method implemented in FEniCS \cite{LoggMardalGarth12} to solve all of the PDEs, with linear elements for the approximation of both the pressure field $u$ and the parameter field $m$  in a uniform mesh of triangles of size $32 \times 32$, which leads to $33^2$ dimension for the discrete parameter.
	We run the optimization algorithm, Algorithm \ref{alg:BFGS4CCOUU}, to solve the chance constrained optimization problem \eqref{eq:optimization}, with Taylor approximation of the constraint function $T_K f$ with $K = 0, 1, 2$, and a sample average approximation (SAA) with $1024$ samples. We set the tolerance for the gradient norm as $10^{-3}$ for the stopping criterion. For the quadratic Taylor approximation, we compute 10 eigenpairs \eqref{eq:GenEigen_f} by the randomized algorithm, Algorithm \ref{alg:RandomizedSVD}.
	The smoothing and penalty parameters are specified as $(\beta, \gamma) = (2^{n+2}, 10^{n+2})$ for four steps $n = 1, 2, 3, 4$. The optimization solution at step $n$ is used as the initial guess for that at step $n+1$. After each optimization step with quadratic approximation $T_2f$, we compute the chance $P(g > 0) =  \bE[\bI_{[0, \infty)}(g)] $ and the smoothed value $\bE[\ell_\beta(g)]$, by SAA approximation in \eqref{eq:ChanceSAA} with different numbers of samples, where $g$ represents the constraint function $f$ or its Taylor approximations $T_K f$ for $K = 0, 1, 2$. The results are shown in Figure \ref{fig:chance}, from which we can observe: (1) the quadratic approximation $T_2 f$ yields more accurate approximation of the chance than that by the linear and constant approximations; (2) as the smoothing parameter $\beta$ increases, the smooth approximation $\ell_\beta(g)$ becomes more accurate and leads to more accurate approximation of the chance for $g = f, T_K f$, $K = 0, 1, 2$; (3) as the penalty parameter $\gamma$ increases, the violation of the chance $P(g > 0)$ exceeding $\alpha = 0.05$ is more strongly penalized, and the chance converges to $0.05$ by the quadratic approximation $T_2 f$, while the chance by linear and constant approximations become smaller than $\alpha = 0.05$, i.e., the violation is over penalized; and (4) with an increasing number of samples for the SAA, the approximate chance becomes more accurate. 
	
	\begin{figure}[!htb]
		\centering
		\includegraphics[width=0.48\textwidth]{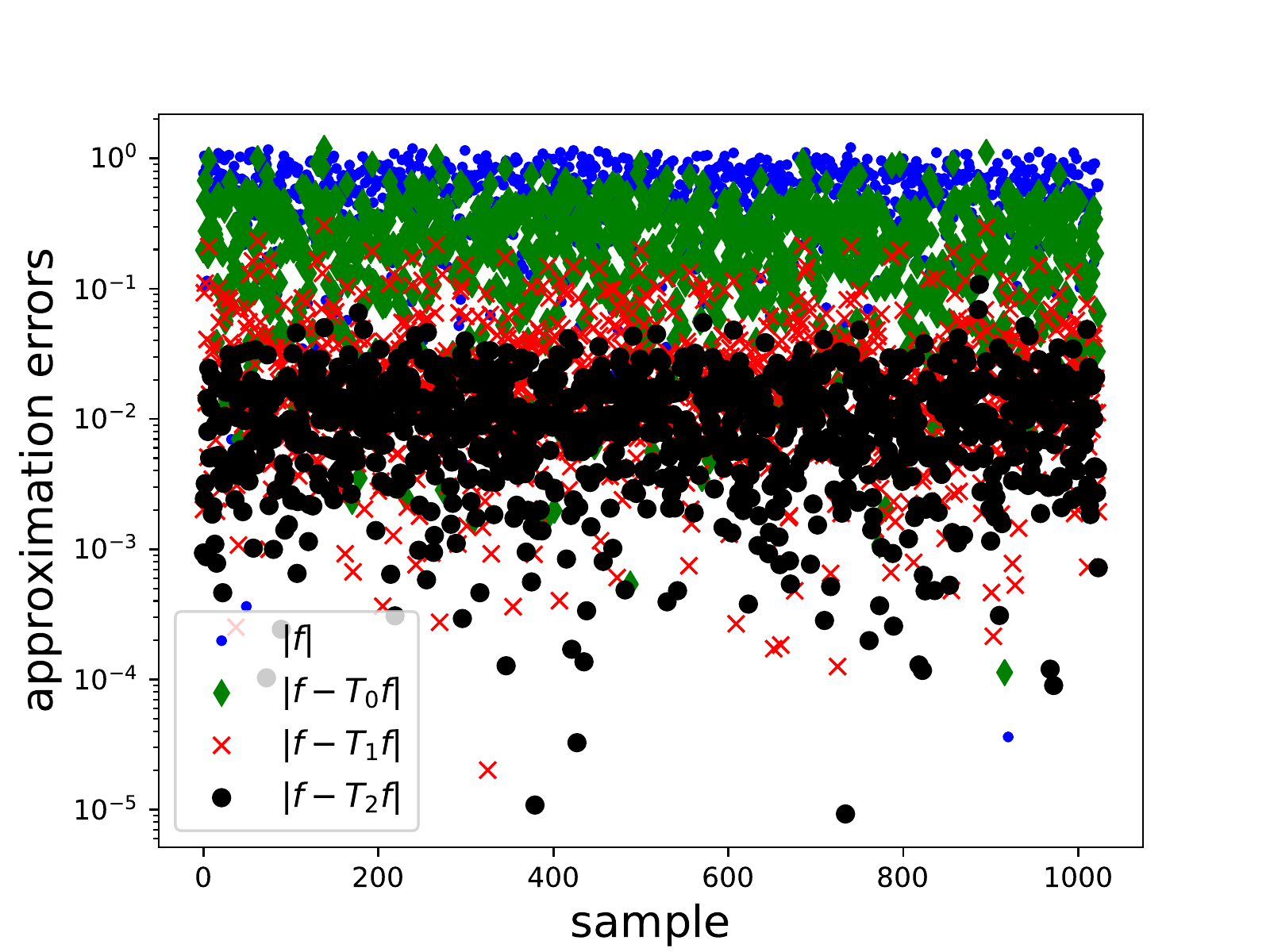}
		\includegraphics[width=0.48\textwidth]{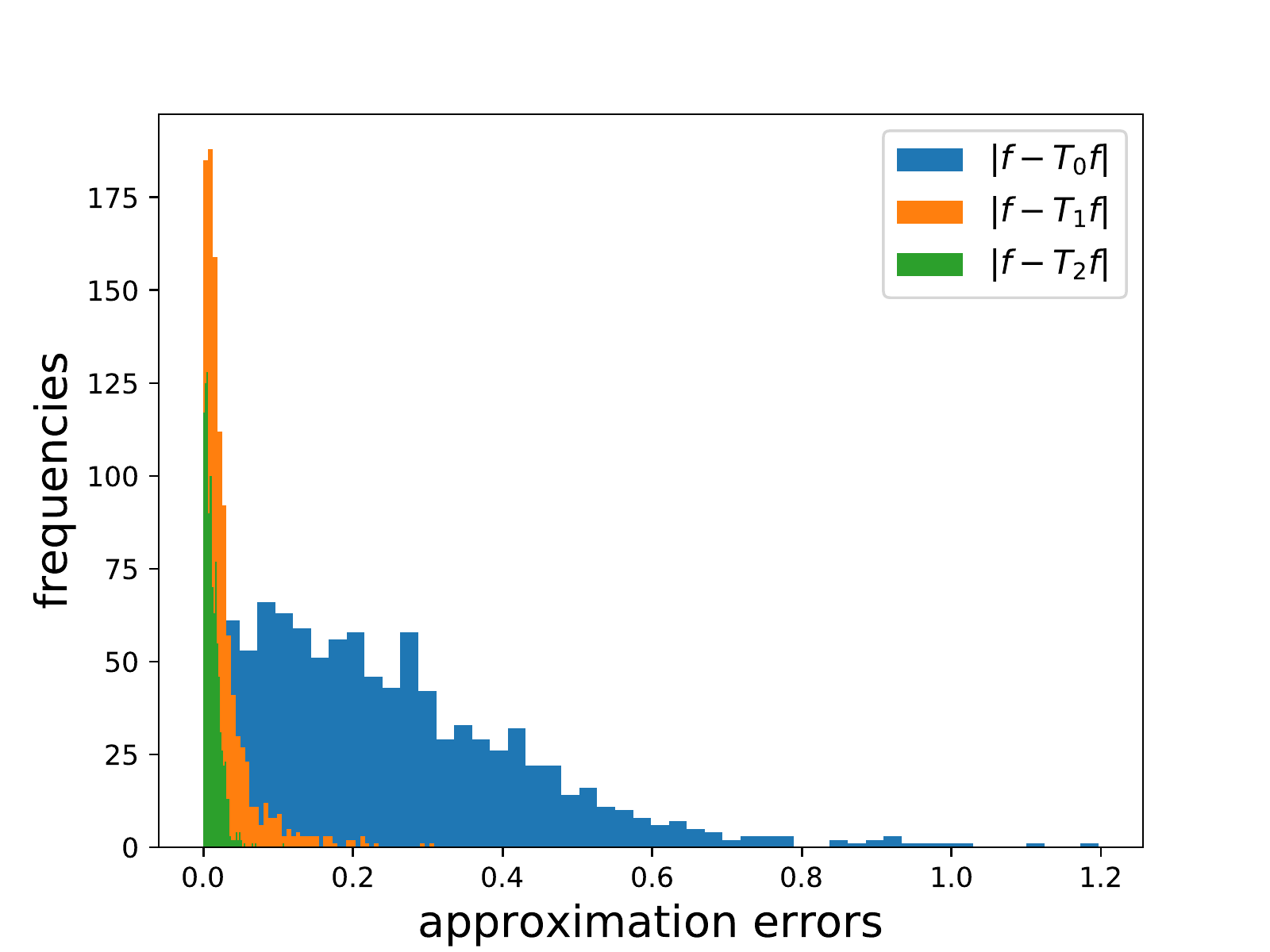}
		\caption{Left: (absolute) values of the constraint function $|f|$ at $1024$ samples and their approximation errors by Taylor approximations $T_K$, $K = 0, 1, 2,$. Right: histograms of the approximation errors $|f - T_K f|$.}\label{fig:Taylor}
	\end{figure}
	
	\begin{table}[!htb]
		\caption{Estimate of the chance $\hat{\ell}_\beta(g)$ for $g = f, T_K f$, $K = 0, 1, 2$, by SAA with 1024 samples, and the estimation bias ($\sqrt{\text{MSE}/1024}$) induced by 1024 samples for SAA  and estimation errors by Taylor approximations.}\label{table:errors}
		\begin{tabular}{|c|c|c|c|c|c|}
			\hline
			step & $\hat{\ell}_\beta(f)$ & SAA bias & $|\hat{\ell}_\beta(f) - \hat{\ell}_\beta(T_0f)|$ &  $|\hat{\ell}_\beta(f) - \hat{\ell}_\beta(T_1f)|$ &  $|\hat{\ell}_\beta(f) - \hat{\ell}_\beta(T_2f)|$\\
			\hline
			0   &  8.25E-1   &  8.47E-3   &  1.52E-1   &  9.81E-3  & 6.14E-3  \\
			\hline
			1   &  3.08E-1   &   1.08E-2  &   1.75E-1  &   1.28E-3 &  6.24E-3\\
			\hline
			2   &  1.14E-1  &    8.31E-3  &    1.14E-1  &   2.04E-2 & 2.91E-3 \\	
			\hline
			3  &   6.02E-2  &   6.77E-3  &   6.02E-2  &  2.26E-2  & 7.61E-4 \\
			\hline 
			4  &   5.00E-2  &   6.47E-3 &    4.98E-2 &   2.21E-2 &  8.60E-5\\
			\hline 
		\end{tabular}
	\end{table}
	
	The observation on the accuracy of the Taylor approximations drawn from Figure \ref{fig:chance} is further demonstrated in Figure \ref{fig:Taylor}, where the approximation errors for $1024$ samples (at the optimal variable by the quadratic approximation) are shown on the left and their histograms are shown on the right, from which we can observe the quadratic approximation is statistically more accurate than the linear and constant approximations. Moreover, in Table \ref{table:errors}, we report the SAA of the chance $\hat{\ell_\beta}(g)$ for $g = f, T_K f$, $K = 0, 1, 2$ with the same 1024 samples at each optimal variable $z^*_n$ obtained by the quadratic approximation with different parameters $(\beta, \gamma) = (2^{n+2}, 10^{n+2})$ at step $n = 1, 2, 3, 4$. We take the initial guess of the optimal variable as $z^*_0 = (18, \dots, 18)$ at step 0. The SAA bias for the estimate of the chance $\hat{\ell}_\beta(f)$ is computed as $\sqrt{\text{MSE}/1024}$, i.e., the square root of the mean square error of $\ell_\beta(f)$ divided by the number of samples, $1024$, which is the bias from the true value induced by a finite number of samples. We can see that the quadratic approximation gives a two orders of magnitude smaller estimation error for the chance than the SAA bias with 1024 samples, while the linear and constant approximations lead to larger estimation errors. In each optimization step, 1024 state PDEs and 1024 adjoint PDEs have to be solved for the direct SAA $\hat{\ell}_\beta(f)$, while 1 state PDE, and $55$ linearized PDEs (see the counts in Algorithm \ref{alg:BFGS4CCOUU}) are solved by the quadratic approximation. A speedup factor of $37 \approx 2048/56$ in terms of PDE solves is achieved. A higher speedup factor is achieved when (1) the number of samples is increased; and (2) the state PDE is more expensive to solve than the linearized PDEs, as is the case for nonlinear state PDEs.
	
\begin{figure}[!htb]
		\centering
		\includegraphics[trim=80 80 80 80, clip, width=0.4\textwidth]{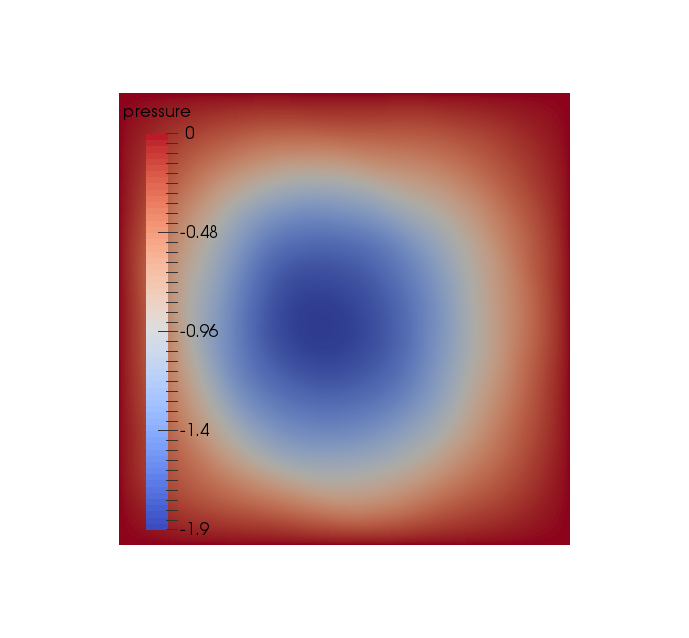}
		\includegraphics[trim=80 80 80 80, clip, width=0.4\textwidth]{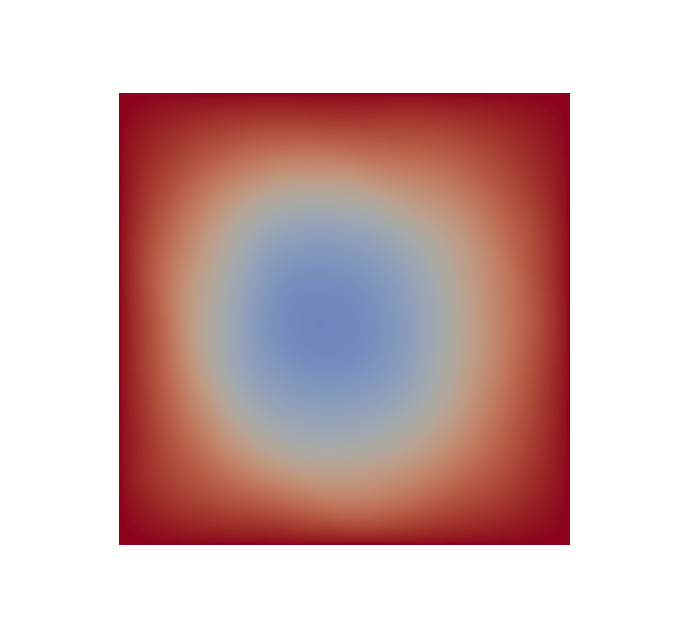}
		\caption{Pressure field at the mean of the log-permeability $m = \bar{m}$ and the target $\bar{z}$ (left), the optimal variable $z^*$  obtained with quadratic approximation (right); see Figure \ref{fig:optimal-controls} for the value of $z^*$. The pressure (in absolute value) at the optimal variable $z^*$ is effectively reduced from that at the target $\bar{z}$.}\label{fig:pressure}
	\end{figure}

	The target $\bar{z} = (18, \dots, 18)$, and the optimal variables obtained by different approximations are shown in Figure \ref{fig:optimal-controls}, from which we observe that the optimal variable obtained by the quadratic approximation of the constraint function, i.e., using SAA for $\ell_\beta(T_2 f)$, is very close to that by SAA for $\ell_\beta(f)$. The distribution of the optimal variable by the quadratic approximation is shown in the right part of Figure \ref{fig:optimal-controls}, with the corresponding pressure field shown in the right part of Figure \ref{fig:pressure}, whose (absolute) value in the region of interest $D_o = (0.25, 0.75)^2$ is effectively reduced from the initial state as displayed in the left part of Figure \ref{fig:pressure}.
	
	%\begin{figure}
	%	\centering
	%	
	%	\includegraphics[width=0.48\textwidth]{figure/Beta_1}
	%	\includegraphics[width=0.48\textwidth]{figure/Beta_4}
	%	
	%%	\includegraphics[width=0.48\textwidth]{figure/Chance_3}
	%%	\includegraphics[width=0.48\textwidth]{figure/Chance_4}
	%	\caption{Convergence with respect to the smoothing parameter $\beta$ at step $n = 1, 4$ corresponding to Figure \ref{fig:chance}.}
	%\end{figure}
	
	\begin{figure}[!htb]
		\centering
		\includegraphics[width=0.45\textwidth]{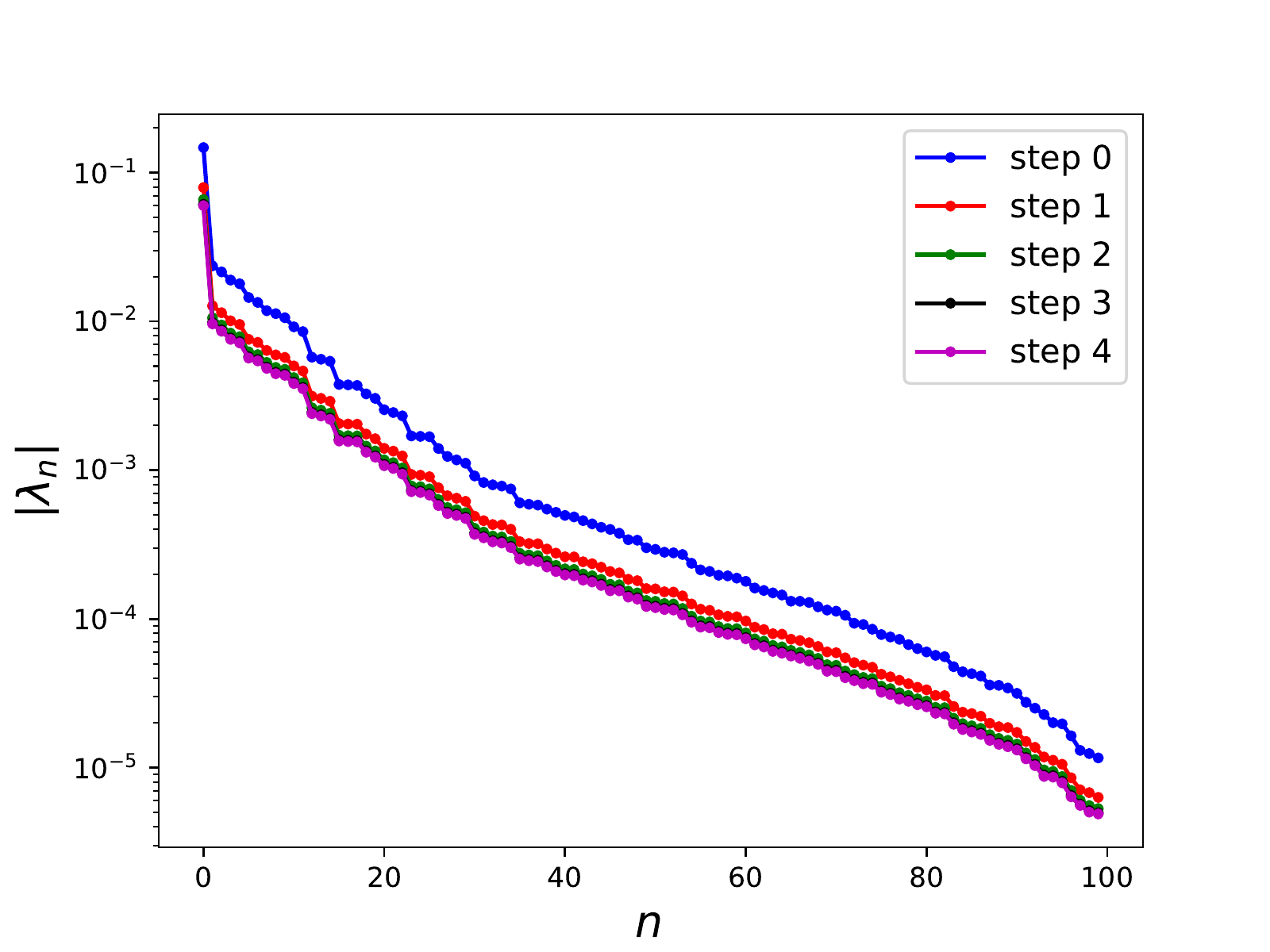}
		\includegraphics[width=0.45\textwidth]{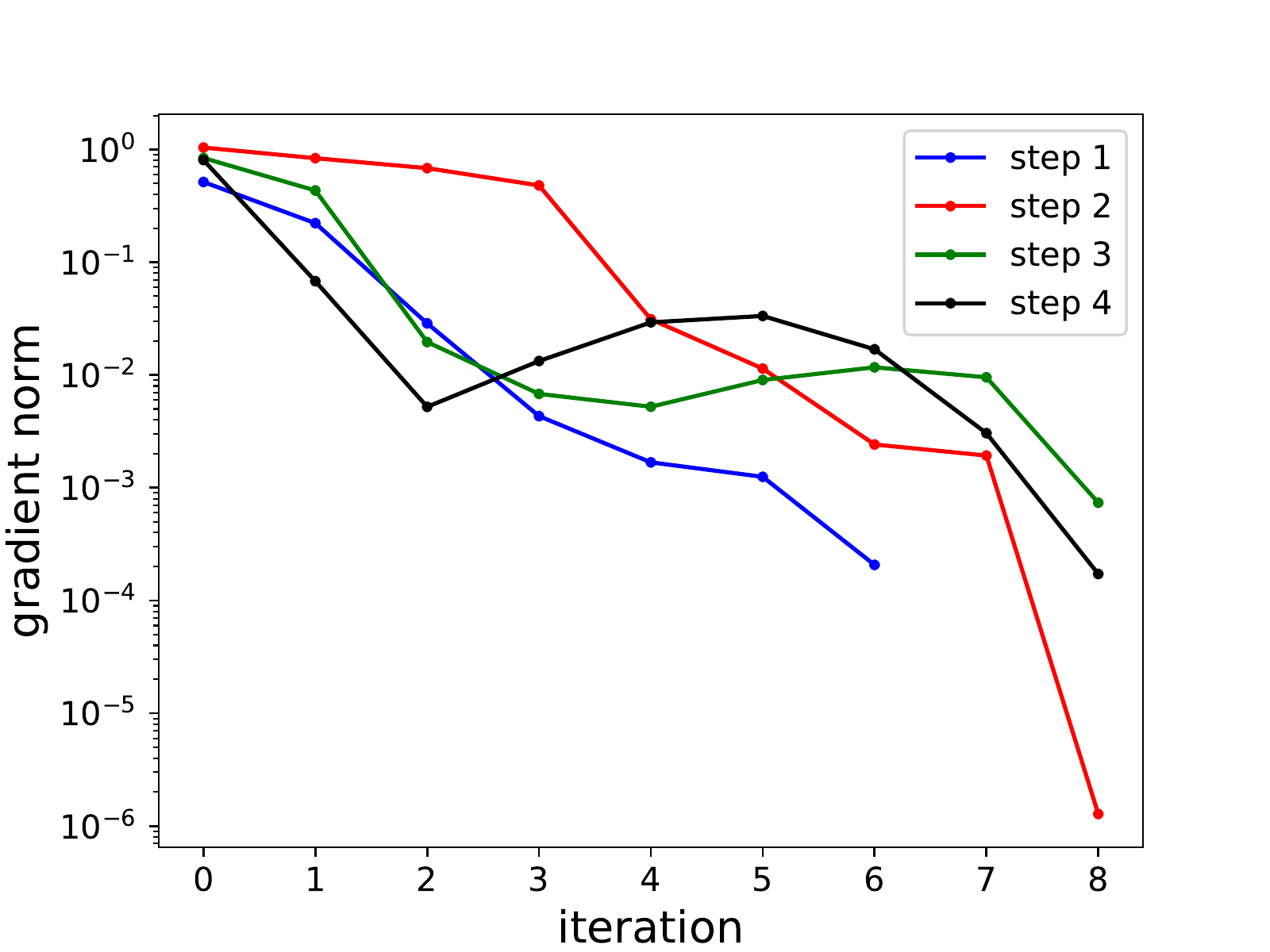}
		
		\includegraphics[width=0.45\textwidth]{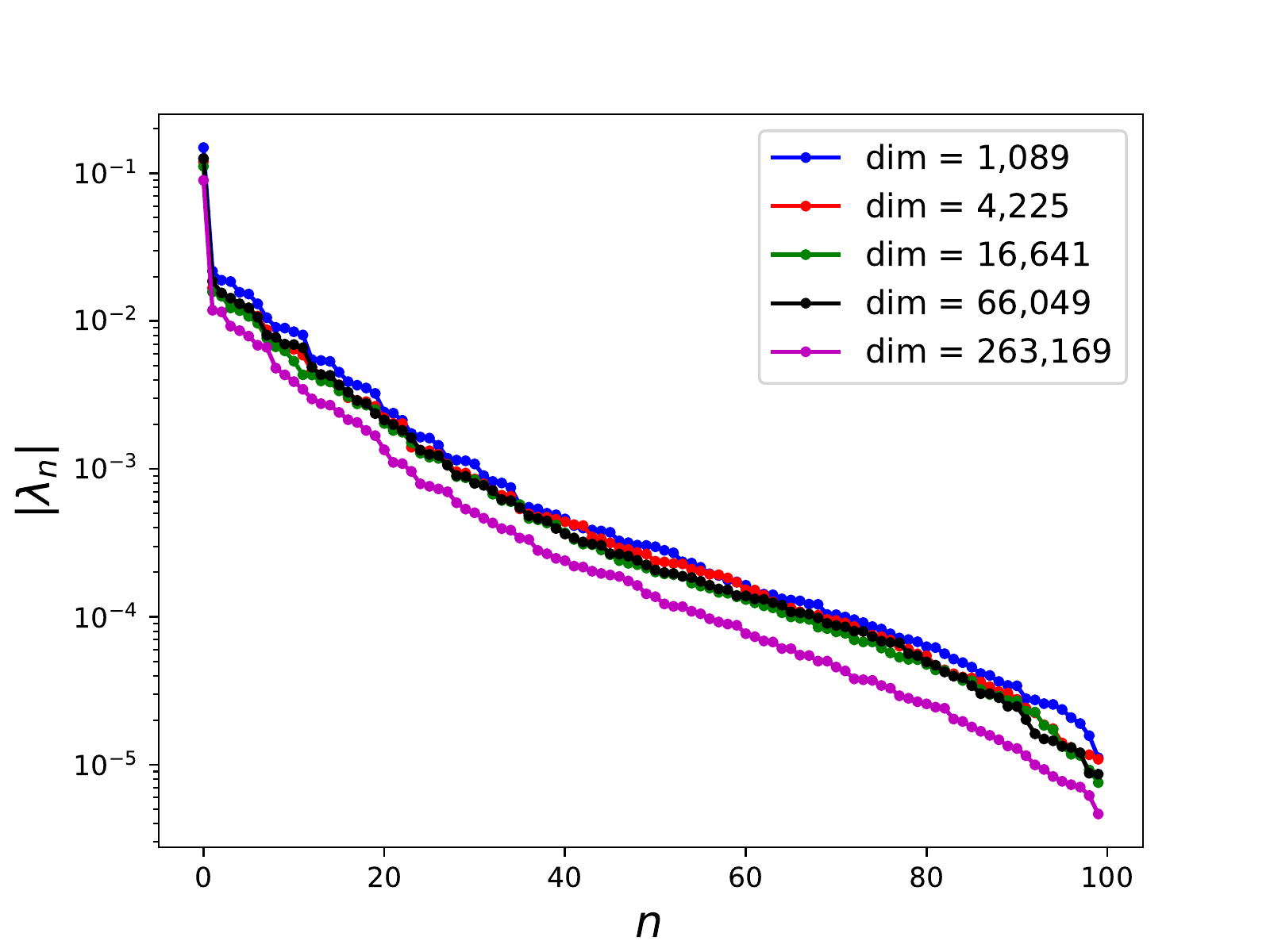}
		\includegraphics[width=0.45\textwidth]{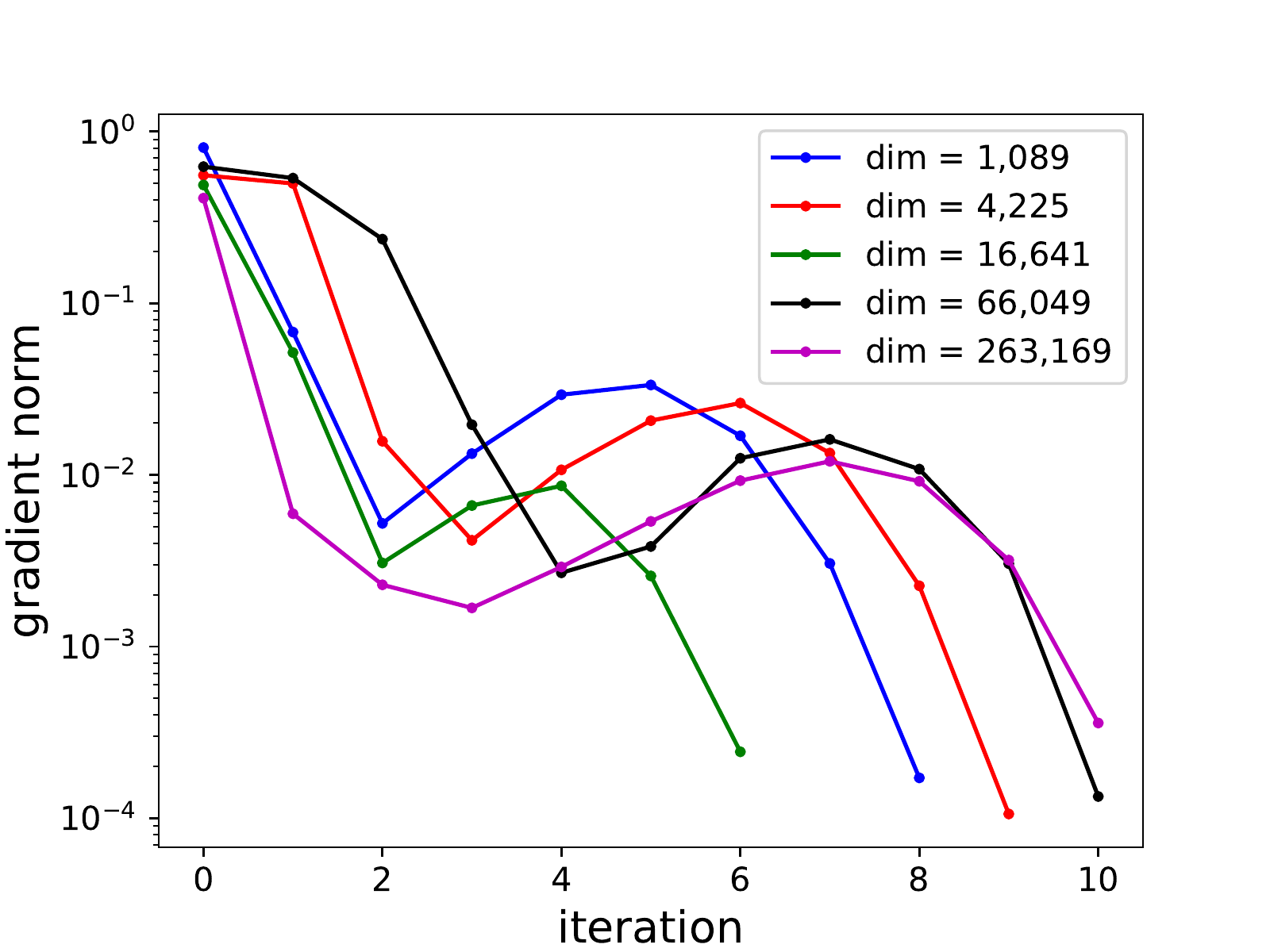}
		\caption{Top: decay of the eigenvalues \eqref{eq:GenEigen_f} for the quadratic approximation $T_2f$ (left) and decay of the gradient norm of the BFGS optimization (right)  at different steps (corresponding to Figure \ref{fig:chance}) with dimension $1,089$ for discrete $m$.
			Bottom: the same plot at optimal variables with different dimensions.}\label{fig:scalability}
	\end{figure}
	
	%\begin{figure}[!htb]
	%	\centering
	%	\includegraphics[width=0.48\textwidth]{figure/scalability-dimension}
	%	\includegraphics[width=0.48\textwidth]{figure/BFGS}
	%	\caption{Decay of the eigenvalues \eqref{eq:GenEigen_f} for the quadratic approximation $T_2 f$ (left) and decay of the the gradient norm of the BFGS optimization (right) at step $4$ with different parameter dimensions.}\label{fig:scalability}
	%\end{figure}
	
	Finally, we plot the decay of the eigenvalues of \eqref{eq:GenEigen_f} for the quadratic approximation $T_2 f$ at different optimization steps and different dimensions of the discrete random parameters in Figure \ref{fig:scalability}. By the similarity of the eigenvalue decay, we can conclude that the quadratic approximation is scalable with respect to the parameter dimensions in that the number of PDE solves are similar with increasing dimension. Moreover, the number of optimization iterations stays similar with increasing smoothing and penalty parameters $(\beta, \gamma)$ as well as increasing dimension, thus demonstrating that the optimization method is also scalable for this example.

	\section{Conclusion}
	\label{sec:conclusion}
	
	We proposed a Taylor approximation based continuation optimization method to solve chance and random PDE constrained optimization problems. We presented the derivation and efficient computation of the Taylor approximations using randomized algorithms. To address the challenges of discontinuous indicator function and inequality constraint, we employed a smooth approximation and a penalty method in a continuation BFGS optimization algorithm. We compared the accuracy of the Taylor (constant, linear, and quadratic) approximations and sample average approximation for both the chance evaluation and the optimal variable, demonstrated the acceleration by the Taylor approximation, reported the convergence of the continuation optimization algorithm in satisfying the chance constraint and minimizing the objective functional, and showed the scalability of the proposed method in that the number of PDE solves is essentially insensitive to increasing dimension of the random parameters. In particular, for the test problem the quadratic Taylor approximation achieves two orders of magnitude higher accuracy than SAA at 37X cheaper cost measured in the number of PDE solves.
	
	The following research directions are of great interest: (1) higher order (beyond quadratic) Taylor approximations may improve the accuracy and efficiency of the proposed method; these rely on efficient low rank tensor decomposition which only requires tensor action \cite{AlgerChenGhattas20}; (2) vector or function valued chance constraint functions \cite{Farshbaf-ShakerHenrionHoemberg18}, e.g., representing pointwise pressure or the whole pressure field, require further development of the optimization method with respect to the Taylor approximations and Lagrangian approach to computing the gradient; (3) to deal with extreme chance with a critical value $\alpha \ll 1$,  importance sampling \cite{PeherstorferKramerWillcox17, PeherstorferKramerWillcox18} with Taylor approximation in the failure region rather than at the mean of the random parameters can be employed; and (4) theoretical analysis of the convergence of the Taylor approximations and continuation optimization remain open.

	\bibliographystyle{plain}
	\bibliography{references}
	
\end{document}

%% file: shared.tex
\usepackage{geometry}
\usepackage[utf8]{inputenc}
\usepackage{amsmath, amssymb, amsfonts}
\usepackage{algorithm, algorithmic}
\usepackage{graphicx}
\usepackage{epstopdf}
\usepackage{hyperref}
\hypersetup{
	colorlinks=true,
	linkcolor=blue,
	filecolor=magenta,      
	urlcolor=cyan,
}

\usepackage{lipsum}
\allowdisplaybreaks
\ifpdf
  \DeclareGraphicsExtensions{.eps,.pdf,.png,.jpg}
\else
  \DeclareGraphicsExtensions{.eps}
\fi

\usepackage{amsopn}

\usepackage{xcolor}
\newcommand{\pc}[1]{{\color{black}{#1}}}

\newcommand{\bE}{{\mathbb{E}}}

\newcommand{\bI}{{\mathbb{I}}}

\newcommand{\bR}{{\mathbb{R}}}

\newcommand{\cC}{\mathcal{C}}

\newcommand{\cE}{\mathcal{E}}

\newcommand{\cH}{\mathcal{H}}
\newcommand{\cJ}{\mathcal{J}}

\newcommand{\cL}{\mathcal{L}}
\newcommand{\cM}{\mathcal{M}}
\newcommand{\cN}{\mathcal{N}}
\newcommand{\cP}{\mathcal{P}}
\newcommand{\cR}{\mathcal{R}}

\newcommand{\cS}{\mathcal{S}}
\newcommand{\cU}{\mathcal{U}}
\newcommand{\cV}{\mathcal{V}}

\newcommand{\cX}{\mathcal{X}}
\newcommand{\cY}{\mathcal{Y}}

\newcommand{\cZ}{\mathcal{Z}}

\newcommand{\beq}{\begin{equation}}
\newcommand{\eeq}{\end{equation}}
\newcommand{\ba}{\begin{array}}
\newcommand{\ea}{\end{array}}

%% file: article.bbl
\begin{thebibliography}{10}

\bibitem{AlexanderianPetraStadlerEtAl14}
Alen Alexanderian, Noemi Petra, Georg Stadler, and Omar Ghattas.
\newblock {A}-optimal design of experiments for infinite-dimensional {B}ayesian
  linear inverse problems with regularized $\ell_0$-sparsification.
\newblock {\em SIAM Journal on Scientific Computing}, 36(5):A2122--A2148, 2014.

\bibitem{AlexanderianPetraStadlerEtAl16}
Alen Alexanderian, Noemi Petra, Georg Stadler, and Omar Ghattas.
\newblock A fast and scalable method for {A}-optimal design of experiments for
  infinite-dimensional {B}ayesian nonlinear inverse problems.
\newblock {\em SIAM Journal on Scientific Computing}, 38(1):A243--A272, 2016.

\bibitem{AlexanderianPetraStadlerEtAl17}
Alen Alexanderian, Noemi Petra, Georg Stadler, and Omar Ghattas.
\newblock Mean-variance risk-averse optimal control of systems governed by
  {PDEs} with random parameter fields using quadratic approximations.
\newblock {\em SIAM/ASA Journal on Uncertainty Quantification},
  5(1):1166--1192, 2017.
\newblock arXiv preprint arXiv:1602.07592.

\bibitem{AlgerChenGhattas20}
Nick Alger, Peng Chen, and Omar Ghattas.
\newblock Tensor train construction from tensor actions, with application to
  compression of large high order derivative tensors.
\newblock {\em arXiv preprint arXiv:2002.06244, to appear in SIAM Journal on
  Scientific Computing}, 2020.

\bibitem{AliUllmannHinze17}
Ahmad~Ahmad Ali, Elisabeth Ullmann, and Michael Hinze.
\newblock Multilevel {M}onte {C}arlo analysis for optimal control of elliptic
  {PDEs} with random coefficients.
\newblock {\em SIAM/ASA Journal on Uncertainty Quantification}, 5(1):466--492,
  2017.

\bibitem{AllaHinzeKolvenbachEtAl19}
Alessandro Alla, Michael Hinze, Philip Kolvenbach, Oliver Lass, and Stefan
  Ulbrich.
\newblock A certified model reduction approach for robust parameter
  optimization with pde constraints.
\newblock {\em Advances in Computational Mathematics}, 45(3):1221--1250, 2019.

\bibitem{BashirWillcoxGhattasEtAl08}
O.~Bashir, K.~Willcox, O.~Ghattas, B.~van Bloemen~Waanders, and J.~Hill.
\newblock Hessian-based model reduction for large-scale systems with initial
  condition inputs.
\newblock {\em International Journal for Numerical Methods in Engineering},
  73:844--868, 2008.

\bibitem{BennerOnwuntaStoll16}
P.~Benner, A.~Onwunta, and M.~Stoll.
\newblock Block-diagonal preconditioning for optimal control problems
  constrained by {PDEs} with uncertain inputs.
\newblock {\em SIAM Journal on Matrix Analysis and Applications},
  37(2):491--518, 2016.

\bibitem{Borzi10}
A.~Borz{\`\i}.
\newblock Multigrid and sparse-grid schemes for elliptic control problems with
  random coefficients.
\newblock {\em Computing and Visualization in Science}, 13(4):153--160, 2010.

\bibitem{BorziSchulzSchillingsEtAl10}
A~Borz{\`\i}, V~Schulz, C~Schillings, and G~Von~Winckel.
\newblock On the treatment of distributed uncertainties in {PDE}-constrained
  optimization.
\newblock {\em GAMM-Mitteilungen}, 33(2):230--246, 2010.

\bibitem{Bui-ThanhBursteddeGhattasEtAl12}
Tan Bui-Thanh, Carsten Burstedde, Omar Ghattas, James Martin, Georg Stadler,
  and Lucas~C. Wilcox.
\newblock {Extreme-scale UQ for Bayesian inverse problems governed by PDEs}.
\newblock In {\em SC12: Proceedings of the International Conference for High
  Performance Computing, Networking, Storage and Analysis}, 2012.

\bibitem{Bui-ThanhGhattas12a}
Tan Bui-Thanh and Omar Ghattas.
\newblock Analysis of the {H}essian for inverse scattering problems. {P}art
  {I}: Inverse shape scattering of acoustic waves.
\newblock {\em Inverse Problems}, 28(5):055001, 2012.

\bibitem{Bui-ThanhGhattas12}
Tan Bui-Thanh and Omar Ghattas.
\newblock Analysis of the {H}essian for inverse scattering problems. {P}art
  {II}: Inverse medium scattering of acoustic waves.
\newblock {\em Inverse Problems}, 28(5):055002, 2012.

\bibitem{Bui-ThanhGhattas13a}
Tan Bui-Thanh and Omar Ghattas.
\newblock Analysis of the {H}essian for inverse scattering problems. {P}art
  {III}: Inverse medium scattering of electromagnetic waves.
\newblock {\em Inverse Problems and Imaging}, 7(4):1139--1155, 2013.

\bibitem{Bui-ThanhGhattas15}
Tan Bui-Thanh and Omar Ghattas.
\newblock A scalable {MAP} solver for {B}ayesian inverse problems with {B}esov
  priors.
\newblock {\em Inverse Problems and Imaging}, 9(1):27--54, 2015.

\bibitem{Bui-ThanhGhattasMartinEtAl13}
Tan Bui-Thanh, Omar Ghattas, James Martin, and Georg Stadler.
\newblock A computational framework for infinite-dimensional {B}ayesian inverse
  problems {P}art {I}: {T}he linearized case, with application to global
  seismic inversion.
\newblock {\em SIAM Journal on Scientific Computing}, 35(6):A2494--A2523, 2013.

\bibitem{ChenMangasarian95}
Chunhui Chen and Olvi~L Mangasarian.
\newblock Smoothing methods for convex inequalities and linear complementarity
  problems.
\newblock {\em Mathematical programming}, 71(1):51--69, 1995.

\bibitem{ChenGhattas19}
P.~Chen and O.~Ghattas.
\newblock Sparse polynomial approximation for optimal control problems
  constrained by elliptic {PDEs} with lognormal random coefficients.
\newblock {\em submitted}, 2019.
\newblock https://arxiv.org/abs/1903.05547.

\bibitem{ChenQuarteroni14}
P.~Chen and A.~Quarteroni.
\newblock Weighted reduced basis method for stochastic optimal control problems
  with elliptic {PDE} constraints.
\newblock {\em SIAM/ASA J. Uncertainty Quantification}, 2(1):364--396, 2014.

\bibitem{ChenQuarteroniRozza16}
P.~Chen, A.~Quarteroni, and G.~Rozza.
\newblock Multilevel and weighted reduced basis method for stochastic optimal
  control problems constrained by {S}tokes equations.
\newblock {\em Numerische Mathematik}, 133(1):67--102, 2016.

\bibitem{ChenVillaGhattas17}
P.~Chen, U.~Villa, and O.~Ghattas.
\newblock Hessian-based adaptive sparse quadrature for infinite-dimensional
  {B}ayesian inverse problems.
\newblock {\em Computer Methods in Applied Mechanics and Engineering},
  327:147--172, 2017.

\bibitem{Chen18}
Peng Chen.
\newblock Sparse quadrature for high-dimensional integration with {G}aussian
  measure.
\newblock {\em ESAIM: Mathematical Modelling and Numerical Analysis},
  52(2):631--657, 2018.

\bibitem{ChenGhattas19a}
Peng Chen and Omar Ghattas.
\newblock Hessian-based sampling for high-dimensional model reduction.
\newblock {\em International Journal for Uncertainty Quantification}, 9(2),
  2019.

\bibitem{ChenGhattas20}
Peng Chen and Omar Ghattas.
\newblock Projected stein variational gradient descent.
\newblock In {\em Advances in Neural Information Processing Systems}, 2020.

\bibitem{ChenHabermanGhattas20}
Peng Chen, Michael Haberman, and Omar Ghattas.
\newblock Optimal design of acoustic cloak under uncertainty.
\newblock {\em arXiv:2007.13252}, 2020.

\bibitem{ChenQuarteroniRozza13}
Peng Chen, Alfio Quarteroni, and Gianluigi Rozza.
\newblock Stochastic optimal {R}obin boundary control problems of
  advection-dominated elliptic equations.
\newblock {\em SIAM Journal on Numerical Analysis}, 51(5):2700--2722, 2013.

\bibitem{ChenVillaGhattas19}
Peng Chen, Umberto Villa, and Omar Ghattas.
\newblock Taylor approximation and variance reduction for {PDE}-constrained
  optimal control under uncertainty.
\newblock {\em Journal of Computational Physics}, 385:163--186, 2019.

\bibitem{ChenWuChenEtAl19a}
Peng Chen, Keyi Wu, Joshua Chen, Thomas O'Leary-Roseberry, and Omar Ghattas.
\newblock Projected {S}tein variational {N}ewton: {A} fast and scalable
  {B}ayesian inference method in high dimensions.
\newblock {\em Advances in Neural Information Processing Systems}, 2019.

\bibitem{ChenWuGhattas20}
Peng Chen, Keyi Wu, and Omar Ghattas.
\newblock Bayesian inference of heterogeneous epidemic models: {Application to
  COVID-19} spread accounting for long-term care facilities.
\newblock {\em arXiv preprint arXiv:2011.01058}, 2020.

\bibitem{CrestelAlexanderianStadlerEtAl17}
Benjamin Crestel, Alen Alexanderian, Georg Stadler, and Omar Ghattas.
\newblock {A}-optimal encoding weights for nonlinear inverse problems, with
  application to the {H}elmholtz inverse problem.
\newblock {\em Inverse Problems}, 33(7):074008, 2017.

\bibitem{DeHamptonMauteEtAl19}
Subhayan De, Jerrad Hampton, Kurt Maute, and Alireza Doostan.
\newblock Topology optimization under uncertainty using a stochastic
  gradient-based approach.
\newblock {\em arXiv preprint arXiv:1902.04562}, 2019.

\bibitem{DickLeGiaSchwab16}
J.~Dick, Q.T. Le~Gia, and Ch. Schwab.
\newblock Higher order quasi--{Monte Carlo} integration for holomorphic,
  parametric operator equations.
\newblock {\em SIAM/ASA Journal on Uncertainty Quantification}, 4(1):48--79,
  2016.

\bibitem{Farshbaf-ShakerHenrionHoemberg18}
M~Hassan Farshbaf-Shaker, Ren{\'e} Henrion, and Dietmar H{\"o}mberg.
\newblock Properties of chance constraints in infinite dimensions with an
  application to {PDE} constrained optimization.
\newblock {\em Set-Valued and Variational Analysis}, 26(4):821--841, 2018.

\bibitem{FlathWilcoxAkcelikEtAl11}
Pearl~H. Flath, Lucas~C. Wilcox, Volkan Ak\c{c}elik, Judy Hill, Bart van
  Bloemen~Waanders, and Omar Ghattas.
\newblock Fast algorithms for {B}ayesian uncertainty quantification in
  large-scale linear inverse problems based on low-rank partial {H}essian
  approximations.
\newblock {\em SIAM Journal on Scientific Computing}, 33(1):407--432, 2011.

\bibitem{GarreisSurowiecUlbrichEtAl19}
S.~Garreis, T.M. Surowiec, and M.~Ulbrich.
\newblock An interior-point approach for solving risk-averse {PDE}-constrained
  optimization problems with coherent risk measures.
\newblock {\em Preprint, submitted, Technical University of Munich}, 2019.

\bibitem{GeiersbachLoayza-RomeroWelker20}
Caroline Geiersbach, Estefania Loayza-Romero, and Kathrin Welker.
\newblock Stochastic approximation for optimization in shape spaces.
\newblock {\em arXiv preprint arXiv:2001.10786}, 2020.

\bibitem{GeiersbachScarinci20}
Caroline Geiersbach and Teresa Scarinci.
\newblock Stochastic proximal gradient methods for nonconvex problems in
  {H}ilbert spaces.
\newblock {\em arXiv preprint arXiv:2001.01329}, 2020.

\bibitem{GeiersbachWollner19}
Caroline Geiersbach and Winnifried Wollner.
\newblock A stochastic gradient method with mesh refinement for {PDE}
  constrained optimization under uncertainty.
\newblock {\em arXiv preprint arXiv:1905.08650}, 2019.

\bibitem{GeletuHoffmannSchmidtEtAl20}
Abebe Geletu, Armin Hoffmann, Patrick Schmidt, and Pu~Li.
\newblock Chance constrained optimization of elliptic {PDE} systems with a
  smoothing convex approximation.
\newblock {\em ESAIM: Control, Optimisation and Calculus of Variations}, 26:70,
  2020.

\bibitem{Gunzburger03}
Max~D. Gunzburger.
\newblock {\em Perspectives in Flow Control and Optimization}.
\newblock SIAM, Philadelphia, 2003.

\bibitem{GunzburgerLeeLee11}
Max~D. Gunzburger, Hyung-Chun Lee, and Jangwoon Lee.
\newblock Error estimates of stochastic optimal {N}eumann boundary control
  problems.
\newblock {\em SIAM Journal on Numerical Analysis}, 49(4):1532--1552, 2011.

\bibitem{GuoXuZhang17}
Shaoyan Guo, Huifu Xu, and Liwei Zhang.
\newblock Convergence analysis for mathematical programs with distributionally
  robust chance constraint.
\newblock {\em SIAM Journal on optimization}, 27(2):784--816, 2017.

\bibitem{HalkoMartinssonTropp11}
Nathan Halko, Per~Gunnar Martinsson, and Joel~A. Tropp.
\newblock Finding structure with randomness: {P}robabilistic algorithms for
  constructing approximate matrix decompositions.
\newblock {\em SIAM Review}, 53(2):217--288, 2011.

\bibitem{HinzePinnauUlbrichEtAl08}
M.~Hinze, R.~Pinnau, M.~Ulbrich, and S.~Ulbrich.
\newblock {\em Optimization with PDE constraints}, volume~23.
\newblock Springer Science \& Business Media, 2008.

\bibitem{HouLeeManouzi11}
L.~S. Hou, J.~Lee, and H.~Manouzi.
\newblock Finite element approximations of stochastic optimal control problems
  constrained by stochastic elliptic {PDE}s.
\newblock {\em Journal of Mathematical Analysis and Applications},
  384(1):87--103, 2011.

\bibitem{IsaacPetraStadlerEtAl15}
Tobin Isaac, Noemi Petra, Georg Stadler, and Omar Ghattas.
\newblock Scalable and efficient algorithms for the propagation of uncertainty
  from data through inference to prediction for large-scale problems, with
  application to flow of the {A}ntarctic ice sheet.
\newblock {\em Journal of Computational Physics}, 296:348--368, September 2015.

\bibitem{KolvenbachLassUlbrich18}
Philip Kolvenbach, Oliver Lass, and Stefan Ulbrich.
\newblock {An approach for robust PDE-constrained optimization with application
  to shape optimization of electrical engines and of dynamic elastic structures
  under uncertainty}.
\newblock {\em Optimization and Engineering}, 19(3):697--731, 2018.

\bibitem{KouriSurowiec16}
D.~P. Kouri and T.~M. Surowiec.
\newblock Risk-averse {PDE}-constrained optimization using the conditional
  value-at-risk.
\newblock {\em SIAM Journal on Optimization}, 26(1):365--396, 2016.

\bibitem{KouriHeinkenschloosVanBloemenWaanders12}
D.P. Kouri, D.~Heinkenschloos, M.~Ridzal, and B.G. Van Bloemen~Waanders.
\newblock A trust-region algorithm with adaptive stochastic collocation for
  {PDE} optimization under uncertainty.
\newblock {\em SIAM Journal on Scientific Computing}, 35(4):1847--1879, 2012.

\bibitem{KouriSurowiec18}
Drew~Philip Kouri and Thomas~M Surowiec.
\newblock Existence and optimality conditions for risk-averse {PDE}-constrained
  optimization.
\newblock {\em SIAM/ASA Journal on Uncertainty Quantification}, 6(2):787--815,
  2018.

\bibitem{KunothSchwab16}
A.~Kunoth and Ch. Schwab.
\newblock Sparse adaptive tensor {G}alerkin approximations of stochastic
  {PDE}-constrained control problems.
\newblock {\em SIAM/ASA Journal on Uncertainty Quantification},
  4(1):1034--1059, 2016.

\bibitem{KunothSchwab13}
Angela Kunoth and Christoph Schwab.
\newblock Analytic regularity and {GPC} approximation for control problems
  constrained by linear parametric elliptic and parabolic {PDE}s.
\newblock {\em SIAM Journal on Control and Optimization}, 51(3):2442--2471,
  2013.

\bibitem{LassUlbrich17}
Oliver Lass and Stefan Ulbrich.
\newblock Model order reduction techniques with a posteriori error control for
  nonlinear robust optimization governed by partial differential equations.
\newblock {\em SIAM Journal on Scientific Computing}, 39(5):S112--S139, 2017.

\bibitem{LiStadler19}
C.~Li and G.~Stadler.
\newblock Sparse solutions in optimal control of {PDEs} with uncertain
  parameters: {T}he linear case.
\newblock {\em SIAM Journal on Control and Optimization}, 57(1):633--658, 2019.

\bibitem{LiWangZhang18}
Jingshi Li, Xiaoshen Wang, and Kai Zhang.
\newblock An efficient alternating direction method of multipliers for optimal
  control problems constrained by random helmholtz equations.
\newblock {\em Numerical Algorithms}, 78(1):161--191, 2018.

\bibitem{LindgrenRueLindstroem11}
Finn Lindgren, H{\aa}vard Rue, and Johan Lindstr{\"o}m.
\newblock An explicit link between {G}aussian fields and {G}aussian {M}arkov
  random fields: the stochastic partial differential equation approach.
\newblock {\em Journal of the Royal Statistical Society: Series B (Statistical
  Methodology)}, 73(4):423--498, 2011.

\bibitem{Lions71}
Jacques~Louis Lions.
\newblock {\em Optimal Control of Systems Governed by Partial Differential
  Equations}, volume 170 of {\em Grundlehren der mathematischen
  Wissenschaften}.
\newblock Springer-Verlag Berlin Heidelberg, 1971.

\bibitem{LoggMardalGarth12}
Anders Logg, Kent-Andre Mardal, and Garth Wells.
\newblock {\em Automated Solution of Differential Equations by the Finite
  Element Method: {T}he {FEniCS} book}, volume~84.
\newblock Springer Science \& Business Media, 2012.

\bibitem{MaLiJiang18}
Lingling Ma, Qiuqi Li, and Lijian Jiang.
\newblock Local--global model reduction method for stochastic optimal control
  problems constrained by partial differential equations.
\newblock {\em Computer Methods in Applied Mechanics and Engineering},
  339:514--541, 2018.

\bibitem{MartinWilcoxBursteddeEtAl12}
James Martin, Lucas~C. Wilcox, Carsten Burstedde, and Omar Ghattas.
\newblock A stochastic {Newton MCMC} method for large-scale statistical inverse
  problems with application to seismic inversion.
\newblock {\em SIAM Journal on Scientific Computing}, 34(3):A1460--A1487, 2012.

\bibitem{MartinNobileTsilifis19}
Matthieu Martin, Fabio Nobile, and Panagiotis Tsilifis.
\newblock A multilevel stochastic gradient method for pde-constrained optimal
  control problems with uncertain parameters.
\newblock {\em arXiv preprint arXiv:1912.11900}, 2019.

\bibitem{MoralesNocedal11}
Jos{\'e}~Luis Morales and Jorge Nocedal.
\newblock {Remark on ``Algorithm 778: L-BFGS-B: Fortran subroutines for
  large-scale bound constrained optimization"}.
\newblock {\em ACM Transactions on Mathematical Software (TOMS)}, 38(1):1--4,
  2011.

\bibitem{NgWillcox14}
L.~Ng and K.~Willcox.
\newblock Multifidelity approaches for optimization under uncertainty.
\newblock {\em International Journal for Numerical Methods in Engineering},
  100(10):746--772, 2014.

\bibitem{NocedalWright06}
Jorge Nocedal and Stephen~J. Wright.
\newblock {\em Numerical Optimization}.
\newblock Springer Verlag, Berlin, Heidelberg, New York, second edition, 2006.

\bibitem{PeherstorferKramerWillcox17}
Benjamin Peherstorfer, Boris Kramer, and Karen Willcox.
\newblock Combining multiple surrogate models to accelerate failure probability
  estimation with expensive high-fidelity models.
\newblock {\em Journal of Computational Physics}, 341:61--75, 2017.

\bibitem{PeherstorferKramerWillcox18}
Benjamin Peherstorfer, Boris Kramer, and Karen Willcox.
\newblock Multifidelity preconditioning of the cross-entropy method for rare
  event simulation and failure probability estimation.
\newblock {\em SIAM/ASA Journal on Uncertainty Quantification}, 6(2):737--761,
  2018.

\bibitem{PetraMartinStadlerEtAl14}
Noemi Petra, James Martin, Georg Stadler, and Omar Ghattas.
\newblock A computational framework for infinite-dimensional {B}ayesian inverse
  problems: {P}art {II}. {S}tochastic {N}ewton {MCMC} with application to ice
  sheet flow inverse problems.
\newblock {\em SIAM Journal on Scientific Computing}, 36(4):A1525--A1555, 2014.

\bibitem{QiSunZhou00}
Liqun Qi, Defeng Sun, and Guanglu Zhou.
\newblock A new look at smoothing newton methods for nonlinear complementarity
  problems and box constrained variational inequalities.
\newblock {\em Mathematical programming}, 87(1):1--35, 2000.

\bibitem{RoaldAndersson17}
Line Roald and G{\"o}ran Andersson.
\newblock Chance-constrained {AC} optimal power flow: {R}eformulations and
  efficient algorithms.
\newblock {\em IEEE Transactions on Power Systems}, 33(3):2906--2918, 2017.

\bibitem{RosseelWells12}
Eveline Rosseel and Garth~N Wells.
\newblock Optimal control with stochastic {PDE} constraints and uncertain
  controls.
\newblock {\em Computer Methods in Applied Mechanics and Engineering},
  213:152--167, 2012.

\bibitem{RuszczynskiShapiro06}
Andrzej Ruszczy{\'n}ski and Alexander Shapiro.
\newblock Optimization of risk measures.
\newblock In {\em Probabilistic and Randomized Methods for Design Under
  Uncertainty}, pages 119--157. Springer, 2006.

\bibitem{SaibabaLeeKitanidis16}
Arvind~K Saibaba, Jonghyun Lee, and Peter~K Kitanidis.
\newblock Randomized algorithms for generalized {H}ermitian eigenvalue problems
  with application to computing {K}arhunen--{L}o{\`e}ve expansion.
\newblock {\em Numerical Linear Algebra with Applications}, 23(2):314--339,
  2016.

\bibitem{SchillingsSchwab13}
Claudia Schillings and Christoph Schwab.
\newblock Sparse, adaptive {S}molyak quadratures for {B}ayesian inverse
  problems.
\newblock {\em Inverse Problems}, 29(6):065011, 2013.

\bibitem{ShapiroDentchevaRuszczynski09}
Alexander Shapiro, Darinka Dentcheva, and Andrezj Ruszczynski.
\newblock {\em Lectures on Stochastic Programming: Modeling and Theory}.
\newblock Society for Industrial and Applied Mathematics, 2009.

\bibitem{Stuart10}
Andrew~M. Stuart.
\newblock Inverse problems: {A B}ayesian perspective.
\newblock {\em Acta Numerica}, 19:451--559, 2010.

\bibitem{TieslerKirbyXiuEtAl12}
Hanne Tiesler, Robert~M Kirby, Dongbin Xiu, and Tobias Preusser.
\newblock Stochastic collocation for optimal control problems with stochastic
  {PDE} constraints.
\newblock {\em SIAM Journal on Control and Optimization}, 50(5):2659--2682,
  2012.

\bibitem{Troeltzsch10}
F.~Tr{\"o}ltzsch.
\newblock {\em Optimal Control of Partial Differential Equations: Theory,
  Methods, and Applications}, volume 112.
\newblock American Mathematical Society, Providence, RI, 2010.

\bibitem{Uryasev13}
Stanislav Uryasev.
\newblock {\em Probabilistic constrained optimization: methodology and
  applications}, volume~49.
\newblock Springer Science \& Business Media, 2013.

\bibitem{AckooijHenrionPerez-Aros19}
Wim van Ackooij, Ren{\'e} Henrion, and Pedro P{\'e}rez-Aros.
\newblock Generalized gradients for probabilistic/robust (probust) constraints.
\newblock {\em Optimization}, pages 1--29, 2019.

\bibitem{VanAckooijMalick19}
Wim Van~Ackooij and J{\'e}r{\^o}me Malick.
\newblock Eventual convexity of probability constraints with elliptical
  distributions.
\newblock {\em Mathematical Programming}, 175(1-2):1--27, 2019.

\bibitem{AckooijPerez-Aros19}
Wim van Ackooij and Pedro P\`erez-Aros.
\newblock Generalized differentiation of probability functions acting on an
  infinite system of constraints.
\newblock {\em SIAM Journal on Optimization}, 29(3):2179--2210, 2019.

\bibitem{VanBarelVandewalle19}
A.~Van~Barel and S.~Vandewalle.
\newblock Robust optimization of {PDEs} with random coefficients using a
  multilevel {Monte Carlo} method.
\newblock {\em SIAM/ASA Journal on Uncertainty Quantification}, 7(1):174--202,
  2019.

\bibitem{WuChenGhattas20}
Keyi Wu, Peng Chen, and Omar Ghattas.
\newblock A fast and scalable computational framework for large-scale and
  high-dimensional {B}ayesian optimal experimental design.
\newblock {\em arXiv preprint arXiv:2010.15196}, 2020.

\bibitem{XuBoyceZhangEtAl17}
Bin Xu, Scott~E Boyce, Yu~Zhang, Qiang Liu, Le~Guo, and Ping-An Zhong.
\newblock Stochastic programming with a joint chance constraint model for
  reservoir refill operation considering flood risk.
\newblock {\em Journal of Water Resources Planning and Management},
  143(1):04016067, 2017.

\bibitem{YangGunzburger17}
Huanhuan Yang and Max Gunzburger.
\newblock Algorithms and analyses for stochastic optimization for turbofan
  noise reduction using parallel reduced-order modeling.
\newblock {\em Computer Methods in Applied Mechanics and Engineering},
  319:217--239, 2017.

\bibitem{ZahrCarlbergKouri19}
Matthew~J Zahr, Kevin~T Carlberg, and Drew~P Kouri.
\newblock An efficient, globally convergent method for optimization under
  uncertainty using adaptive model reduction and sparse grids.
\newblock {\em SIAM/ASA Journal on Uncertainty Quantification}, 7(3):877--912,
  2019.

\end{thebibliography}
